%% file: nts-main.tex
\documentclass{article}
\usepackage[pagebackref=true]{hyperref}%A: to have links to equations and references
\hypersetup{colorlinks=true,linkcolor=blue,citecolor=red,pagebackref=true}
\usepackage[a4paper,margin=1in]{geometry} %A to change easily format
\usepackage{amssymb}
\usepackage{amsmath}
\usepackage{mathabx} %A
\usepackage{graphicx}
\newtheorem{theorem}{Theorem}
\newtheorem{definition}{Definition}
\newtheorem{lemma}[theorem]{Lemma}
\newtheorem{rem}[theorem]{Remark}

\DeclareMathOperator{\inter}{int}
\DeclareMathOperator{\sgn}{sgn}
\DeclareMathOperator{\im}{Im}%A
\DeclareMathOperator{\re}{Re}%A
\newcommand{\abs}[1]{\lvert#1\rvert} %A
\newcommand{\norm}[1]{\lVert#1\rVert} %A
\newcommand{\cover}[1]{\stackrel{#1}{\Longrightarrow}}
\newcommand{\T}{\mathcal{T}}
\newcommand{\J}{\mathcal{J}}
\newcommand{\out}{\mathrm{out}}
\newcommand{\inn}{\mathrm{in}} %A

\newcommand{\outcomment}[1]{}

\def\qed{{\hfill{\vrule height5pt width3pt depth0pt}\medskip}}

\begin{document}
\begin{center} {\bf \LARGE  Shadowing of non-transversal heteroclinic chains in lattices}\\
 \vskip 0.5cm
%{\large Amadeu Delshams}\footnote{Partially supported by the Spanish MINECO/FEDER Grant PID2021-123968NB-I00}\\
%Dept. de Matem\`atiques, Univ. Polit\`ecnica de Catalunya (UPC), Barcelona, Spain\\
%Lab of Geometry and Dynamical Systems, UPC, Barcelona, Spain\\
%Centre de Recerca Matem\`atica, Barcelona, Spain\\
%Amadeu.Delshams@upc.edu
{\large Amadeu Delshams}\footnote{Partially supported by the Spanish MINECO/FEDER Grant PID2021-123968NB-I00}\\
Lab of Geometry and Dynamical Systems and IMTech, Universitat Politècnica de Catalunya (UPC) and Centre de Recerca Matem{\`a}tica (CRM),
Barcelona, Spain\\
\texttt{Amadeu.Delshams@upc.edu}
 \vskip\baselineskip
   {\large Piotr Zgliczy\'nski}\footnote{Research has been supported by polish NCN grant
  2016/22/A/ST1/00077 %Maestro
  }     \\
 Jagiellonian University, Institute of Computer Science and Computational Mathematics, \\
{\L}ojasiewicza 6, 30--348  Krak\'ow, Poland \\
\texttt{umzglicz@cyf-kr.edu.pl}

\vskip 0.5cm
 \today

\end{center}

\begin{abstract}
We deal with dynamical systems on complex lattices possessing chains of non-transversal heteroclinic connections between several periodic orbits.
The systems we consider are inspired by the so-called \emph{toy model systems} (TMS) used to prove the existence of energy transfer from low to high frequencies in the \emph{nonlinear cubic Schr\"odinger equation} (NLS) or generalizations. Using the geometric properties of the complex projective space as a base space, we generate in a natural way collections of such systems containing this type of chains, both in the Hamiltonian and in the non-Hamiltonian setting. On the other hand, we characterize the property of block diagonal dynamics along the heteroclinic connections that allows these chains to be shadowed, a property which in general only holds for transversal heteroclinic connections. Due to the lack of transversality, only finite chains are shadowed, since there is a dropping dimensions mechanism in the evolution of any disk close to them. The main shadowing technical  tool used in our work is the notion of covering relations as introduced by one of the authors.
\end{abstract}

\textbf{keywords:} shadowing, heteroclinic chain, non-transversal intersection, covering relations

\textbf{MSC:} 37C50, %Approximate trajectories (pseudotrajectories, shadowing, etc.) in smooth dynamics
  37J40 %

\tableofcontents

\input intro.tex

\input examples.tex

\input genToyModel.tex

\input covrel.tex

\input enclosure.tex

\input scheme.tex

\input ref.tex
\input appendix.tex

\end{document}

%% file: intro.tex
\section{Introduction and main results}

In this paper we deal with dynamical systems on some complex lattices, called \emph{toy model systems} (just TMS from now on), as they are inspired by those used in~\cite{CK,GK,GHP,Gu,Gi} to prove the existence of energy transfer from low to high frequencies in the \emph{nonlinear cubic Schrödinger equation} or generalizations (just NLS from now on).
From the point of view of dynamical systems, there is a remarkable feature in these papers, which is that the authors were able to shadow a chain of highly degenerate non-transversal heteroclinic connections between several periodic orbits. The length of the non-transversal heteroclinic chain is arbitrary, but finite, which means that we will only deal with \emph{finite} lattices of arbitrary length, indeed finite one-dimensional complex lattices.

Once a TMS is introduced for an NLS, the transfer of energy from low to high frequencies in the NLS translates to the existence of trajectories that pass near a chain of non-transversal heteroclinic connections between a finite sequence of periodic orbits located along different modes of the TMS. It is worth noticing that the standard \emph{shadowing lemmas} near heteroclinic connections between periodic orbits, or more generally invariant tori, are only valid for \emph{transversal} heteroclinic connections (this is the \emph{obstruction property} introduced by Arnold~\cite{Ar64}, see also~\cite[Chapter~11]{DLS06}), and therefore cannot be applied to any of the previously mentioned TMS.

In none of the aforementioned references we were able to find a geometrical description clear enough to allow us to understand how this shadowing property is achieved, so that it could be easily applied to other systems. In our previous work~\cite{DSZ} we presented a mechanism, which we believe gives a geometrical explanation of what is happening. In~\cite{DSZ} we emphasized the basic idea of \emph{dropping} dimensions, or equivalently directions, by assuming a geometric property, which we call \emph{block diagonal dynamics along the heteroclinic connections}, to avoid meeting again dropped dimensions. We applied this dropping dimensions mechanism to a very special simplification of TMS, although we already claimed that this simplified model contained all the essential difficulties of TMS.

In this paper we generalize TMS, simply by assuming sufficiently clear and understandable geometrical conditions, which allow placing naturally the TMS in the complex projective space and producing a non-transversal heteroclinic chain along periodic orbits. Thanks to the fulfillment of the block diagonal dynamics along the heteroclinic connections, we can prove that the periodic orbits are subsequently shadowed by the trajectories of the system connecting neighborhoods of the first and last periodic orbit. This geometric property is automatically satisfied in Hamiltonian systems, but it nevertheless needs to be imposed for non-Hamiltonian TMS, for which no previous theory of non-transversal shadowing existed, theory which could therefore be applied to non-Hamiltonian PDEs.

In this paper we strive to establish an abstract framework, which we hope will make it easier to apply this technique to other systems, both PDEs and ODEs, in questions related to the existence of diffusing orbits. The term \emph{diffusing orbit} relates to the Arnold's diffusion~\cite{Ar64} for nearly integrable Hamiltonian systems. Throughout the paper  an orbit shadowing a heteroclinic chain will be called a \emph{diffusing orbit}, and occasionally the existence of such an orbit will be referred to as a diffusion.

In our picture we think of evolving a disk of  dimension $k$  along a heteroclinic  chain and when a given
transition is not transversal, then we ``drop'' one or more dimensions of our disk, i.e., we select a lower dimensional subdisk
``parallel to expanding directions in future transitions''.  After at most $k$ transitions, our disk is a single point
and we cannot proceed further. We will refer to this phenomenon as the \emph{dropping dimensions mechanism}.
While thinking about disks has some geometrical appeal,
we will instead consider in our construction a thickened disk called h-set in the terminology of~\cite{ZGi} and our approach will be purely topological (just like the one presented in~\cite{CK}).

The main technical shadowing tool used in our work is the notion of \emph{covering relations} as introduced in~\cite{ZGi}, which differs from the notion used under the same name
in~\cite{CK}. Similar ideas about dropping exit dimensions are also implicit in the works~\cite{BM+,WBS}.

We will now recall precisely what the toy model systems (TMS) are  in~\cite{CK,GK,GHP}, and then present our generalized TMS for the Hamiltonian and non-Hamiltonian setting, as well as the results.

\subsection{TMS (toy model systems)}
\label{sec:genToyModel}
The study of growth of Sobolev norms for some solutions of Hamiltonian PDE's started in 2000, when Jean Bourgain~\cite{B}
posed the following question: are there solutions of the cubic nonlinear Schr\"odinger equation (NLS)
\[
-i u_t +\Delta u = \abs{u}^2u
\]
in $\mathbb{T}^2$ such that $\norm{u(t)}_{H^s} \longrightarrow \infty$ as $t \to\infty$?

The first significant answer was given ten years later by the I-team~\cite{CK}, proving an arbitrary growth for some solutions after a large time $T$.
One essential ingredient of their proof was to deal with a \emph{toy model system}
\[
%\label{GenNLSeq}
\dot b_\ell = -i \abs{b_\ell}^2 b_\ell + 2i\overline{b}_\ell\left(b_{\ell-1}^2 + b_{\ell+1}^2\right),
\]
defined on a finite, but arbitrarily large, set of complex Fourier modes ($1\leq\ell\leq N$) of the solutions of the NLS equation.

Estimates of the time $T$ of instability were later provided by Guardia and Kaloshin~\cite{GK} and then generalized to more general NLS systems~\cite{GHP} such as
\[
-i u_t +\Delta u = 2d\abs{u}^{2(d-1)}u + G'(\abs{u}^2)u, \quad d=2,3,\dots
\]
In the next section we will generalize these TMS, simply taking into account their common geometrical features.

\subsection{TMS in Hamiltonian setting}
%\label{subsec:more-general-Hamiltonian}

Let us consider a smooth Hamiltonian function $H(b)$,
also denoted as $H(b,\overline{b})$ where $\overline{b}$ denotes the complex conjugate to $b$, defined on the phase space
$\{b=(b_1,\dots,b_n)\in \mathbb{C}^n\}$ equipped with the standard complex symplectic structure
\[
 \Omega=\frac{i}{2}\sum_{l=1}^n db_l \wedge d\overline{b}_l,
\]
which for a Hamiltonian function $H$ defines the associated vector field $X_H$ by
\[
\label{eq:XH}
   \Omega(X_H(b),\eta)=dH(b)\eta, \quad \forall \eta \in \mathbb{C}^n.
\]
The associated equations of motion $\dot{b}=X_H(b)$ are expressed in components as
\[
\label{eq:HamMotion}
 \dot{b}_l = -2i \frac{\partial H}{\partial \overline{b}_l}, \quad 1\leq l\leq n.
\]

We will assume several hypotheses about $H$. Let us start with the first three.
\begin{enumerate}
\item[\textbf{H1}:] $H$ is a \emph{real smooth} Hamiltonian, that is, $H(b) \in \mathbb{R}$ for any $b\in\mathbb{C}^n$, and $H$ is $C^\infty$.
\item[\textbf{H2}:] $H$ is a \emph{phase invariant}:
\begin{equation}
\label{HPhaseInvariant}
H(\mathrm{e}^{i\theta} b, \mathrm{e}^{-i\theta} \overline{b})=H(b,\overline{b}), \quad \theta\in\mathbb{R}, b\in\mathbb{C}.
\end{equation}
\item[\textbf{H3}:] $H$ has \emph{invariant complex coordinate hyperplanes}, that is, each complex hyperplane $W_l=\{b\in \mathbb{C}^n:b_l=0\}$ is \emph{invariant} under the equations of motion of $H$ for any $l=1,\dots, n$.
\end{enumerate}
We will refer to any $b_l$-coordinate as a \emph{mode} (this terminology is inspired by the fact that in the context of NLSE these $b_l$-coordinates are the Fourier coefficients of the solutions).

By Hypothesis~\textbf{H1}, $\overline{H(b,\overline{b})}=H(b,\overline{b})$ and therefore
$\displaystyle \overline{\frac{\partial H}{\partial \overline{b}_l}}=\frac{\partial H}{\partial b_l}$ so that
\[
 \dot{\overline{b}}_l = 2i \frac{\partial H}{\partial b_l}, \quad 1\leq l\leq n.
\]

For simplicity, we are assuming that $H$ is $C^\infty$, although it would suffice to assume that $H$ is $C^r$ for some $r$ large enough to carry out a suitable reduction to normal form (see section~\ref{sub:polynomialNFc}).

Hypothesis~\textbf{H2} is equivalent to assuming that if $b(t)$ is a solution of~(\ref{eq:HamMotion}),
then $\mathrm{e}^{i\theta} b(t)$ is also a solution for any $\theta\in\mathbb{R}$, and is also equivalent to assuming that the
%After differentiation of eq.~(\ref{HPhaseInvariant}) with respect to $\theta$ we get that the
\emph{total mass} $\displaystyle M=\sum_{l=1}^n |b_l|^2$ is preserved by the equations of motion of the Hamiltonian~$H$, as a consequence of Noether's theorem.

Therefore, we will be able to restrict, in a first step, the phase space to the
unit $2n-1$ sphere $S^{2n-1}\subset\mathbb{C}^n$ given by $M=1$.  In a second step we can take into account
in this unit sphere the space of orbits for the action $(b_1,\dots,b_n) \mapsto (e^{i\theta}b_1,\dots,e^{i\theta}b_n)$
of the circle group $U(1)$,  and therefore restrict the phase space to the complex projective space
$\mathbb{C}\mathbb{P}^{n-1}=S^{2n-1}/U(1)$ that we will just denote by $\mathcal{M}$.
Let us notice that this is the underlying space for the toy models in~\cite{CK,GK,GHP}.

Hypothesis~\textbf{H3} is equivalent to assuming that $b_l=0\Longrightarrow \dot{b}_l =0$, for $l=1,\dots,n$. This condition was met in previous toy models~\cite{CK,GK,GHP},
where all the monomials in the Hamiltonian have even degree in each of the modes $(b_l,\overline{b}_l)$.
The invariance of the hyperplanes $W_l$ implies that for any $1\leq j\leq n$ the \emph{complex line}
\[
V_j = \{b\in \mathbb{C}^n:b_l =0,\ l\neq j \}=\bigcap_{l\neq j} W_l
\]
is also invariant. When we restrict it to $\mathcal{M}$ we get an invariant circle
\[
\mathbb{T}_j = V_j\cap \mathcal{M}=\{b:b_l =0, \ l \neq j,\ |b_j|=1 \}.
\]
The invariance of the hyperplanes $W_l$ also implies that for any pair of different indices $j\neq k$,  $1\leq j,k\leq n$, the \emph{complex plane}
\[
V_{j,k} = V_{k,j} = \{b:b_l =0, \ l \neq j,k \}=\bigcap_{l\neq j,k}W_l
\]
is invariant.
When we consider the restriction of Hamiltonian $H(b)$ just to the two modes $b_j$, $b_{k}$ generating the complex plane $V_{j,k}$, we get a real, phase invariant,
two-degrees-of-freedom Hamiltonian, which for $j<k$ reads as
\begin{equation}
\label{eq:Hjk}
H_{j,k}(b_j,b_{k})=H(0,\dots,0,b_j,\dots,0,\dots, b_{k},\dots,0), \quad (b_j,b_k)\in\mathbb{C}^2,
\end{equation}
whose equations of motion preserve the two complex lines $V_j$ and $V_{k}$
and the mass $M_{j,k}=|b_j|^2+|b_{k}|^2$ of the two modes $b_j$, $b_{k}$.
Therefore, $H_{j,k}$ is an \emph{integrable} Hamiltonian.
When we intersect the two complex lines $V_j$ and $V_{k}$
with the $3$-sphere $S^{3}\subset\mathbb{C}^2$ we get the two invariant circles
\[
\mathbb{T}_j = \{(b_j,0): |b_j|=1\} \text{ and }
\mathbb{T}_{k} = \{(0,b_{k}): |b_{k}|=1\},
\]
which in the next hypotheses will be assumed to be saddle periodic orbits within the complex plane  $V_{j,k}$ for adjacent modes $k=j-1, j+1$,
and elliptic periodic orbits for far (i.e., non-adjacent) modes $\abs{k-j}>1$.
\begin{enumerate}
\item[\textbf{H4}:] $\mathbb{T}_j$ and $\mathbb{T}_{j+1}$ are saddle periodic orbits
%in the same level of energy
of $H_{j,j+1}$ for $1\leq j\leq n-1$, and with the same characteristic exponents $\pm\lambda$ for all $j$.
\item[\textbf{H5}:] There exists at least one heteroclinic orbit of $H_{j,j+1}$ joining $\mathbb{T}_j$ to $\mathbb{T}_{j+1}$, for $1\leq j\leq n-1$.
% are joined by 4 heteroclinic orbits of $H_{j,j+1}$  for $1\leq j\leq n-1$.
\item[\textbf{H6}:] $\mathbb{T}_j$ and $\mathbb{T}_{k}$ are elliptic periodic orbits of $H_{j,k}$  for $\abs{k-j}>1$.

\end{enumerate}

Finally, as the use of an appropriate normal form near the periodic orbits will result in a shorter and more well-structured proof, we will further assume that the Hamiltonian has a \emph{sign symmetry}, which is equivalent to assuming that $H$ is even in each variable separately.

\begin{enumerate}
\item[\textbf{H7}:]
\emph{Sign symmetry}: $H$ is even in each variable separately: $H(s_1 b_1,\dots,s_n b_n)=H(b)$, where $s_l=\pm 1$ is an arbitrary sign for each $l=1,\dots,n$.
%or, equivalently, if $b(t)=(b_1(t),\dots,b_n(t)$ is a solution of the field Hamiltonian, so is $(s_1 b_1(t),\dots,s_n b_n(t)$, where for each $l=1,\dots,n$ $s_l=\pm 1$ is an arbitrary sign.
\end{enumerate}

After introducing these hypotheses, we can state our main result in this Hamiltonian setting:
\begin{theorem}
\label{thm:toymodel}
 Under Hypotheses \textbf{H1-H7}, for any $n\geq 3$  and for all $\varepsilon>0$ there exists a point $x_1$ close to $\mathbb{T}_1$ whose trajectory is $\varepsilon$-close
 to the  chain of heteroclinic connections $\mathbb{T}_1 \to \mathbb{T}_2 \to \cdots \to \mathbb{T}_n$.
\end{theorem}

\begin{rem}
From Hypotheses~\textbf{H4}-\textbf{H6}, it follows that all the circles $\mathbb{T}_j$ are periodic orbits in the same level of energy of the \emph{whole} Hamiltonian $H$.
By Hypothesis~\textbf{H4}, each intermediate periodic orbit $\mathbb{T}_j$, $1<j<n$, behaves as a saddle periodic orbit when we restrict it to the 2 complex planes $V_{j-1,j}$ or $V_{j,j+1}$, while by Hypothesis~\textbf{H6}, $\mathbb{T}_j$ behaves as an elliptic periodic orbit when we restrict it to any of the $n-3$ complex planes $V_{1,j},\dots, V_{j-2,j}, $ or $V_{j,j+2},\dots, V_{j,n}$. Therefore, as an orbit of the full Hamiltonian $H$, $\mathbb{T}_j$  is a saddle$^2$-center$^{n-3}$ periodic orbit.
\end{rem}
\begin{rem}
The assumption about the equality of the characteristic exponents in Hypothesis~\textbf{H4} holds in a natural way when $H_{j,j+1}$ is a symmetric function of its arguments: $H_{j,j+1}(b_j,b_{j+1})=H_{j,j+1}(b_{j+1},b_{j})$. As will be explained in next section~\ref{sec:Vjk}, Hypotheses~\textbf{H4}-\textbf{H6} hold naturally in the case of short-range interactions between a mode and its two \emph{adjacent} modes, that is, when $\displaystyle\frac{\partial H}{\partial \overline{b}_l}$ only depends on the modes $b_{l-1}$, $b_l$ and $b_{l+1}$, as is the case in~\cite{CK,GK,GHP}.
\end{rem}
\begin{rem}
Assumption~\textbf{H7} about sign symmetry is equivalent to assuming that if $b(t)$ is a solution of the Hamiltonian vector field, so is
$(s_1 b_1(t),\dots,s_n b_n(t))$ where $s_l=\pm 1$ is an arbitrary sign for each $l=1,\dots,n$. It is
also satisfied in the models of~\cite{CK,GK,GHP,Gu,Gi}.
Hypothesis~\textbf{H7} and the equality of the characteristic exponents in Hypothesis~\textbf{H4}, are assumptions of a more technical nature and more oriented towards simplifying the proofs than the rest of the hypotheses, and they are also assumed in this work because they are fulfilled in all known TMS that have been derived from a Hamiltonian PDE.
\end{rem}

The dynamical behavior of the periodic orbits $\mathbb{T}_j$ assumed in Hypotheses~\textbf{H4-H6} can be easily checked in
the complex planes $V_{j,k}$ where the Hamiltonian becomes integrable
(see conditions~(\ref{AdjacentModesCondition}-\ref{FarModesCondition})), as well as the existence of the heteroclinic connections between adjacent periodic orbits (see Figure~\ref{fig:4HeteroclinicOrbits}). We dedicate section~\ref{sec:Vjk} to this description.

\subsection{TMS in non-Hamiltonian systems}
%\label{subsec:more-general}
We now introduce a more general setting valid in the non-Hamiltonian context.
Let us consider a smooth $n$-dimensional complex vector field $f(b)$ which will also be denoted as $f(b,\overline{b})$, and with associated equations of motion
\begin{equation}
\label{fEquationsMotion}
\dot{b}=f(b,\overline{b}), \qquad b\in\mathbb{C}^n.
\end{equation}

We will assume several hypotheses about $f$. The first three are coming now.
\begin{enumerate}
\item[\textbf{f1}:] $f$ is a \emph{phase invariant}  smooth vector field:
\begin{equation}
\label{fPhaseInvariant}
f\left(\mathrm{e}^{i\theta}b,\mathrm{e}^{-i\theta}\overline{b}\right)=\mathrm{e}^{i\theta} f(b,\overline{b}), \qquad \quad \theta\in\mathbb{R},  b\in\mathbb{C}^n.
\end{equation}
%\item[\textbf{f2}:] The total mass $\displaystyle M=\sum_{l=1}^n |b_l|^2$ is preserved by the equations of motion of the vector field $f$.
\item[\textbf{f2}:] $f$ has \emph{invariant unit sphere}, that is, the unit sphere given by $M=1$ is preserved by the equations of motion of the vector field $f$, where
$\displaystyle M=\sum_{l=1}^n |b_l|^2$ is the total mass.
%\textbf{PZ: it is enough to have just $M=1$ to be invariant.}
\item[\textbf{f3}:] $f$ has \emph{invariant complex coordinate hyperplanes}, that is, each complex hyperplane $W_l=\{b\in \mathbb{C}^n:b_l=0\}$ is invariant under the equations of motion of $f$ for any $l=1,\dots, n$.
\end{enumerate}

Hypothesis~\textbf{f1} is equivalent to saying that if $b(t)$ is a solution of system~(\ref{fEquationsMotion}),
then $\mathrm{e}^{i\theta}b(t)$ is also a solution for any $\theta$.
For simplicity, we are assuming that $f$ is $C^\infty$, although it would suffice to assume that $f$ is $C^r$ for some $r$ large enough to carry out a suitable reduction to normal form (see section~\ref{sub:polynomialNFc}).

Since the vector field $f$ is no longer conservative, in principle it has no conserved quantity associated to its phase invariance~\eqref{fPhaseInvariant}, so we need to add the assumption that the total mass $M$ is preserved, at least for $M=1$, in Hypothesis~\textbf{f2} in order to have a geometric setting analogous to the Hamiltonian case.

As a consequence of Hypotheses~\textbf{f1} and~\textbf{f2} we are able to restrict ourselves first
to the unit $2n-1$ sphere $S^{2n-1}\subset\mathbb{C}^n$ given by $M=1$ and then to
$\mathcal{M}=\mathbb{C}\mathbb{P}^{n-1}=S^{2n-1}/U(1)$.

The invariance of the hyperplanes $W_l$ stated in condition~\textbf{f3} implies that any complex line
$V_j = \{b\in \mathbb{C}^n:b_l =0,\ l\neq j \}$ is invariant for $1\leq j\leq n$ as well as any
complex plane $V_{j,k} = V_{k,j} =  \{b:b_l =0, \ l \neq j,k \}$ for any two different indices $j\neq k$.
When we restrict $V_j$ to $M=1$ we get an invariant circle $\mathbb{T}_j =\{b:b_l =0, \ l \neq j,\ |b_j|=1 \}$ which becomes an equilibrium point in
$\mathcal{M}$.

When we consider the restriction of the vector field $f(b,\overline{b})$ to the two different modes $b_j$, $b_{k}$, $j\neq k$, generating $V_{j,k}$
just by imposing $b_l=0$ for $l\neq j, k$, we get the vector field
$f_{j,k}(b_j,b_{k})=f_{j,k}(b_j,b_{k},\overline{b}_j,\overline{b}_{k})$ with phase invariant equations of motion
\[
\dot{b}_j=f_j(b_j,b_{k}), \qquad \dot{b}_{k}=f_{k}(b_j,b_{k}),  \quad (b_j,b_k)\in\mathbb{C}^2,
\]
preserving the two complex lines $V_j$ and $V_{k}$
and the mass $M_{j,k}=|b_j|^2+|b_{k}|^2$ of the two modes $b_j$, $b_{k}$.
When we intersect the two complex lines $V_j$ and $V_{k}$
with the $3$-sphere $S^{3}\subset\mathbb{C}^2$ we get the two invariant circles
\[
\mathbb{T}_j = \{(b_j,0): |b_j|=1\} \text{ and }
\mathbb{T}_{k} = \{(0,b_{k}): |b_{k}|=1\},
\]
which in the next hypotheses will be assumed to be saddle periodic orbits for adjacent modes $k=j-1, j+1$,
and elliptic periodic orbits for far modes $\abs{k-j}>1$.
\begin{enumerate}
\item[\textbf{f4}:] $\mathbb{T}_j$ and $\mathbb{T}_{j+1}$ are saddle periodic orbits of the vector field $f_{j,j+1}$
with the same characteristic exponents of the form $\pm\lambda$ with $\lambda>0$.
\item[\textbf{f5}:] There exists at least one heteroclinic orbit of the vector field $f_{j,j+1}$ joinig $\mathbb{T}_j$ to $\mathbb{T}_{j+1}$.
\item[\textbf{f6}:] $\mathbb{T}_j$ and $\mathbb{T}_{k}$ are elliptic periodic orbits of the vector field $f_{j,k}$ for $\abs{k-j}>1$.
\item[\textbf{f7}:]
\emph{Sign symmetry}: $f_l(s_1 b_1,\dots,s_n b_n)=s_l f_l(b)$, where $s_l=\pm 1$ is an arbitrary sign for each $l=1,\dots,n$.
%or, equivalently, if $b(t)=(b_1(t),\dots,b_n(t)$ is a solution of the field Hamiltonian, so is $(s_1 b_1(t),\dots,s_n b_n(t)$, where for each $l=1,\dots,n$ $s_l=\pm 1$ is an arbitrary sign.
\end{enumerate}

We can now state our main result in this general setting:
\begin{theorem}
\label{thm:toymodel-gen}
 Assuming Hypotheses \textbf{f1-f7},
%  as well as block diagonal dynamics along the heteroclinic connections (\ref{eq:add-hete-in}-\ref{eq:add-hete-out}),
for any $n\geq 3$  and for all $\varepsilon>0$ there exists a point $x_1$ close to $\mathbb{T}_1$ whose trajectory is $\varepsilon$-close
 to the  chain of heteroclinic connections $\mathbb{T}_1 \to \mathbb{T}_2 \to \cdots \to \mathbb{T}_n$.
\end{theorem}

\begin{rem}
From Hypotheses~\textbf{f4} and~\textbf{f6}, it follows that all the circles $\mathbb{T}_j$ are periodic orbits of the whole vector field $f$.
Adding Hypothesis~\textbf{f2}, $\mathbb{T}_j$  is a saddle$^2$-center$^{n-3}$ periodic orbit of the full vector field $f$.
\end{rem}

\begin{rem}
The assumption about the equality of the characteristic exponents in Hypothesis~\textbf{f4} holds in a natural way when $f_{j,j+1}$ is a symmetric function of its arguments. As will be explained later in section~\ref{sec:Vjk}, Hypotheses~\textbf{f4}-\textbf{f6} may hold naturally in the case of short-range interactions between a mode and its two adjacent modes.
\end{rem}

\begin{rem}
Assumption~\textbf{f7} is equivalent to assuming that if $b(t)$ is a solution of the vector field $f$, so is
$(s_1 b_1(t),\dots,s_n b_n(t))$ where $s_l=\pm 1$ is an arbitrary sign for each $l=1,\dots,n$. 
\end{rem}

The dynamical behavior of the periodic orbits $\mathbb{T}_j$ assumed in Hypotheses~\textbf{f4-f6} can be easily checked in
the complex planes $V_{j,k}$ where the vector field has the first integral $M_j$.
We dedicate the next section to the description of motion in the phase space $V_{j,k}$.

After introducing the hypotheses and results, it is worth remarking that thanks to the hypotheses, a geometrical property called 
% a prointroducing these seven hypotheses, completely analogous to those assumed in the Hamiltonian setting, we still have to assume an %extra geometric condition, which is called
 the \emph{block diagonal dynamics along the heteroclinic connections} is satisfied, property which is necessary to apply the dropping dimensions mechanism already introduced in~\cite{DSZ}. This geometric condition amounts to the existence of a block diagonal structure in the variational equations of the vector field along the heteroclinic orbits, where the blocks are matrices over the complex lines $V_j$, as will be explained in section~\ref{subsec:Blockdiagonal} and computed explicitly in~\eqref{eq:decompTM} or~(\ref{eq:add-hete-in}-\ref{eq:add-hete-out}). This condition is equivalent to the invariance of the splitting of the tangent space in the complex lines $V_j$ along the heteroclinic orbits, see~\eqref{eq:decompTM}.  As shown in section~\ref{subsec:HamBlockDiagonal}, this geometric condition is automatically satisfied in Hamiltonian systems, and also holds automatically in the general systems thanks to Hypothesis~\textbf{f7}.

We finish this introduction detailing the structure of this paper.

In Section~\ref{sec:Examples} several examples of TMS satisfying the assumptions of the main theorems are introduced. In the Hamiltonian case, emphasis is placed on the simplicity and generality of TMS satisfying Hypotheses~\textbf{H1}-\textbf{H6}. Particularly the existence of heteroclinic connections of Hypothesis~\textbf{H5}, thanks to the consideration of the geometric properties of the Hamiltonians defined on the projective space complex. To the best of our knowledge, the presented examples contain all the previously considered TMS. In the non-Hamiltonian case, the existence of the heteroclinic connections of Hypothesis~\textbf{f5} is not satisfied in general, and specific conditions must be added to guarantee it. These examples of non-Hamiltonian TMS are entirely new.

Section~\ref{sec:Vjk} is devoted to explaining the geometrical and dynamical consequences of Hypotheses \textbf{H1}-\textbf{H6} and \textbf{f1}-\textbf{f6} of the main Theorems~\ref{thm:toymodel} and~\ref{thm:toymodel-gen}, respectively. Among them we stress the characterization of the complex projective space $\mathbb{C}\mathbb{P}^{n-1}$ as the phase space due to phase invariance, as well as the dynamical characterization of the periodic orbits $\mathbb{T}_j$ and existence of non-transversal heteroclinic connections between them, thanks to the invariance of the complex planes $V_{j,j+1}$ containing two adjacent periodic orbits $\mathbb{T}_j$ and $\mathbb{T}_{j+1}$.

Section~\ref{sec:Phase-space-local} contains a detailed study of the phase space near the invariant complex planes, thanks to the introduction of suitable coordinates, called \emph{$j$-charts}, centered on the periodic orbits $\mathbb{T}_j$.
The assumption about the block diagonal dynamics along the heteroclinic connections is introduced in several equivalent ways, and it is checked that it holds in the Hamiltonian setting. This condition will be necessary for shadowing heteroclinic connections when we use the dropping dimensions mechanism. The transition map between consecutive $j$-charts is also computed and an adequate polynomial normal form, which preserves the invariant hyperplanes as well as the block diagonal dynamics along the heteroclinic connections, is introduced around each periodic orbit. This normal form will be very convenient in later sections to obtain good estimates of shadowing trajectories.

The purpose of Section~\ref{sec:covrel} is to review the notions of h-sets and covering relations already introduced in~\cite{ZGi}, to state the main topological theorem~\ref{thm:top-shadowing} for finding shadowing trajectories using the dropping dimensions mechanism.

Section~\ref{sec:encl-full-system} is dedicated to presenting accurate estimates, called \emph{enclosures}, for the trajectories close to the periodic orbits. Its main result is Theorem~\ref{thm:ful-system}.

Finally, using all these tools developed in the previous sections: block diagonal structure of the coordinate changes involved, estimates for covering relations, and a polynomial normal form, the main Theorems~\ref{thm:toymodel} and~\ref{thm:toymodel-gen} are proved in section~\ref{sec:scheme}, thanks to the construction of a sequence of suitable covering relations along the heteroclinic chain $\mathbb{T}_1 \to \mathbb{T}_2 \to \cdots \to \mathbb{T}_n$ of the TMS presented in the Introduction.

%% file: examples.tex
\section{Several concrete examples}
\label{sec:Examples}

We now present a whole collection of simple systems that satisfy the above assumptions, and for which Theorems~\ref{thm:toymodel} or~\ref{thm:toymodel-gen} can be applied. We will detail in this section the most direct checks, leaving some of the more technical ones for the next section. It is important to note that the collection of Hamiltonian systems we present contains all previous examples known to the authors, while the collection of non-Hamiltonian systems is completely new.

\subsection{A quartic Hamiltonian}
Consider first a quartic homogeneous polynomial
\begin{equation}
\label{quartic}
H(b)=H(b,\overline{b})=\frac{1}{4}\sum_{\ell,m=1}^{n} \overline{b}_\ell^2 a_{\ell m} b_m^2,
\end{equation}
with associated Hamiltonian equations
\[
\dot b_\ell = -2i\frac{\partial H}{\partial \overline{b}_\ell}
=-i \left(\sum_{m=1}^n a_{\ell m}b_m^2\right)\overline{b}_\ell, \quad \ell=1,\dots,n,
\]
generated by a \emph{Hermitian} complex matrix $A=(a_{\ell m})_{1\leq j,\ell\leq n}$
with supra-sub diagonal elements more dominant that the diagonal ones, and these more dominant than the rest:
\begin{equation}
\label{eq:ExCond}
a_{\ell m}=\overline{a_{m\ell}} \text{ for } 1\leq \ell,m\leq n, \quad\text{ and }\quad
\abs{a_{\ell m}}<\abs{a_{\ell\ell}}<\abs{a_{\ell,\ell+1}} \text{ for } \abs{\ell-m}\geq 2.
\end{equation}
The Hermitian character of the matrix $A$ guarantees that $H(b)=H(\overline{b})$ is a real Hamiltonian, that is, Hypothesis~\textbf{H1}. It is clear that $H$ is phase invariant, has invariant complex coordinate hyperplanes and is even in each variable separately, so it satisfies Hypotheses \textbf{H2}, \textbf{H3} and \textbf{H7}.

Take any pair of different indices $j\neq k, 1 \leq j, k\leq n$. When we consider the restriction of Hamiltonian $H(b)$ to the
two modes $b_j$, $b_k$ generating the complex plane $V_{j,k}$, we get the real, phase
invariant, two-degrees-of-freedom Hamiltonian
\[
H_{j,k}(b_j,b_k)=\frac{1}{4}\left(\overline{b}_k^2 a_{kk} b_k^2+\overline{b}_k^2 a_{kj} b_j^2
+\overline{b}_j^2 \overline{a_{kj}} b_k^2 +\overline{b}_j^2 a_{jj} b_j^2  \right).
\]

%The change~(\ref{eq:j-chart-polar}) to polar variables $b_j=r e^{i\theta}$ is smooth only for $|b_j|=r>0$, i.e., for $| c|<1$, so it is not smoothly defined over $\mathbb{T}_k$ for $k\neq j$. However, the change~(\ref{Hj-chart}) can be extended continuously so that it blows up $\mathbb{T}_{k}$ to the circle $|c|=1$.
%In other words, as $H_{j,k}$ is defined on the complex projective line $\mathbb{C}\mathbb{P}^1=S^{3}/U(1)=S^2$, we then have that $H(c)$ is defined on the Riemann sphere, with the particularity that its two poles, for example, are saddle points of $H(c)$ in the same energy level, say $H=0$.

To reduce one degree of freedom, we now use the coordinates $b_j=r \mathrm{e}^{i\theta}$,  $b_k=c \mathrm{e}^{i\theta}$ of the $j$-chart introduced in Section~\ref{sec:Vjk},
restricted to the complex plane $V_{j,k}$, with $r^2=1-|c|^2$ to restrict also the phase space to
the sphere $M_{j,k}:=|b_j|^2+|b_k|^2=1$. The standard polar coordinates $r$, $\theta$ for the $j$-th mode are no longer present in the Hamiltonian $H_{j,k}$ in the variable $c$,
which takes the form
\begin{align}
H(c)&=\frac{1}{4}\left(a_{jj}\left(1-\abs{c}^2\right)^2+2\left(1-\abs{c}^2\right)\re\left(a_{kj}c^2\right)+a_{kk}\abs{c}^4\right)\notag\\
&= \frac{a_{kk}}{4}-\left(1-\abs{c}^2\right)H_2(c),\label{eq:H(c)}
%&\frac{\alpha}{4}-\frac{\alpha|c|^2}{2}+\frac{\mathrm{Re}(a_{kj} c^2)}{2}
%  +\frac{|c|^2}{2} \left(\alpha|c|^2-\mathrm{Re}(a_{kj} c^2)\right) \notag\\
%  &=\frac{\alpha}{4}-(1-|c|^2)\left(\frac{\alpha|c|^2}{2}-\frac{\mathrm{Re}(a c^2)}{2}\right) \label{eq:H(c)}\\
%  &=\frac{\alpha}{4}-(1-c \overline{c})H_2(c), \notag
\end{align}
where we have introduced, except for an additive constant, the general real quadratic homogeneous polynomial~\eqref{eq:H2}
\[
H_2(c)=H_2(c,\overline{c})=
\frac{a_{kk}-a_{jj}}{4}+\frac{a_{kk}+a_{jj}}{4}\abs{c}^2-\frac{\re\left(a_{kj}c^2\right)}{2}.
\]
Notice that the origin $c=0$ (which is the reduction of the periodic orbit $\mathbb{T}_j$ of $H_{j,k}$),
as well as the circle $\abs{c}=1$ (which is the reduction of the periodic orbit $\mathbb{T}_k$) are invariant under
the Hamiltonian $H(c)$ in the level sets $H=a_{kk}/4$ and $H=a_{jj}/4$, respectively.
%\frac{\alpha|c|^2}{2}-\frac{\mathrm{Re}(a c^2)}{2}
%=\frac{\alpha c \overline{c}}{2}-\frac{a c^2}{4}-\frac{\overline{a}\,\overline{c}^2}{4},\ \alpha\in\mathbb{R},\ a\in\mathbb{C}.
%\]

For two adjacent modes $k=j-1, j+1$, $H(c)$ takes the form, say for $k=j+1$,
\begin{equation}
\label{H2(c)j,j+1}
H(c)=\frac{a_{j+1,j+1}}{4}-(1-c \overline{c})H_2(c),
\end{equation}
where
\[
H_2(c)=
\frac{a_{j+1,j+1}-a_{jj}}{4}+\frac{a_{j+1,j+1}+a_{jj}}{4}\abs{c}^2-\frac{\re\left(a_{j+1,j}c^2\right)}{2}.
\]
By condition~\eqref{eq:ExCond} we have that
\[
\abs{a_{j+1,j}}^2>\left(\frac{a_{j+1,j+1}+a_{jj}}{2}\right)^2,
\]
which is equivalent to the condition~\eqref{AdjacentModesCondition} used in Section~\ref{subsec:HamiltonianCase} to prove that the origin is a saddle equilibrium of $H(c)$, with real characteristic exponents $\pm \lambda$ with
\[
\lambda=\sqrt{\abs{a_{j+1,j}}^2-\left(\frac{a_{j+1,j+1}+a_{jj}}{2}\right)^2}.
\]
Consequently, $\mathbb{T}_j$ is a saddle periodic orbit of the Hamiltonian $H_{j,j+1}$.

The level curves $\displaystyle H(c)=\frac{a_{j+1,j+1}}{4}$ in~equation~\eqref{H2(c)j,j+1} consist of the circle $|c|=1$ ($\mathbb{T}_{j+1}$) and the two straight lines $H_2(c)=0$. They will pass through the origin and intersect transversally the circle $|c|=1$, as long as $a_{j+1,j+1}=a_{j,j}=:\alpha$. In this case, all the periodic orbits will have the same characteristic exponent $\lambda$ as long that $a_{j+1,j}=:a$, which will imply both Hypotheses~\textbf{H4} and~~\textbf{H5}, with four straight heteroclinic orbits between $\mathbb{T}_j$ and $\mathbb{T}_{j+1}$, which jointly with the circle $|c|=1$ enclose four disjoint open regions.
Summarizing, we had to add the extra conditions that the three main diagonals of the matrix $A$ are constant
\begin{equation}
\label{eq:ExCondExtra}
a_{j,j}=:\alpha, \qquad a_{j+1,j}=a,
\end{equation}
to conditions~\eqref{eq:ExCond} to have Hypotheses~\textbf{H4} and~\textbf{H5} fulfilled.

\begin{rem}
The critical points of $H(c)$ are given by the cubic complex equation $\dfrac{\partial H}{\partial \overline{c}}=0$
which is equivalent to searching for the common zeros of two cubic real polynomial equations in the real variables $\mathrm{Re}\,c$ and $\mathrm{Im}\,c$, and by B\'ezout's theorem the number of common zeros is in general equal to 9, the product of the degrees. Since we already have 5 saddle points along the two straight lines $H_2(c)=0$, there is only room for 4 more critical points, each of them inside the 4
open regions between the circle $|c|=1$ and the two straight lines. Recalling that $H_{j,j+1}$ is defined on the complex projective line $\mathbb{C}\mathbb{P}^1=S^{3}/U(1)=S^2$, we then have that $H(c)$ is defined on the Riemann sphere and, by Morse theory or
Poincar{\'e}--Hopf theorem, these four extra critical points of $H(c)$ must be extrema points.
\end{rem}

For \emph{far} modes $|k-j|\geq 2$, the origin of $H(c)$ is an elliptic equilibrium point due to conditions~\eqref{eq:ExCond} and \eqref{eq:ExCondExtra},
which imply $|a_{kj}|<|\alpha|$, equivalent to~\eqref{FarModesCondition}, as explained in Section~\ref{subsec:HamiltonianCase},
so that $\mathbb{T}_j$ is an elliptic periodic orbit of $H_{j,j+1}$ and Hypothesis~\textbf{H6} is satisfied.

We finish this example by noticing that the toy model considered in~\cite{CK,GK,Gu,Gi} is a particular case of this example, with
a tridiagonal matrix $A$ satisfying $a_{\ell,\ell}=1$, $a_{\ell,\ell+1}=-2$ and $a_{\ell,m}=0$ for $\abs{\ell-m}\geq 2$.
See Fig.~\ref{fig:4HeteroclinicOrbits} for an illustration.

\subsection{Perturbed quartic Hamiltonian}
We can also add real perturbations to $H$ in~\eqref{quartic} as long as they preserve Hypotheses \textbf{H2}, \textbf{H3} and \textbf{H7}.
For instance, consider now the Hamiltonian
\[
H(b)=H(b,\overline{b})=\frac{1}{4}\sum_{\ell,m=1}^{n} \overline{b}_\ell^2 a_{\ell m} b_m^2
+\varepsilon\mathcal{P}\left(b,\bar{b},\varepsilon\right),
\]
with the matrix $(a_{\ell m})$ satisfying conditions~\eqref{eq:ExCond} and~\eqref{eq:ExCondExtra}, where $\mathcal{P}$ is real ($\mathcal{P}\left(b,\bar{b},\varepsilon\right)\in\mathbb{R}$ for all $b$),  polynomial in its variables $(b,\bar{b})$,
phase invariant, all the monomials in $\mathcal{P}$ are of even degree in each $(b_j,\bar{b}_j)$ (adding degrees for $b_j$ and $\bar{b}_j$), to
guarantee that $H$ has invariant complex coordinate hyperplanes $V_j$, and on each two adjacent invariant complex planes
$V_{j,j+1}$, the restricted Hamiltonian $H_{j,j+1}(b_j,b_{j+1})$ is a symmetric function of its arguments.

Then for $|\varepsilon|$ small enough the dynamical properties of the Hamiltonians $H_{j,\ell}$ are still preserved,
and Hypotheses~\textbf{H1}-\textbf{H7} hold.

One could also consider
\[
H(b)=H(b,\overline{b})= \left(\sum_{i=1}^n |b_i|^2 \right)^{d-2}\frac{1}{4}\sum_{\ell,m=1}^{n} \overline{b}_\ell^2 a_{\ell m} b_m^2
+\frac{1}{N}\mathcal{P}\left(b,\bar{b},\frac{1}{N}\right),
\]
with $N=2^{n-1}$, $d\in\mathbb{N}$, $d\geq 2$, as in~\cite{GHP}, and deal only with $N\gg 1$ to have a small perturbation of the quartic
Hamiltonian~\eqref{quartic}. This is the case considered in~\cite{GHP} for the same values
$a_{\ell,\ell}=1$, $a_{\ell,\ell+1}=-2$ and $a_{\ell,m}=0$ for $\abs{\ell-m}\geq 2$
as in~\cite{CK,GK,Gu,Gi}.
Notice that the factor $M^{d-2}=\left(\sum_i |b_i|^2 \right)^{d-2}$ of $H(b)$ does not introduce any change in the restricted Hamiltonian $H_{j,k}$ since our phase space is restricted to $M=1$.

\subsection{Non-Hamiltonian vector field}
Consider a vector field $f(b)=f(b,\overline{b})$ given in components by
\begin{equation}
\label{Non-HamiltonianVectorField}
\dot b_\ell = f_\ell(b,\overline{b})=
-i \left(\sum_{m=1}^n a_{\ell m}b_m^2\right)\overline{b}_\ell
+\rho \left(\sum_{j=1}^n C_{j\ell}b_j \overline{b}_j\right)b_\ell + \rho R b_\ell,
\end{equation}
generated by a complex matrix $(a_{\ell m})$ satisfying conditions~\eqref{eq:ExCond} and~\eqref{eq:ExCondExtra}, a \emph{real} matrix $(C_{j\ell})$, $(1\leq j,\ell\leq n)$,  a real parameter $\rho\geq 0$, and with
$R=R(b,\overline{b})=R(b \mathrm{e}^{i\theta},\overline{b}\mathrm{e}^{-i\theta})$ a phase invariant real function,
to be chosen to guarantee that $M=1$ is invariant.

It is clear that the vector field $f$ is phase invariant and has invariant complex coordinate hyperplanes, so it satisfies Hypotheses~\textbf{f1} and~\textbf{f3}.
%Notice that for $\rho=0$ and $R=0$, the vector field $f$ is the Hamiltonian vector field of the quartic real Hamiltonian $H$ in~\eqref{quartic}, and then the total mass %$M=\sum_{\ell=1}^n b_\ell \overline{b}_\ell$ is preserved by the vector field $f$.

A straightforward computation in system~\eqref{Non-HamiltonianVectorField} yields
\[
\dot M= 2\rho RM+ 2\rho \sum_{j,\ell=1}^n C_{j\ell}\,\overline{b}_j b_j\overline{b}_\ell b_\ell
 =2\rho R(M-1)+2\rho R+2\rho\,\sum_{j,\ell=1}^n C_{j\ell}\overline{b}_j b_j\overline{b}_\ell b_\ell.
\]
Therefore, choosing $R$ as the real quartic homogeneous polynomial
\[
R= -\sum_{j,\ell=1}^n C_{j\ell}\overline{b}_j b_j\overline{b}_\ell b_\ell,
\]
we have that $\dot M=2\rho R(M-1)$ so that $M=1$ is invariant and Hypothesis~\textbf{f7} is satisfied.
Notice that for $C_{j\ell}=0$ then $R=0$ and we recover a Hamiltonian system associated to Hamiltonian~\eqref{quartic}.
The same happens when the non-Hamiltonian parameter $\rho$ vanishes.
To have a genuine non-Hamiltonian system, we are going to assume $\rho>0$ and that $C_{j\ell}$
do \emph{not} vanish for all $j,\ell$.

For any pair of different indices $j\neq k, 1 \leq j, k\leq n$, when we consider the restriction $f_{jk}$ of the vector field $f$ to the
two modes $b_j$, $b_k$ generating the complex plane $V_{j,k}$, we get the two complex equations
\begin{align*}
\dot b_j&=-i\left(a_{jj} b_j^2 + a_{kj} b_k^2\right)\overline{b}_j+\rho\left(C_{jj} b_j\overline{b}_j + C_{kj} b_k \overline{b}_k\right)b_j +\rho Rb_j,\\
\dot b_k&=-i\left(a_{jk} b_j^2 + a_{kk} b_k^2\right)\overline{b}_k+\rho\left(C_{jk} b_j\overline{b}_j + C_{kk} b_k \overline{b}_k\right)b_k+\rho Rb_k,
\end{align*}
where
\[
R= -\left(C_{jj}b_j^2\overline{b}_j^2+\left(C_{jk}+C_{kj}\right)b_j\overline{b}_j b_k \overline{b}_k+C_{kk}b_k^2\overline{b}_k^2\right).
\]

If we use again the coordinates $b_j=r \mathrm{e}^{i\theta}$, $b_k=c \mathrm{e}^{i\theta}$ of the $j$-chart introduced in Section~\ref{sec:Vjk}, restricted to the complex plane $V_{j,k}$, plus $r^2=1-c\overline{c}$ to restrict also the phase space to
the sphere $M_{j,k}:=b_j\overline{b}_j+b_k\overline{b}_k=1$, the equations become
\begin{align}
\dot r +i r\dot\theta&=e_j+\rho Rr,\ e_j(r,c)=-i\left(a_{jj} r^2 + a_{kj} c^2\right)r+\rho\left(C_{jj} r^2 + C_{kj} c\overline{c}\right)r ,\label{rthetaEq}\\
\dot c +i c\dot\theta&=e_k+\rho Rc,\ e_k(r,c)=-i\left(a_{jk}r^2+a_{kk} c^2\right)\overline{c}+\rho\left(C_{jk} r^2 + C_{kk} c\overline{c}\right)c,\label{cEq}
\end{align}
where $R=R(c,\overline{c})$ is the real polynomial
\[
R= -\left(C_{jj}r^4+\left(C_{jk}+C_{kj}\right)r^2 c\overline{c}+C_{kk}c^2\overline{c}^2\right), \quad r^2= 1-c\overline{c}.
\]
Writing the real and imaginary part of equation~\eqref{rthetaEq} we get the o.d.e. for~$r$ and~$\theta$
\[
\dot r=\frac{e_j+\overline{e}_j}{2}+\rho Rr,\qquad \dot\theta=\frac{e_j-\overline{e}_j}{2i},
\]
from which we see that $r=0$ or, equivalently, $c\overline{c}=1$ ($\mathbb{T}_k$), is invariant for $f_{jk}$.

Replacing the expression above of $\dot\theta$ in equation~\eqref{cEq} we finally get a
complex \emph{reduced} vector field $f(c)$ for~$c$
\[
%\dot c= f(c)=f(c,\overline{c}):=ir^2\left(\alpha_j c-\hat{a}_{jk}\,\overline{c}\right)
%        +i\left(\mathrm{Re}\left(\hat{a}_{kj}c^2\right)-\alpha_k c\overline{c}\right)c+Rc,
\dot c= f(c)=f(c,\overline{c})=-\frac{e_j-\overline{e}_j}{2r}c+e_k+\rho Rc.
\]

Thanks to the real character of the coefficients $C_{jk}$, it turns out that $f(c)$ is the sum of a Hamiltonian system with a Hamiltonian $H(c)$ as in~\eqref{eq:H(c)} plus a \emph{radial} vector field of the form $g_k\, c$ with a real $g_k$ function given by
\[
g_k(c,\overline{c})=\rho\left(C_{jk} r^2 + C_{kk} c\overline{c}+R\right), \ r^2=1-c\overline c,
\]
so that \emph{any} invariant line of the Hamiltonian $H(c)$ is preserved.
Particularly for two adjacent modes $j,k$, say $k=j+1$, the two straight lines $H_2(c)=0$ passing through the origin ($\mathbb{T}_j$) and therefore intersecting transversally the circle $|c|=1$, which give rise to four straight heteroclinic orbits between $\mathbb{T}_j$ and $\mathbb{T}_{j+1}$. To keep the saddle characteristic exponents $\pm\lambda$ of the case $\rho=0$, one has to assume $C_{j,j}=C_{j,j+1}$. The four nodes inside the regions bounded by the straight heteroclinic orbits and the unit circle are typically attractor nodes for $\rho>0$.

\begin{figure}[hbt!]
\centering
\includegraphics[width=0.7\textwidth]{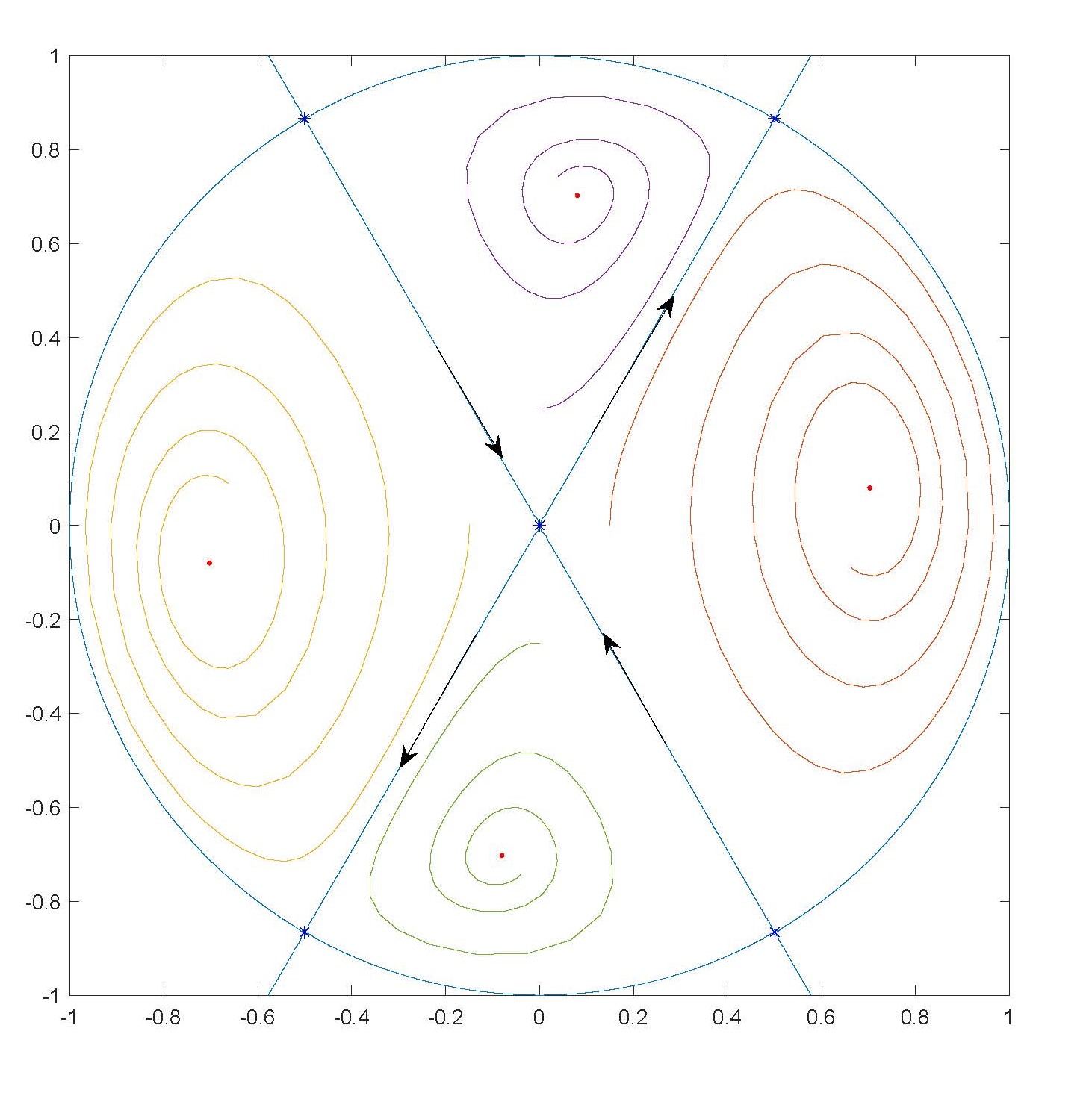}
\caption{Phase space $(\mathrm{Re}\,c, \mathrm{Im}\,c)$ of the reduced vector field $f(c)$ of the vector field $f_{j,j+1}$ derived from system~\eqref{Non-HamiltonianVectorField}, exhibiting 4 heteroclinic orbits. The values of the coefficients of the matrix are
$a_{\ell,\ell}=1$, $a_{\ell,\ell+1}=-2$ and $a_{\ell,m}=0$ for $\abs{\ell-m}\geq 2$ (as in the Hamiltonian considered in~\cite{CK,GK}), and $\rho=0.03$, $C_{j,j}=C_{j,j+1}=(-1)^{j}$, $C_{j+1,j}=2$ and $C_{j,k}=0$ for $\abs{j-k}\geq 2$, which has 5 saddle points and 4 attractor nodes on the Riemann sphere.}
\label{fig:NonHamiltonianExampleFigure}
\end{figure}

\begin{rem}
Due to the radial term $\rho Rc$ with $R$ a \emph{real} quartic polynomial in $c,\overline{c}$, the critical points of $f(c)$ are given by a quintic complex equation $f(c)=f(c,\overline{c})=0$, which is equivalent to searching for the common zeros of a system formed by a quintic and a cubic real polynomial equation in the real variables $\mathrm{Re} c$ and $\mathrm{Im} c$, and by B\'ezout's theorem the number of common zeros is in general equal to 15, the product of the degrees.
Nevertheless, for $\rho=0$ we are in the previous Hamiltonian setting with 5 saddle equilibria connected by heteroclinic orbits and 4 elliptic equilibria, so that for $\rho$ small enough, there will be also 5 saddle points and 4 nodes.
These situation remains for larger values, at least according to several numerical experiments. Finally notice that once we assume that $f(c)$ has 9 non-degenerate equilibria points in $c\overline{c}\leq 1$, since 5 of them are saddle points,the remaining four have to be nodes, according to Poincar\'e–Hopf theorem.
\end{rem}

Finally, we notice that due to the fact that all the monomials in the non-Hamiltonian part of the $\ell$-component of the vector field~\eqref{Non-HamiltonianVectorField} are of even degree in each $(b_j,\bar{b}_j)$ (adding degrees for $b_j$ and $\bar{b}_j$) for any $j\neq \ell$, it is clear that Hypothesis~\textbf{f7} holds. From this it is easy to check (see Section~\ref{subsec:HamBlockDiagonal}) that block diagonal dynamics
along the straight heteroclinic connections holds.

We have now all the conditions for a non-Hamiltonian vector field~\eqref{Non-HamiltonianVectorField} to satisfy Hypotheses~\textbf{f1}-\textbf{f7}, so the hypotheses of Theorem~\ref{thm:toymodel-gen} are satisfied.
See Fig.~\ref{fig:NonHamiltonianExampleFigure} for an illustration for $a_{\ell,\ell}=1$, $a_{\ell,\ell+1}=-2$ and $a_{\ell,m}=0$ for $\abs{\ell-m}\geq 2$ (the Hamiltonian considered in~\cite{CK,GK}), and $\rho=0.3$, $C_{j,j}=C_{j,j+1}=(-1)^{j}$, $C_{j+1,j}=2$ and $C_{j,k}=0$ for $\abs{j-k}\geq 2$.

%% file: genToyModel.tex
\section{Phase space in the invariant complex planes}
\label{sec:Vjk}
The phase space $\mathcal{M}$ is covered by the patches $\{b: b_j\neq 0\}$ where the $j$-th mode does not vanish, $1\leq j \leq n$.
To work out explicitly the (symplectic)  reduction from
$\mathbb{C}^n$ to $\mathcal{M}=\mathbb{C}\mathbb{P}^{n-1}=S^{2n-1}/U(1)$ we introduce, as in \cite{CK,GK,GHP}, coordinates centered on the $j$-th mode, defined by  rotating $b_j$ so that $b_j \in \mathbb{R}_+$, i.e., its phase angle is equal to $0$.
To be more precise, we  define coordinates on $\mathcal{M} \setminus \{b_j =0\}$, the $j$-\emph{chart},  as follows.

On $\{ b\in \mathbb{C}^{n}:b_j\neq 0\}$ we first introduce standard polar coordinates $r_j$, $\theta_j$ for the $j$-th mode
\begin{equation}
\label{eq:j-chart-polar}
  b_j=r_je^{i\theta_j}, \quad r_j > 0,
\end{equation}
and then the new coordinates $\{c_k\}_{k \neq j}$ for the rest of the modes are given by
\begin{equation}
  b_k=c_k e^{i \theta_j}.  \label{eq:j-chart-ck}
\end{equation}
Observe that  in this $j$-chart the variable $r_j$ corresponding to the $j$-th node
`disappears' but  it can be recovered from
\begin{equation}
  r_j = \sqrt{1 - \sum_{k\neq j} |c_k|^2}.
\label{eq:rj}\end{equation}
In particular the circle $\mathbb{T}_j$ defined by $r_j=1$ and $b_k=0$ for $k \neq j$ is mapped in the $j$-chart to the origin, i.e., $c_k=0$ for $k \neq j$.
Moreover, note that for $ 1 \leq k \neq j \leq n $ each invariant complex hyperplane $ W_k $ introduced in Hypothesis~\textbf{H3}
is defined simply by the equations $ c_k = 0 $ for $ k \neq j $.

\subsection{Hamiltonian case}
\label{subsec:HamiltonianCase}
After the change (\ref{eq:j-chart-polar}) and (\ref{eq:j-chart-ck}), the symplectic form reads as
\[
%\label{Omegac}
\Omega=\frac{i}{2}\sum_{k\neq j}dc_k \wedge d\overline{c}_k
+ \frac{1}{2}\sum_{k\neq j}\left(\overline{c}_k dc_k + c_k d\overline{c}_k \right)\wedge d\theta_j + r_j dr_j \wedge d\theta_j,
\]
and the associated Hamiltonian equations for the new Hamiltonian $H(c,\theta_j,r_j)$ are
\begin{eqnarray*}
\dot{c}_k &=& -i\dot\theta_j c_k -2i \frac{\partial H}{\partial \overline{c}_k},\\
\dot{\overline{c}}_k &=& \hphantom{-}i\dot\theta_j \overline{c}_k +2i \frac{\partial H}{\partial c_k},\\
\dot{\theta_j}&=&-\frac{1}{r_j}\frac{\partial H}{\partial r_j},\\
\dot{r_j}&=&\frac{1}{r_j}\left(\frac{\partial H}{\partial \theta_j}-\frac12 \sum_{k\neq j}\left( \overline{c}_k \dot{c}_k+c_k \dot{\overline{c}}_k\right)\right).
\end{eqnarray*}
Indeed, one can easily check that the vector field
\[
X_H= \sum_{k\neq j}\dot{c}_k \frac{\partial}{\partial c_k} + \sum_{k\neq j}\dot{\overline{c}}_k \frac{\partial}{\partial \overline{c}_k}
+\dot \theta_s \frac{\partial}{\partial \theta_j} + \dot r_j\frac{\partial}{\partial r_j}
\]
satisfies eq.~(\ref{eq:XH}).
If we now impose the phase invariance~(\ref{HPhaseInvariant}) of Hamiltonian~$H$,
the total mass $M$ is preserved and $\dfrac{\partial H}{\partial \theta_j}=0$.
Besides, if we restrict ourselves to $M=1$, the equation for $r_j$ can be clearly replaced by $r_j=r_j(c)$ as in~(\ref{eq:rj}).
As a consequence, in terms of the new Hamiltonian $H(c)=H(c,\theta_j,r_j(c))=H(c,r_j(c))$, the $n-1$ complex o.d.e. equations for ${c}_k$ simply read as
\[
\dot{c}_k = -2i \frac{\partial H}{\partial \overline{c}_k}, \quad 1\leq k\neq j \leq n,
%\label{eq:ck}
\]
that is, they are just the associated Hamiltonian equations for $H(c)=H(c,\overline{c})$ in the standard complex symplectic form on $\mathbb{C}^{n-1}$
\[
%\label{eq:Omegak}
\Omega=\frac{i}{2}\sum_{\genfrac{}{}{0pt}{1}{k=1}{k\neq j}}^{n} dc_k \wedge d\overline{c}_k.
\]

To see more clearly the local and global behavior around the periodic orbits,
let us fix one of them $\mathbb{T}_j$, take any other mode $k\neq j$ and write Hamiltonian $H_{j,k}$ introduced in eq.~(\ref{eq:Hjk}) in the coordinates
\[
(r=r_j\geq 0, \theta=\theta_j\in\mathbb{T},c=c_{k}\in\mathbb{C})
\]
of the $j$-chart (\ref{eq:j-chart-polar}-\ref{eq:j-chart-ck}) restricted to the complex plane $V_{jk}$ given by
\begin{equation}
\label{Hj-chart}
b_j=r e^{i\theta},  \quad b_{k} =c e^{i \theta}.
\end{equation}
As $H_{j,k}$ is phase invariant, $M_{j,k}=|b_j|^2+|b_{k}|^2$ is a first integral  and $\displaystyle \frac{\partial H_{j,k}}{\partial \theta}=0$ so $H_{j,k}$ depends only on the variables $r$ and $c$. Using $r^2=1-|c|^2$, we get the reduced Hamiltonian $H(c)=H(c,\overline{c})$ with associated equation of motion  $\displaystyle \dot{c} = -2i \frac{\partial H}{\partial \overline{c}}$ that preserves the origin $c=0$ ($\mathbb{T}_j$) and the circle
$|c|=1$ ($\mathbb{T}_{k}$). Notice that the invariance of $r=0$ and $c=0$ is equivalent to the invariance of $b_j=0$ and $b_k=0$, respectively, and a consequence of the invariance of the coordinate complex hyperplanes of $H$ provided by Hypothesis~\textbf{H3}.

\begin{rem}
The change~(\ref{eq:j-chart-polar}) to polar variables $b_j=r e^{i\theta}$ is smooth only for $|b_j|=r>0$, i.e., for $| c|<1$, so it is not smoothly defined over $\mathbb{T}_k$ for $k\neq j$. However, the change~(\ref{Hj-chart}) can be extended continuously so that it blows up $\mathbb{T}_{k}$ to the circle $|c|=1$.
In other words, as $H_{j,k}$ is defined on the complex projective line $\mathbb{C}\mathbb{P}^1=S^{3}/U(1)=S^2$, we then have that $H(c)$ is defined on the Riemann sphere, with the particularity that its two poles, for example, are saddle points of $H(c)$ in the same energy level, say $H=0$. In the $c$-coordinate, the north pole, for example, is at the origin $c=0$ and the south pole is at the circle $|c|=1$.

\end{rem}

As the origin ($\mathbb{T}_j$) of the reduced Hamiltonian $H(c)$ is an equilibrium point, we can write $H(c)=H_2(c)+O_3(c)$, where $H_2(c)$ is a real quadratic form in the variables $c, \overline{c}$:
\begin{equation}
\label{eq:H2}
H_2(c)= \frac{\alpha |c|^2}{2} - \frac{\text{Re}\,(a c^2)}{2}= \frac{\alpha c \overline{c}}{2} - \frac{a c^2+ \overline{a}\,\overline{c}^2}{4},
\quad \alpha\in\mathbb{R},\ a \in\mathbb{C}.
\end{equation}

The linear equation of motion associated to $H_2(c)$ is
\begin{equation}
\label{eq:cdot}
\dot c = -2i \frac{\partial H_2}{\partial \overline{c}}=-i(\alpha c - \overline{a}\, \overline{c})
\end{equation}
or, taking the real and imaginary part of $c$, as
\begin{equation}
\label{eq:B}
\frac{\mathrm{d}}{\mathrm{d}t}\left(\begin{array}{l}\text{Re}\, c\\ \text{Im}\,c\end{array}\right)
=\left(\begin{array}{cc}\text{Im}\,a&\alpha+\text{Re}\,a\\
                        -\alpha+\text{Re}\,a&-\text{Im}\,a\end{array}\right)
\left(\begin{array}{l}\text{Re}\,c\\ \text{Im}\, c\end{array}\right)
=B \left( \begin{array}{l}\text{Re}\,c\\ \text{Im}\, c\end{array}\right),
\end{equation}
with $\text{tr}\,B=0$, $\text{det}\,B=\alpha^2 -|a|^2$, so that the eigenvalues of the matrix $B$ are $\pm\lambda$ where
\[
%\label{eq:lambda}
\lambda=\sqrt{|a|^2-\alpha^2},
\]
and the stability of this system is totally determined by the sign of $|a|^2-\alpha^2$.

For two adjacent modes $k=j-1, j+1$ the condition
\begin{equation}
\label{AdjacentModesCondition}
|a|^2-\alpha^2>0
\end{equation}
provides real eigenvalues $\pm \lambda$ with $\lambda=\sqrt{|a|^2-\alpha^2}$ which give rise to  a saddle equilibrium at the origin of $H(c)$  and, consequently, a saddle periodic orbit $\mathbb{T}_j$ of Hamiltonian $H_{j,j+1}$. Therefore, condition~\eqref{AdjacentModesCondition} is the one needed to guarantee the first part of hypothesis~\textbf{H4}.

For far modes $\abs{k-j}>1$ the condition
\begin{equation}
\label{FarModesCondition}
|a|^2-\alpha^2<0
\end{equation}
provides pure imaginary eigenvalues $\pm i \nu$,
with $\nu=\sqrt{\alpha^2-|a|^2}>0$, which give rise to a center equilibrium at the origin of the reduced Hamiltonian $H(c)$ and,
consequently, an elliptic periodic orbit $\mathbb{T}_j$ of Hamiltonian $H_{j,k}$. Therefore, condition~\eqref{FarModesCondition} is the one needed to guarantee hypothesis~\textbf{H6}.

\begin{rem}
In the quadratic Hamiltonian $H_2(c)$ given in~(\ref{eq:H2}), the real coefficient $\alpha$
gives a measure of the inner energy of the $j$-th mode in the complex plane $V_{jk}$, whereas the complex parameter $a$ gives a measure of the interaction energy between the $j$-th mode and the $k$-th mode. In the so-called \emph{short-range interactions}, $a=0$ for $k$ away from $j$, so the $j$-th mode behaves as a center in the direction pointing to $k$-th mode. On the other hand, if the interaction energy between $j$-th mode and an adjacent mode $k=j-1,j+1$ is large enough compared with the inner energy of the $j$-th mode $(|a|>|\alpha|)$, the $j$-mode possesses a saddle behavior in the direction of the adjacent mode $k=j-1,j+1$. Notice that if the interaction energy between the $j$-th mode and an adjacent mode $k=j-1,j+1$ is not large enough $(|a|<|\alpha|)$, the $j$-mode still behaves like a center in the direction of the $k$-mode.
\end{rem}
\begin{rem}
The coefficients $\alpha$ and $a$ depend on $j$ and $k$. For adjacent modes $k=j+1$, when $H_{j,j+1}$ is a symmetric function of its arguments:
$H_{j,j+1}(b_j,b_{j+1})=H_{j,j+1}(b_{j+1},b_{j})$, $\alpha$ and $a$ are the same for all the adjacent modes, and consequently the same happens to the
characteristic exponent of all the periodic orbits $\mathbb{T}_j$ in the direction of the next mode $j+1$.
\end{rem}

In any case, when $\alpha^2\neq |a|^2$ the linear equation of motion~(\ref{eq:cdot}) associated to $H_2(c)$ can be easily diagonalized.
For the adjacent modes $k=j-1,j+1$ where the condition $|a|^2>\alpha^2$ holds, the matrix $B$ has real eigenvalues $\pm \lambda$ with $\lambda=\sqrt{|a|^2-\alpha^2}>0$
and one can try directly eigenvectors $\omega$, $\bar{\omega}$ of eq.~(\ref{eq:cdot}) of the form
\[
%\label{eq:omega}
\omega=\mathrm{e}^{i\vartheta}, \text{ with }\omega^2=\mathrm{e}^{2i\vartheta}=\frac{i\bar{a}}{\lambda+i\alpha}.
\]
%\textbf{PZ: in the formula above (\ref{eq:omega}) I would change $\theta$ to something else, for example $\vartheta$ or $\eta$, to avoid collision with $\theta$ being the %angular coordinate. The formula for $\omega^2$ is wrong. It should be $\omega^2=\frac{-i\bar{\beta}}{\lambda+i\alpha}$}
Then the conformal change of variables
\begin{equation}
\label{Conformal}
c=\omega x +\bar \omega y
\end{equation}
transforms the complex variable $c$ to 2 real variables $z=(x,y)$ in such a way that
diagonalizes eq.~(\ref{eq:cdot})
\begin{eqnarray*}
\dot x&=&\hphantom{-}\lambda x,\\
\dot y&=&-\lambda y,
\end{eqnarray*}
which is a Hamiltonian system with associated Hamiltonian
$H_2=\lambda xy$ with respect to the symplectic form
$\Omega = dx \wedge dy$.

For far modes $|k-j|>1$ where the condition $\alpha^2>|a|^2$ holds, the matrix $B$ has pure imaginary eigenvalues $\pm i\nu$ with $\nu=\sqrt{\alpha^2-|a|^2}>0$.
%and  \textbf{PZ: I would remove the formula below, as it is not immediate to see it, and in fact it is not necessary, we should write just the final result of the coordinate change}
%\[
%w=\left(\begin{array}{c}\text{Im}\beta\\-\alpha-\text{Re} \beta\end{array}\right)
%+i \left(\begin{array}{r}\nu\\0\end{array}\right)=: u + i v
%\]
%is an associated eigenvector to the eigenvalue $i\nu$: $Bw=i\nu w$, as well as any other complex  vector $\gamma w$ for any non-zero complex number $\gamma$.
Denoting $w=u+iv$ an associated eigenvector to the eigenvalue $i\nu$ of the matrix $B$: $Bw=i\nu w$,
and introducing the $2\times 2$ real matrix $U=(u | v)$ formed by the real vectors $u$ and $v$ as columns, the real change of variables
 \begin{equation}
\label{NonConformal}
\left(\begin{array}{l}\text{Re}\, c\\ \text{Im}\,c\end{array}\right)=U\left(\begin{array}{c}X\\Y\end{array}\right)
\end{equation}
transforms eq.~(\ref{eq:B}) to
\begin{eqnarray*}
\dot X&=&-\nu Y,\\
\dot Y&=&\hphantom{-}\nu X,
\end{eqnarray*}
which introducing the complex variable $C=X+iY$ can be just written as
\[
\dot C= i\nu C
\]
so that the linear equation~(\ref{eq:cdot}) has been diagonalized. Choosing adequately the factor $\gamma$, the change~(\ref{NonConformal}) transforms
Hamiltonian $H_2(c)$ to $\displaystyle H_2(C)=-\frac{\nu |C|^2}{2} $ and satisfies $\displaystyle \frac{i}{2}dc\wedge d\bar c=\frac{i}{2}dC\wedge d\bar C$.

\begin{figure}[hbt!]
\centering
\includegraphics[width=\textwidth]{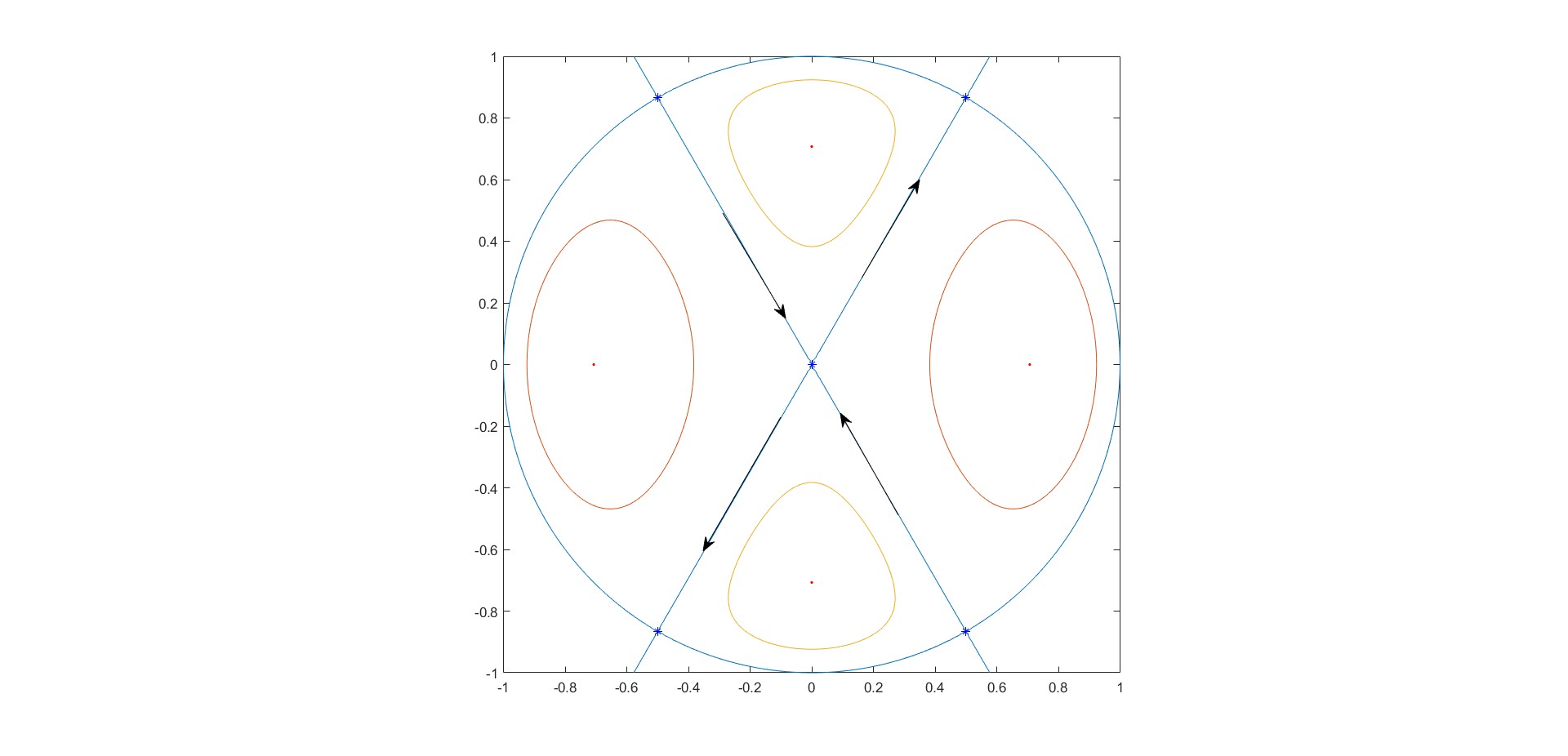}
\caption{4 heteroclinic orbits for the reduced Hamiltonian $H(c)$ of $H_{j,j+1}$, which has 5 saddle points and 4 extrema on the Riemann sphere.}
\label{fig:4HeteroclinicOrbits}
\end{figure}

\begin{figure}[hbt!]
\centering
\includegraphics[width=0.8\textwidth]{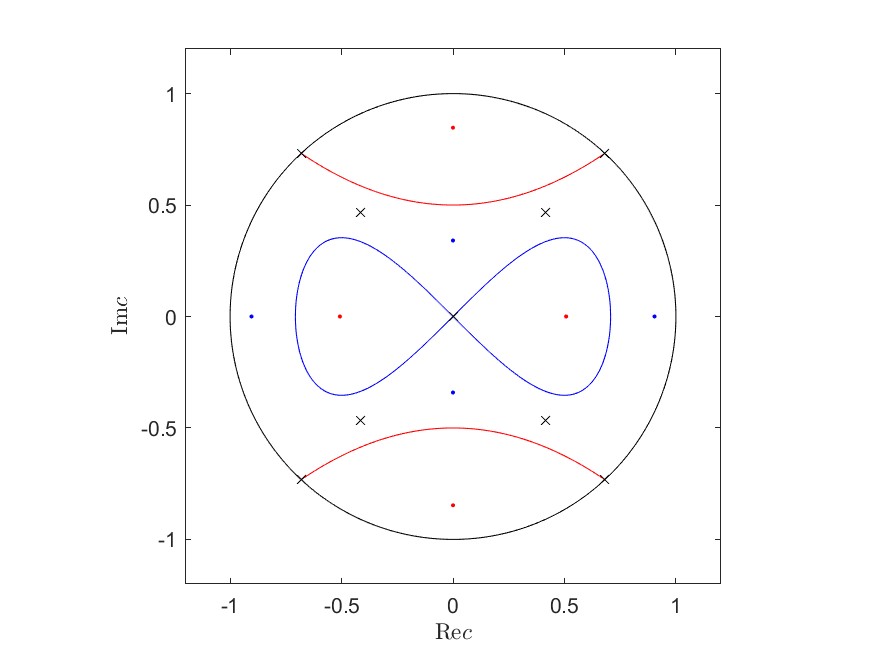}
\caption{2 possible heteroclinic orbits for the reduced Hamiltonian $H(c)$ of $H_{j,j+1}$, which has at least 5 saddle points and may have more than 4 extrema on the Riemann sphere.}
\label{fig:2HeteroclinicOrbits}
\end{figure}

According to hypothesis~\textbf{H5}, the phase portrait for the oriented levels of energy of $H_{j,j+1}$ is as the one depicted in Fig.~\ref{fig:4HeteroclinicOrbits} or Fig.~\ref{fig:2HeteroclinicOrbits} where there exist four or two heteroclinic orbits between $\mathbb{T}_j$ and $\mathbb{T}_{j+1}$, all of them sharing the same energy.

\begin{figure}[hbt!]
\centering
\includegraphics[width=0.8\textwidth]{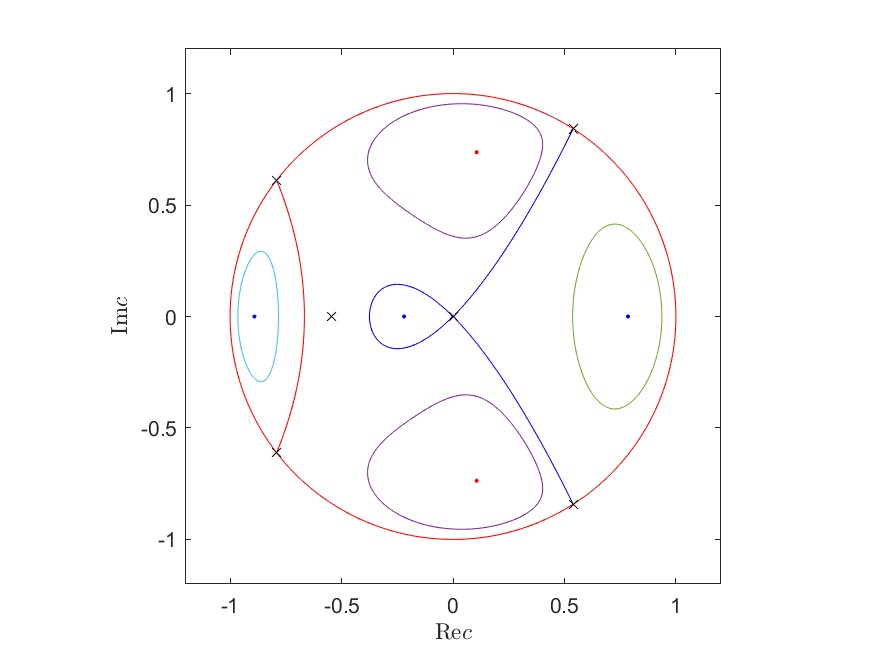}
\caption{Homoclinic orbits for the reduced Hamiltonian $H(c)$ of $H_{j,j+1}$, which has at least 5 saddle points and more than 4 extrema on the Riemann sphere. }
\label{fig:4HomoclinicOrbits}
\end{figure}

\begin{rem}
A different phase portrait is depicted in Fig.~\ref{fig:4HomoclinicOrbits} with homoclinic but no heteroclinic orbits (which is the standard case when $\mathbb{T}_j$ and $\mathbb{T}_{j+1}$ lie on \emph{different} levels of energy).
This phase portrait cannot take place neither in \cite {CK, GK} where the heteroclinic orbits are straight segments, nor in \cite {GHP} which is only a perturbation of the previous ones.
\end{rem}

The heteroclinic orbits of $H_{j,j+1}$ lie along two different invariant curves, the \emph{unstable} and \emph{stable} invariant curves $C^{\mathrm{u}}$, $C^{\mathrm{s}}$ of the saddle $\mathbb{T}_j$, that can be straightened in adequate symplectic coordinates. Indeed, by the conformal change of variables~(\ref{Conformal})
$H(c)=H_2(c)+O_3(c)$ is transformed to $H(x,y)=\lambda x y + O_3(x,y)$ with associated real canonical equations of motion $z=J\nabla H(z)$ where
$\renewcommand{\arraystretch}{0.5}
J=J_2=\left(\begin{array}{@{}r@{\hskip 2pt}r@{}}0&1\\-1&0\end{array}\right)$.
\begin{rem}
\label{HamStraighten}
By the unstable invariant manifold Theorem, in these $(x,y)$-coordinates the unstable curve can be locally written
as $y=\xi^{\mathrm{u}}(x)$, where $\xi^{\mathrm{u}}(x)=\nabla_x F^{\mathrm{u}}(x)$, since the unstable manifold is \emph{isotropic}, and extended globally as it is contained in the stable curve to the saddle $\mathbb{T}_{j+1}$. Equivalently, the generating function $S^{\mathrm{u}}(x,Y)=xY+F^{\mathrm{u}}(x)$ generates a canonical change $(x,y)\mapsto (X,Y)$ through the equations $y=\nabla_x S^{\mathrm{u}}(x,Y)$, $X=\nabla_Y S^{\mathrm{u}}(x,Y)=x$. Notice that $H(x,\nabla_x S^{\mathrm{u}}(x,Y))=0$. In the canonical variables
$(x,Y)=(x,y-\xi^{\mathrm{u}}(x))$, the local unstable curve is contained in the line $Y=0$. A latter canonical change $(X,Y)=(x-\xi^{\mathrm{s}}(Y),Y)$ brings also the
local stable curve to the line $X=0$. Consequently, the heteroclinic orbits of $H_{j,j+1}$ lie along two different straight lines in the canonical variables $(X,Y)$.
\end{rem}
\begin{rem}
%\label{BirkhoffNF}
Alternatively, the Birkhoff normal form around a saddle equilibrium of an \emph{analytical} one-degree-of-freedom Hamiltonian is convergent~\cite{M}, so there exist a canonical transformation $(x,y)=(X,Y)+O_3(X,Y)$ such that $H(X,Y)=F(J)=\lambda J +O_2(J)$, where $J=XY$. This transformation is well defined in a neighborhood of the saddle equilibrium point $(\mathbb{T}_{j})$ and can be extended along the unstable curve~\cite{SOD} to an arbitrary small distance of the circle $|c|=1$ $(\mathbb{T}_{j+1})$. The unstable invariant manifold is given by $Y=0$ and the stable invariant manifold by $X=0$.
\end{rem}

We summarize the effects of the changes performed in this section.
Introducing the subindex $_-$ for the previous $j-1$ mode and $_+$ for the next $j+1$ mode in the $j$-chart~(\ref{eq:j-chart-polar}-\ref{eq:j-chart-ck}), after diagonalizing by the linear changes (\ref{Conformal}) for the adjacent modes and the linear changes (\ref{NonConformal}) for the rest of the modes, and stretching the heteroclinic orbits by Remark~\ref{HamStraighten}, in the canonical coordinates
$(X,Y,C_*)=(X_-,X_+,Y_-,Y_+,C_1,\dots,C_{j-2}, C_{j+2},\dots,C_n)\in\mathbb{R}^2\times\mathbb{R}^2\times\mathbb{C}^{n-3}$
with respect to the symplectic form
\[
\Omega = dX_- \wedge dY_- +dX_+ \wedge dY_+ + \frac{i}{2}\sum_{k \neq j-1,j,j+1} dC_k \wedge d\overline{C}_k,
%\label{eq:Omega-loc},
\]
the Hamiltonian reads as
\[
H=\lambda X_- Y_- + \lambda X_+ Y_+ -\sum_{k \neq j-1,j,j+1}\frac{\nu_k |C_k|^2}{2} + O_3(X,Y,C_*),
%\label{eq:H-Diag}
\]
and has $n-3$ invariant complex hyperplanes $W_k$ given by the equations $c_k=0$ for $k=1,\dots,j-2,j+2,\dots,n$.
In the next section we will work in these variables and just move from uppercase to lowercase letters.

\subsection{Non-Hamiltonian case}
To see more clearly the local and global behavior around the periodic orbits in the non-Hamiltonian case,
fix $\mathbb{T}_j$, take any other mode $k\neq j$ and write vector field $f_{j,k}$ in the coordinates
\[
(r=r_j\geq 0, \theta=\theta_j\in\mathbb{T},c=c_{k}\in\mathbb{C})
\]
of the $j$-chart given in~(\ref{Hj-chart}).
As $f_{j,k}$ is phase invariant, the equations of motion depend only on the variables $r$ and $c$.
As the unit sphere $r^2+\abs{c}^2=\abs{b_j}^2 + \abs{b_k}^2=1$ is invariant, we can isolate $r^2=1-|c|^2$, to get the reduced vector field $f(c)=f(c,\overline{c})$ with associated equation of motion
\begin{equation}
\label{eq:cField}
\dot{c} = f(c,\overline{c})
\end{equation}
that preserves the origin $c=0$ ($\mathbb{T}_j$) and the circle
$|c|=1$ ($\mathbb{T}_{k}$).

%\begin{figure}[hbt!]
%\centering
%\includegraphics[width=0.8\textwidth]{4GeneralHeteroclinicOrbits}
%\caption{4 heteroclinic orbits for the reduced vector field $f(c)$ of the vector field $f_{j,j+1}$, which has 2 saddle points and 4 sources or sinks on the Riemann sphere.}
%\label{fig:4GeneralHeteroclinicOrbits}
%\end{figure}
Due to Hypothesis~\textbf{f6}, for $|k-j|>1$ the origin is a center equilibrium point of $f_{j,k}$ with characteristic exponents $\pm i \nu$ with $\nu>0$.
By Hypothesis~\textbf{f4}, for $k=j-1,j+1$  the origin is a saddle equilibrium point with characteristic exponents $\pm\lambda$ with $\lambda>0$.
A totally analogous computation can be performed around $\mathbb{T}_{j+1}$, providing a saddle equilibrium point with characteristic exponents $\pm\lambda$ with $\lambda>0$ on
$\mathbb{T}_{j+1}$. In the coordinates of the $j$-chart, this gives rise to 4 saddle equilibrium points in the circle $|c|=1$.
According to hypothesis~\textbf{f5}, the phase portrait for the heteroclinic orbits of $f_{j,j+1}$ is as the one depicted in
%Fig.~\ref{fig:4GeneralHeteroclinicOrbits}
Fig.~\ref{fig:NonHamiltonianExampleFigure} where there exist four heteroclinic orbits between $\mathbb{T}_j$ and $\mathbb{T}_{j+1}$.

The heteroclinic orbits of $f_{j,j+1}$ lie along two different invariant curves, the \emph{unstable} and \emph{stable} invariant curves $C^{\mathrm{u}}$, $C^{\mathrm{s}}$ of the saddle $\mathbb{T}_j$, that can be straightened in adequate coordinates. Indeed, one performs a linear change of variables
$z=(x,y)\in\mathbb{R}^2 \mapsto c=\omega x +\overline{\omega} y\in\mathbb{C}^2$,
which transforms equation~(\ref{eq:cField}) to $\dot x =\lambda x+O_2(x,y)$, $\dot y =-\lambda y+O_2(x,y)$. In other words, this change transforms $f(c)$
to $f(x,y)=(x+O_2(x,y),-y+O_2(x,y))$.
\begin{rem}
\label{Straighten}
By the unstable invariant manifold Theorem, in $(x,y)$-coordinates the local unstable curve can be locally written
as $y=\xi^{\mathrm{u}}(x)$ and extended globally as it is contained in the stable curve to the saddle $\mathbb{T}_{j+1}$. In the variables
$(x,Y)=(x,y-\xi^{\mathrm{u}}(x))$, the local unstable curve is contained in the line $Y=0$. A latter change $(X,Y)=(x-\xi^{\mathrm{s}}(Y),Y)$ brings also the
local stable curve to the line $X=0$. Consequently, the heteroclinic orbits of $H_{j,j+1}$ lie along two different straight lines in the variables $(X,Y)$.
\end{rem}

\section{Phase space close to the invariant complex planes}
\label{sec:Phase-space-local}
\subsection{Motion in  the $j$-chart}
%\label{subsec:TMS-motion-j-chart}

For any $j=1,\dots,n$, in the $j$-chart~(\ref{eq:j-chart-polar}-\ref{eq:j-chart-ck})  we introduce the subindex $_-$ for the previous $j-1$ mode and $_+$ for the next $j+1$ mode. After diagonalizing by the linear changes (\ref{Conformal}) for the adjacent modes and the linear changes (\ref{NonConformal}) for the rest of the modes,
%In the chart centered on the $j$-th mode, which is symplectic in the Hamiltonian case
%\begin{equation}
%\label{coordinatesj}
%\begin{array}{l}b_j=r_j e^{i\theta_j},\quad b_k=c_k e^{i\theta_j} \text{ for } k\neq j,\text{ and }\\[3pt]
%c_{j-1}=\omega x_- +\overline{\omega} y_-,\quad  c_{j+1}=\omega x_+ +\overline{\omega} y_+,
%\end{array}
%\end{equation}
%where $(\omega,\overline{\omega})\in\mathbb{C}^2$ are introduced in~eq.(\ref{Conformal}),
the equations of motion take the following local form in the real-complex coordinates
$(x,y,c_*)=(x_-,x_+,y_-,y_+,c_1,\dots,c_{j-2}, c_{j+2},\dots,c_n)\in\mathbb{R}^2\times\mathbb{R}^2\times\mathbb{C}^{n-3}$:
\begin{eqnarray}
\dot x_-&=&\hphantom{-}\lambda x_-+O_2(x,y,c_*),\notag\\
\dot y_-&=&-\lambda y_-+O_2(x,y,c_*),\notag\\
\dot x_+&=&\hphantom{-}\lambda x_+ +O_2(x,y,c_*),\label{eq:sys-res-0-full}\\
\dot y_+&=&-\lambda y_+ +O_2(x,y,c_*),\notag\\
\dot{c}_k &=& \ i\nu_k c_k +O_2(x,y,c_*), \quad \text{ for } k\neq j-1,j,j+1,\notag % \label{eq:sr-ck}
\end{eqnarray}
where $\lambda$ and $\nu_k$ are non-zero real numbers.
Notice that the linear changes~(\ref{Conformal}-\ref{NonConformal}) preserve the sign symmetry~\textbf{H7} or~\textbf{f7}.

In the Hamiltonian case equations~\eqref{eq:sys-res-0-full} are associated to the Hamiltonian
\[
H=\lambda x_- y_- + \lambda x_+ y_+ -\sum_{k \neq j-1,j,j+1}\frac{\nu_k |c_k|^2}{2} + O_3(x,y,c_*),
%\label{eq:H-diag}
\]
with respect to the symplectic form
\[
  \Omega = dx_- \wedge dy_- +dx_+ \wedge dy_+ + \frac{i}{2}\sum_{k \neq j-1,j,j+1} dc_k \wedge d\overline{c}_k.
%\label{eq:omega-loc}
\]

Moreover, as a consequence of the fact that the complex hyperplanes $W_k$ of Hypothesis~\textbf{H3} or~\textbf{f3} are invariant,
the real-complex hyperplanes $W_-=\{(x,y,c_*): x_+=y_+=0\}$ and $W_+=\{(x,y,c_*): x_-=y_-=0\}$,
as well as the complex hyperplanes $W_k=\{(x,y,c_*): c_k=0\}$, for $k \neq j-1,j,j+1$, are invariant. Taking intersections, the real planes
$V_-=\{(x,y,0): x_+=y_+=0\}$, $V_+=\{(x,y,0): x_-=y_-=0\}$ and the complex lines $V_k=\{(0,0,c_*): c_l=0 \text{ for } l\neq k\}$ are also invariant.

From equations~(\ref{eq:sys-res-0-full}), it is clear that the origin ($\mathbb{T}_j$) is an equilibrium point of type $(\text{saddle})^2 \times (\text{center})^{n-3}$,
with saddle characteristic exponents $\pm \lambda$ and center characteristic exponents $\pm i\nu_k$.
Indeed the unstable and stable invariant surfaces of the origin, $W^{\mathrm{u}}$ and $W^{\mathrm{s}}$, are tangent at the origin to the real planes $V^{\mathrm{u}}=\{(x,0,0)\}$ and $V^{\mathrm{s}}=\{(0,y,0)\}$.
By Hypothesis~\textbf{H5} or~\textbf{f5}, they contain at least one ingoing heteroclinic orbit arriving at $\mathbb{T}_{j}$ from $\mathbb{T}_{j-1}$ and one outgoing heteroclinic orbit departing from $\mathbb{T}_{j}$ to $\mathbb{T}_{j+1}$. By the straightening (Remarks~\ref{HamStraighten} and~\ref{Straighten}), they can be considered included in straight lines, at least close to the origin.
In other words, in a neighborhood of the origin $(x,y,c_*)=(0,0,0)$, there is an ingoing heteroclinic orbit arriving at $\mathbb{T}_{j}$ from $\mathbb{T}_{j-1}$ contained in the half-line
\begin{equation}
\label{Hin}
H_{\inn}=\{(0,y,0):  y_- >0, y_+=0\}=\{(0,0,y_-(t),0),t\in \mathbb{R}\}
\end{equation}
with $y_-(t)\to 0$ as $t\to \infty$, and an outgoing heteroclinic orbit departing from $\mathbb{T}_{j}$ to $\mathbb{T}_{j+1}$ contained in the half-line
\begin{equation}
\label{Hout}
H_{\out}=\{(x,0,0):  x_-=0, x_+>0\}=\{(0,x_+(t),0,0),t\in \mathbb{R}\}
\end{equation}
with $x_+(t)\to 0$ as $t\to -\infty$.

\subsection{Block diagonal dynamics along the heteroclinic orbits}
\label{subsec:Blockdiagonal}

%\textbf{PZ: the existence of uniform bounds on the constants for $\T$
%should be one of our assumptions if we will look for infinite chain!!!!}

In the $j$-th chart with dynamics~\eqref{eq:sys-res-0-full} we have the incoming and the outgoing heteroclinic orbits (\ref{Hin}-\ref{Hout}).
%Taking into account Hypothesis~\textbf{H3} or~\textbf{f3},
%the real planes $V_-=\{(x,y,c_*): x_+=y_+=0\}$, $V_+=\{(x,y,c_*): x_-=y_-=0\}$
%and the complex hyperplanes $W_k=\{(x,y,c_*): c_k=0\}$ for $k \neq j-1,j,j+1$, are also invariant.

We now introduce what we call \emph{block diagonal dynamics along heteroclinic orbits}, which is the property that the variational equations of the system~\eqref{eq:sys-res-0-full} along these orbits heteroclinics have a \emph{block diagonal} structure:
\begin{eqnarray}
   \dot{z}_- &=&B^{\mathrm{v}}_-(t)z_-,\notag \\
    \dot{z}_+ &=&B^{\mathrm{v}}_+(t)z_+,\label{BlockDiagonal} \\
    \dot c_k&=&i\nu_k c_k,\quad k\neq j-1,j,j+1\notag
\end{eqnarray}
where  $z_-=(x_-,y_-)$, $z_+=(x_+,y_+)$, $\mathrm{v} \in \{\inn,\out\}$,  $B^{\mathrm{v}}_{\pm}(t)$ are some matrices, $t \in (-\infty,\infty)$.
Notice that the blocks are matrices on the complex lines $V_-,V_+,V_k$.
The dependence on $t$ comes from the heteroclinic orbit and
$\inn,\out$ stand for ingoing and outgoing heteroclinic orbit, respectively.

Observe that the block diagonal dynamics along the heteroclinic orbits will hold if along the ingoing and outgoing heteroclinic orbit
the following splitting of the tangent space
\begin{equation}
  T_p \mathcal{M}= p + \left( V_- \oplus V_+ \bigoplus_{k \neq j\pm 1, k\neq j} V_k\right)   \label{eq:decompTM}
\end{equation}
remains \emph{invariant} by the flow on each $j$-chart, $1\leq j \leq n$.

%\subsubsection{Invariant subspaces}
%\label{subsubsec:inv-subspaces}

In the $(z_-,z_+,c_*)$-coordinates of the $j$-chart, where $z_-=(x_-,y_-)$, $z_+=(x_+,y_+)$, we can rewrite  (\ref{eq:sys-res-0-full}) in a more compact way as
\begin{eqnarray*}
  \dot z_-&=&g_{-}(z_-,z_+,c_*),  \notag \\
  \dot z_+&=&g_{+}(z_-,z_+,c_*), \\ \label{eq:z-z+c*}
  \dot c_*&=&g_{*}(z_-,z_+,c_*). \notag
\end{eqnarray*}

As we are assuming that the hyperplanes $W_\pm$ and $W_k$ are invariant, then also $V_-=\{(z_-,0,0)\}$, $V_+=\{(0,z_+,0)\}$ and $V_k$ are invariant.
The ingoing and outgoing heteroclinic orbits $H_{\inn}$ and $H_{\out}$ introduced in~\eqref{Hin} and~\eqref{Hout} take now the form in the
$(z_-,z_+,c_*)$-coordinates of system~\eqref{eq:z-z+c*}:
\begin{align*}
H_{\inn}&=\{(z_-,0,0): z_-=(x_-=0,y_->0)\}\subset V_-,\\%\label{Hinm}\\
H_{\out}&=\{(0, z_+,0): z_+=(x_+>0,y_+=0)\}\subset V_+.%\label{Houtp}
\end{align*}
%The ingoing heteroclinic orbit $H_{\inn}$ to the $j$-mode introduced in~(\ref{Hin})
%is contained in the subspaces  $V_-$, $W_+$ and $W_k$, whereas the outgoing heteroclinic orbit $H_{\out}$ from the $j$-mode introduced in~(\ref{Hout})
%is contained in the subspaces $V_+$, $W_-$ and $W_k$.

%The invariance of $V^j_-$ means that $z_+'$ and $c_k'$ vanish for points $z_- \in \mathbb{R}^2,z_+=0,c_k=0$, hence
%\begin{equation}
%  \frac{\partial g_{+}}{\partial z_-}(z_-,0,0)=0, \quad \frac{\partial g_{k}}{\partial z_-}(z_-,0,0)=0.  \label{eq:inv-sub-var-zera-z-}
%\end{equation}
The invariance of $W_+$  implies that $\dot z_+=0$ if $z_+=0$  for arbitrary $z_-$ and $c_*$, therefore
in particular for any $(z_-,0,0)\in V_-$
\begin{equation}
   \frac{\partial g_{+}}{\partial z_-}(z_-,0,0)=0, \quad \frac{\partial g_{+}}{\partial c_k}(z_-,0,0)=0.  \label{eq:isvar-Wj+}
\end{equation}
The invariance of $W_k$  implies that $\dot c_k=0$ if $c_k=0$ for arbitrary $z_\pm$ and $c_l$ for $l \neq k$, therefore
\begin{equation}
   \frac{\partial g_{k}}{\partial z_-}(z_-,0,0)=0, \quad \frac{\partial g_{k}}{\partial z_+}(z_-,0,0)=0,
     \quad \frac{\partial g_{k}}{\partial c_l}(z_-,0,0)=0, \quad l \neq k  .  \label{eq:isvar-Wjk}
\end{equation}

Now, if we look at the variational equations along the ingoing heteroclinic orbit $H_{\inn}$ which is contained in $V_-$ we obtain
(the Jacobian matrix is evaluated at a point $(z_-,0,0)$ and we put zeros according to (\ref{eq:isvar-Wj+},\ref{eq:isvar-Wjk}))
\begin{equation}
\label{VariationalEq}
  \left[ \begin{array}{c}
               \dot{z}_- \\
               \dot{z}_+ \\
               \dot{c}_*
          \end{array}\right] =
          \left[\begin{array}{ccc}
               \frac{\partial g_{-}}{\partial z_-}, & \frac{\partial g_{-}}{\partial z_+} & \frac{\partial g_{-}}{\partial c_*} \\
               0, & \frac{\partial g_{+}}{\partial z_+} & 0   \\
                0, & 0 & \frac{\partial g_{*}}{\partial c_*}
          \end{array}
          \right]
           \left[ \begin{array}{c}
               z_- \\
               z_+ \\
               c_*
          \end{array}\right]
\end{equation}
where $\frac{\partial g_{*}}{\partial c_*}$ is a diagonal matrix because
from (\ref{eq:isvar-Wjk}) we have $\frac{\partial g_{k}}{\partial c_l}(z_-,0,0)=0$ for $l \neq k$.

Therefore to have a block diagonal Jacobian matrix in the variational equations~\eqref{VariationalEq} or, equivalently, the block diagonal property of system~\eqref{BlockDiagonal}, we have to assume that the off-diagonal entries vanish, that is
\begin{equation}
  \frac{\partial g_{-}}{\partial z_+}(z_-,z_+=0,c_*=0)=0, \quad  \frac{\partial g_{-}}{\partial c_*}(z_-,z_+=0,c_*=0)=0.  \label{eq:add-hete-in}
\end{equation}

Analogously for the outgoing heteroclinic orbit we will need to assume that
\begin{equation}
  \frac{\partial g_{+}}{\partial z_-}(z_-=0,z_+,c_*=0)=0, \quad  \frac{\partial g_{-}}{\partial c_*}(z_-=0,z_+,c_*=0)=0.  \label{eq:add-hete-out}
\end{equation}

Conditions (\ref{eq:add-hete-in}-\ref{eq:add-hete-out}) will be called \emph{block diagonal dynamics along the heteroclinic orbits} and they are essential for the design of our diffusion, consisting on dropping some dimensions, or equivalently directions, when passing close to a periodic orbit. Notice that they are equivalent to the invariance of splitting (\ref{eq:decompTM}) and to the block diagonal structure of equations~(\ref{BlockDiagonal}).

\subsection{Block diagonal dynamics along the heteroclinic orbits holds in the Hamiltonian and general setting}
\label{subsec:HamBlockDiagonal}
It turns out that assumptions  (\ref{eq:add-hete-in}-\ref{eq:add-hete-out}) follow from the Hamiltonian form of our vector field and the Hypothesis~\textbf{H3} about invariance of $W_\pm$ and $W_k$. Indeed the equations of motion~(\ref{eq:z-z+c*}) become
\begin{eqnarray*}
  \dot z_-&=&g_{-}(z_-,z_+,c_*)=J \partial_{z_-}H (z_-,z_+,c_*), \\
  \dot z_+&=&g_{+}(z_-,z_+,c_*)=J \partial_{z_+}H (z_-,z_+,c_*), \\
  \dot c_*&=&g_{*}(z_-,z_+,c_*)=-2i \partial_{\overline{c}_*}H (z_-,z_+,c_*),
  \end{eqnarray*}
where
{\renewcommand{\arraystretch}{0.5}
$J=J_2=\left(\begin{array}{@{}r@{\hskip 2pt}r@{}}0&-1\\-1&0\end{array}\right)$}, and for instance along the ingoing  heteroclinic orbit
\begin{eqnarray*}
\frac{\partial g_{-}}{\partial z_+}(z_-,0,0)&=&J\partial^2_{z_+ z_-} H(z_-,0,0)= \frac{\partial g_{+}}{\partial z_-}(z_-,0,0)=0,\\
\frac{\partial g_{-}}{\partial c_*}(z_-,0,0)&=&J\left(-2i\partial _{c_* z_-}H(z_-,0,0)\right)=J\frac{\partial g_{*}}{\partial z_-}(z_-,0,0)=0.
\end{eqnarray*}
Analogously for the outgoing heteroclinic orbit from the $j$-th mode.

Alternatively, we can take into account the sign symmetry~\textbf{H7}, or~\textbf{f7} for the general setting, to easily prove
conditions~(\ref{eq:add-hete-in}-\ref{eq:add-hete-out}). Notice, however, that the sign symmetry~\textbf{H7} is not necessary in the Hamiltonian case.

%We are now in a position to state our main result in the general setting.
%\begin{theorem}
%\label{thm:toymodel-gen}
% Under Hypotheses \textbf{f1-f6} and dropping conditions (\ref{eq:add-hete-in}-\ref{eq:add-hete-out}),
% \textbf{absence of bad terms is no longer needed}
% plus the absence of some bad terms in system~(\ref{eq:sys-res-0-full}) (see Theorem~\ref{thm:enclo-hyp-dir}),
% for any $n\geq 3$ and for all $\varepsilon>0$ there exists a point $x_1$ close to $\mathbb{T}_1$ whose trajectory is $\varepsilon$-close
% to the  chain of heteroclinic connections $\mathbb{T}_1 \to \mathbb{T}_2 \to \cdots \to \mathbb{T}_n$.
%\end{theorem}

\subsection{The transition map between consecutive charts }
\label{subsubsec:cc-charts}

Consider the $j$-th chart with coordinates $(z_-,z_+,c_*)$, where $z_\pm=(x_\pm,y_\pm)$ and $c_*=(c_k)_{k\neq j-1,j,j+1}$, and the next $(j+1)$-th chart with ``tilde'' coordinates $(\tilde z_-,\tilde z_+,\tilde c_*)$, and denote by $\J=\J_{j \to j+1}$ the transition map between these consecutive charts: $\J(z_-,z_+,c_*)=(\tilde{z}_-,\tilde{z}_+,\tilde{c}_*)$.

The assumption about the existence of the heteroclinic connection departing from $\mathbb{T}_{j}$ to $\mathbb{T}_{j+1}$  means
that the outgoing heteroclinic orbit $H_{\out}$~\eqref{Hout} contained in the segment $(x_+,y_+=0,z_-=0,c_*=0)$ parameterized by $x_+$ is mapped by the transition map to the ingoing heteroclinic orbit $\tilde H_{\inn}$~\eqref{Hin} contained in the segment $(\tilde{y}_-,\tilde{x}_-=0,\tilde{z}_+=0,\tilde{c}_*=0)$ parameterized by $\tilde{y}_-$: $\J(H_{\out})=\tilde H_{\inn}$. Observe that we assumed that the heteroclinics are locally straightened.

Let $p=(x_+,y_+=0,z_-=0,c_*=0)\in H_{\out}$ for some $x_+$. Then its image $\J(p)$ is of the form $(\tilde{z}_+=0,\tilde{y}_-,\tilde{x}_-=0,\tilde{c}_*=0)\in \tilde H_{\inn}$
and $D\J(p)$ is a  map from the tangent space at $p$ which is $V^j_- \oplus V^j_+ \bigoplus_{k \neq j,j\pm 1} V^j_{k} $
to the tangent space at $\J(p)$ which is $V^{j+1}_- \oplus V^{j+1}_+ \bigoplus_{k \neq j+1,j,j+2} V^{j+1}_{k}$. We will use these decompositions
of the tangent space at $p$ and $\J(p)$ to define the blocks of $D\J(p)$. We  assume that  the only non-zero blocks in $D\J(p)$ are (and are isomorphisms)
\begin{itemize}
\item  $\frac{\partial \tilde{c}_k}{\partial c_k}$  ,  $k \leq j-2$ (past modes) and $k \geq j+3$ (future modes),
\item  $\frac{\partial \tilde{c}_{j-1}}{\partial (x_-,y_-)}$, the saddle directions $(x_-,y_-)$ become a past mode,
\item  $\frac{\partial (\tilde{x}_-,\tilde{y}_-)}{\partial (x_+,y_+)}$, the saddle directions $(x_+,y_+)$ are identified with the saddle directions $(\tilde{x}_-,\tilde{y}_-)$,
\item $\frac{\partial (\tilde{x}_+,\tilde{y}_+)}{\partial c_{j+2}}$, the new saddle directions $(\tilde{x}_+,\tilde{y}_+)$ are ``created" from $(j+2)$-th future mode.
\end{itemize}

By computing explicitly the transition map $\J$, we are going to see in the next section that the above decomposition for $D\J(p)$ holds (see~\eqref{eq:BlockStructure}), as well as that the transition map is symplectic and preserves the coordinate invariant hyperplanes.

\subsubsection{The transition between  charts in the global model}
%\label{subsubsec:tran-chart}

In the symplectic chart $z=\Psi_j(z_-,z_+,c_*,r_j,\theta_j)$ centered on the $j$-th torus we have (see ( \ref{eq:j-chart-polar},\ref{eq:j-chart-ck},\ref{eq:rj}))
\[
\begin{array}{l}z_j=r_j e^{i\theta_j},\quad z_k=c_k e^{i\theta_j} \text{ for } k\neq j,\text{ and }\\[3pt]
c_{j-1}=\omega x_- +\bar{\omega} y_-,\quad  c_{j+1}=\omega x_+ +\bar{\omega} y_+,
\end{array}
\]
where $\omega$, $\bar\omega$ are introduced in~(\ref{Conformal}).
In the symplectic chart  \\
$z=\Psi_{j+1}(\tilde z_-,\tilde z_+,\tilde{c_k},r_{j+1},\theta_{j+1})$
centered on the $(j+1)$-th torus, where we introduce tildes in the coordinates to distinguish them form the previous ones, we have
\[
\begin{array}{l}
z_{j+1}=r_{j+1} e^{i\theta_{j+1}},\quad z_k=\tilde{c}_k e^{i\theta_{j+1}} \text{ for } k\neq j+1, \text{ and }\\[5pt]
\tilde{c}_{j}=\omega \tilde{x}_- +\bar{\omega} \tilde{y}_-, \quad \tilde{c}_{j+2}=\omega \tilde{x}_+ +\bar{\omega} \tilde{y}_+.
\end{array}
\]
%\begin{eqnarray}
%  r_j^2&=&1 - \sum_{k \neq j} |b_k|^2= 1 - \sum_{k \neq j} |c_k|^2=1 - \sum_{k \neq j} r_k^2, \\
%  b_j&=&r_j e^{i \theta_j}, \quad b_k=c_k e^{i \theta_j}, k \neq j  \\
%  c_{j+1}&=&\omega y_+ + \omega^2 x_+, \quad  c_{j-1}=\omega y_- + \omega^2 x_-.
%\end{eqnarray}

We look for the symplectic relation
\[
(\tilde z_-,\tilde z_+,\tilde{c}_*,r_{j+1},\theta_{j+1})=\left(\Psi^{-1}_{j+1}\circ \Psi_j\right)(z_-,z_+,c_*,r_j,\theta_j)
\]
between coordinates in the $j$-th and $(j+1)$-th charts.
%First we will prove the following simple lemma for $\omega$, $\bar\omega$ introduced in~(\ref{Conformal}).
\begin{lemma}
If $c=\omega x + \bar{\omega} y$ for $\omega$ introduced in~(\ref{Conformal}) and real $x,y$, then
\begin{equation*}
 |c|^2=x^2-2a_1 x y +y^2, \quad  \frac{1}{c}= \frac{1}{|c|^2}\left(\bar{\omega}x + \omega y\right), \text{ where }  a_1 = - 2 \Re \omega^2 .
\end{equation*}
\end{lemma}
\textbf{Proof:}
As $\bar{c}=\bar{\omega} x + \omega y$, the assertion follows from
$\displaystyle \frac{1}{c}=\frac{\overline{c}}{|c|^2}$.
\qed

From $z_{j+1}=c_{j+1} e^{i\theta_j}=r_{j+1} e^{i\theta_{j+1}}$ we get
\begin{eqnarray*}
  r_{j+1}&=&|c_{j+1}|=\sqrt{x_+^2-2a_1 x_+ y_+ +y_+^2}, \\
   e^{i \theta_{j+1}}&=& \frac{c_{j+1}}{|c_{j+1}|} e^{i \theta_j},  \\
   e^{i (\theta_j - \theta_{j+1})}&=&\frac{|c_{j+1}|}{c_{j+1}}=\frac{\bar{\omega} x_+ + \omega y_+}{|c_{j+1}|}
   =\frac{\bar{\omega} x_+ + \omega y_+}{\sqrt{x_+^2-2a_1 x_+ y_+ +y_+^2}},
\end{eqnarray*}
and therefore
\begin{equation}
\label{eq:transitionj}
\begin{split}
  r_j &= \sqrt{1 - (x_+^2 -2a_1 x_+ y_+ + y_+^2) - (x_-^2 -2a_1 x_- y_- + y_-^2) - \sum_{k \neq j,j\pm 1} c_k \bar{c}_k} \\
 \tilde{c}_k &= e^{i(\theta_j - \theta_{j+1})}c_k = \frac{|c_{j+1}|}{c_{j+1}}c_k , \quad  k < j-1 \  \mbox{or} \ k > j+2, \\
 \tilde{c}_{j-1} &= \frac{|c_{j+1}|}{c_{j+1}}c_{j-1} = \frac{|c_{j+1}|}{c_{j+1}}\left(\omega x_- +\bar{\omega} y_-\right) , \\
 \tilde{c}_j&=r_j e^{i(\theta_j - \theta_{j+1})} = r_j \frac{|c_{j+1}|}{c_{j+1}} =
  \frac{r_j}{|c_{j+1}|} \left(\bar{\omega}x_+ + \omega y_+\right), \\
 \tilde{x}_-&= \frac{r_j y_+}{|c_{j+1}|}, \qquad \tilde{y}_-= \frac{r_j x_+}{|c_{j+1}|},\\
 \omega\tilde x_+ +\bar\omega\tilde y_+ &=\tilde{c}_{j+2}= \frac{|c_{j+1}|}{c_{j+1}}c_{j+2}.
\end{split}
\end{equation}

%\begin{eqnarray*}
% \tilde{c}_k &=& e^{i(\theta_j - \theta_{j+1})}c_k = \frac{|c_{j+1}|}{c_{j+1}}c_k , \quad  k < j-1 \  \mbox{or} \ k > j+2, \\
% \tilde{c}_{j-1} &=& \frac{|c_{j+1}|}{c_{j+1}}c_{j-1} = \frac{|c_{j+1}|}{c_{j+1}}\left(\omega x_- +\bar{\omega} y_-\right) , \\
% \tilde{c}_j&=&r_j e^{i(\theta_j - \theta_{j+1})} = r_j \frac{|c_{j+1}|}{c_{j+1}} =
%  \frac{r_j}{|c_{j+1}|} \left(\bar{\omega}x_+ + \omega y_+\right), \\
% \tilde{x}_-&=& \frac{r_j y_+}{|c_{j+1}|}, \qquad \tilde{y}_-= \frac{r_j x_+}{|c_{j+1}|},\\
% \omega\tilde x_+ +\bar\omega\tilde y_+ &=&\tilde{c}_{j+2}= \frac{|c_{j+1}|}{c_{j+1}}c_{j+2}.
%\end{eqnarray*}
Observe that for the above coordinate change the invariant hyperplanes $W_l$ are preserved, as well as the block diagonal dynamics along the heteroclinic connections,
because we have
\begin{equation}
\setlength{\arraycolsep}{1.4pt}\renewcommand{\arraystretch}{1.4}
\label{eq:BlockStructure}
\begin{array}{rcl}
  \tilde{c}_k &=& g_k(z_+) c_k,  \qquad  k < j-1 \  \mbox{or} \ k > j+2, \\
 \tilde{c}_{j-1} &=& g_{j-1}(z_-,z_+) z_-, \\
  \tilde{z}_- &=& g_-(r_j(z_+,z_-,c_*), z_+) z_+, \\
  \tilde{z}_+&=& g_+(z_+) c_{j+2},
\end{array}
\end{equation}
where $g_i$ are smooth functions. Moreover, this transition map is symplectic in the Hamiltonian case.

\subsection{Polynomial normal form for the saddle variables}
\label{sub:polynomialNFc}
We can take one more step in the simplification of the system close to the partially saddle circle $\mathbb{T}_j$,  consisting of writing system~\eqref{eq:sys-res-0-full} in the $j$-chart in resonant polynomial normal form with respect to the saddle variables in appropriate coordinates. This normal form step was also carried out in~\cite{GK,GHP}, and will make the calculations easier in order to find adequate estimates for the transition of the trajectories of the TMS near $\mathbb{T}_j$.

In the $j$-th chart with coordinates $z=(z_-,z_+,c_*)$, with $z_\pm=(x_\pm,y_\pm)$ and $c_*=(c_\ell)_{\ell\neq j-1,j,j+1}$,
for any index $m=(m_{x_-},m_{y_-},m_{x_+},m_{y_+},m_{c}) \in \mathbb{N}^5$ we introduce the notation for a monomial as
\[
z^m=x_-^{m_{x_-}} y_-^{m_{y-}} x_+^{m_{x_+}} y_+^{m_{y+}} c_*^{m_c},
\]
where the symbol $c_*^{m_c}$ denotes any monomial $c_*^{m_*}\cdot\overline{c}_*^{\overline{m}_*}$ of \emph{center degree}
$\abs{m_*}+\abs{\overline{m}_*}=m_c$ in the center variables $c_\ell,\overline{c}_k$. Analogously we will call $m_s:=m_{x_-}+m_{y_-}+m_{x_+}+m_{y_+}$ the \emph{saddle degree}
so that the total degree of the monomial satisfies $\abs{m}=m_s+m_c$.
Notice that the monomial $z^m$ is completely determined in the saddle variables $z_\pm$, but not fully specified in the center variable $c_*$, because its concrete expression in the variable $c_*$ will not be necessary.

%We will identify the term $z^J$ with an element in $\mathbb{N}^5$.

Taking into account the linear terms of our model~\eqref{eq:sys-res-0-full} in the $j$-th chart we say that $z^m$ is a \emph{resonant monomial for the saddle variables} if
for any saddle variable $v \in \{x_-,y_-,x_+,y_+\}$
\begin{equation}
\label{eq:resonant-monomial}
  \lambda_v = m_{x_-}-m_{y_-}+ m_{x_+} -m_{y_+},
\end{equation}
where $\lambda_{x_\pm}=1$ and $\lambda_{y_\pm}=-1$ in~\eqref{eq:sys-res-0-full}.

The following lemma follows from results established in \cite{DZ} (see Section \emph{Sign symmetry})
\begin{lemma}
\label{lem:PolyNormForm}
For any $k\geq 1$ if system~\eqref{eq:sys-res-0-full} is $C^r$ with $r$ sufficiently
 large, then there exists a $C^k$ change of variables in a neighborhood of the origin,
 transforming system~\eqref{eq:sys-res-0-full} to the system
 \begin{equation}
\label{eq:NF}
\begin{split}
  \dot{x}_-&= \lambda x_- +  N_{x_-}(z), \\
  \dot{y}_-&= -\lambda y_- +  N_{y_-}(z),  \\%\label{eq:NFy-}   \\
  \dot{x}_+&=  \lambda x_+ + N_{x_+}(z), \\ % \label{eq:NFx+} \\
  \dot{y}_+&= -\lambda y_+ + N_{y_+}(z), \\
 \dot{c}_\ell &=  i\nu_l c_\ell +O_2(z), \quad \ell \neq j, j \pm 1, % \label{eq:NFcl}
\end{split}
\end{equation}
where for any saddle variable $v \in \{x_-,y_-,x_+,y_+\}$ we have
\[
  N_v(z)=\sum_{m \in M_{1,v}}  g_{v,m}(c_*)z^m  + \sum_{m \in M_{2,v}} g_{v,m}(z)z^m,  %\label{eq:fNF-full-enclo}
\]
where   $g_{v,m}$ are continuous functions,
$M_{1,v},M_{2,v}$ are \emph{finite} sets of indices,
and any $z^m$ is a resonant monomial for the saddle variables~\eqref{eq:resonant-monomial},
satisfying on the one hand $m_s:=m_{x_-}+m_{y_- }+m_{x_+ } +m_{y_+}\geq 3$ and $m_c =0$ if $m=\left(m_{x_-},m_{y_- },m_{x_+},m_{y_+},m_c\right) \in M_{ 1, v }$,
and on the other hand $m_s=1$ and $m_c \geq 3$ if $m\in M_{2,v }$.

%\noindent\textbf{A: I wrote $m_s\geq 3$ because there are no resonant terms of order 2.}

%\noindent\textbf{A: I wonder if we can put $=O_3(z)$ in $\dot{c}_\ell$ as in~\cite{CK}  Check resonant conditions}

This change of coordinates preserves the invariant subspaces $W_l$'s as well as the sign symmetry~\textbf{f7} and therefore
the block diagonal structure along $V_\pm$'s.
\end{lemma}

The main part of the above theorem (without invariant subspaces and block diagonal structure) has been established in \cite{BK95,BLW,BK96}. Part of the technique developed in \cite{BLW} was used in the proof in \cite{DZ}
%, but apparently our estimate for $r$, $k$  are not as good as given there. It appears that  this difference  is a consequence of the need to handle the block %diagonal structure.
where it is proved that $r=11$ is enough to have $k=2$, which is minimum smoothness of coordinate change required for our construction
of h-sets to shadow the non-transversal heteroclinic chain for TMS.

\begin{rem}
%\label{rem:PNF-straightMan}
  In the polynomial normal form variables the local stable and unstable manifolds are straight as well as the ingoing and outgoing heteroclinic orbits.
  \end{rem}

%The main part of the above theorem (without invariant subspaces and block diagonal structure) has been established in \cite{BK95,BLW,BK96}. Part of the technique developed in %\cite{BLW} was used in the proof
%in \cite{DZ}, but apparently our estimate for $r$, $k$  are not as good as given there. It appears that  this difference  is a consequence of the need to handle the block %diagonal structure.

%% file: covrel.tex
\section{h-sets, covering relations and dropping dimensions}
\label{sec:covrel}

The goal of this section is to recall from \cite{ZGi} the notions of h-sets and covering relations, and state the theorem
about the existence of points realizing the chain of covering relations. This will be the main technical tool to prove the existence of the orbits shadowing the  heteroclinic chain in the following sections.

\subsection{h-sets and covering relations}
%\label{subsec:covrel}

By $B_n$ we will denote a unit ball in $\mathbb{R}^n$ and by $B_n(c,r)$ we denote a ball in $\mathbb{R}^n$ with the center at $c$ and radius $r$.
The norm used in definition of $B_n$ or $B_n(c,r)$ is either known from the context or is arbitrary, but fixed.

\begin{definition} \cite[Definition 1]{ZGi}
\label{def:covrel} An $h$-set $N$ is a quadruple
$\left(\abs{N},u(N),s(N),c_N\right)$ such that
\begin{itemize}
 \item $|N|$ is a compact subset of ${\mathbb R}^n$
 \item $u(N),s(N) \in \{0,1,2,\dots,n\}$ are such that $u(N)+s(N)=n$
 \item $c_N:{\mathbb R}^n \to
   {\mathbb R}^n={\mathbb R}^{u(N)} \times {\mathbb R}^{s(N)}$ is a
   homeomorphism such that
      \begin{displaymath}
        c_N(|N|)=\overline{B_{u(N)}} \times
        \overline{B_{s(N)}}.
      \end{displaymath}
\end{itemize}
We set
\begin{eqnarray*}
   \dim(N) &:=& n,\\
   N_c&:=&\overline{B_{u(N)}} \times \overline{B_{s(N)}}, \\
   N_c^-&:=&\partial B_{u(N)} \times \overline{B_{s(N)}}, \\
   N_c^+&:=&\overline{B_{u(N)}} \times \partial B_{s(N)}, \\
   N^-&:=&c_N^{-1}(N_c^-) , \quad N^+=c_N^{-1}(N_c^+).
\end{eqnarray*}
\end{definition}

Hence a $h$-set $N$ is a product of two closed balls in some
coordinate system. The numbers $u(N)$ and $s(N)$ are called the
\emph{exit} and \emph{entry} dimensions, respectively.
The subscript $c$ refers to the new coordinates given by
the homeomorphism $c_N$. Observe that if $u(N)=0$, then
$N^-=\emptyset$ and if $s(N)=0$, then $N^+=\emptyset$. In the
sequel to make the notation less cumbersome we will often drop the
bars in the symbol $|N|$ and we will use $N$ to denote both the
h-sets and its support.

We will call $N^-$ \emph{the exit set of N} and $N^+$
\emph{the entry set of $N$}. These names are motivated by the
Conley index theory \cite{C,MM} and the role that these sets will play in the
context of covering relations.

\begin{definition}\cite[Definition 6]{ZGi}
\label{def:covw} Assume that $N,M$ are $h$-sets, such that
$u(N)=u(M)=u$ and $s(N)=s(M)=s$. Let $f:N \to {\mathbb R}^n$ be a
continuous map and $f_c:= c_M \circ f \circ c_N^{-1}: N_c \to
{\mathbb R}^u \times {\mathbb R}^s$. Let $w$ be a nonzero integer.
We say that \emph{$N$ $f$-covers $M$ with degree $w$}, in symbols
\begin{displaymath}
  N\cover{f,w} M,
\end{displaymath}
iff the following conditions are satisfied
\begin{description}
\item[1.] There exists a continuous homotopy $h:[0,1]\times N_c \to {\mathbb R}^u \times {\mathbb R}^s$,
   such that the following conditions are met:
   \begin{eqnarray}
      h_0&=&f_c,\notag  \\ %\label{eq:hc1}
      h([0,1],N_c^-) \cap M_c &=& \emptyset ,  \label{eq:hc2} \\
      h([0,1],N_c) \cap M_c^+ &=& \emptyset .\label{eq:hc3}
   \end{eqnarray}
\item[2.] If $u >0$, then there exists a  map $A:{\mathbb R}^u \to {\mathbb
R}^u$ such that
   \begin{eqnarray}
    h_1(p,q)&=&(A(p),0), \mbox{ for $p \in \overline{B_u}(0,1)$ and $q \in
    \overline{B_s}(0,1)$,}\label{eq:hc4}\\
      A(\partial B_u(0,1)) &\subset & {\mathbb R}^u \setminus
      \overline{B_u}(0,1).  \label{eq:mapaway}
   \end{eqnarray}
  Moreover, we require that
\begin{equation}
  \deg(A,\overline {B_u}(0,1),0)=w, \label{eq:deg-A}
\end{equation}
\end{description}

We will call condition \eqref{eq:hc2} \emph{the exit condition}
and condition \eqref{eq:hc3} \emph{the entry
condition}.

\end{definition}
Note that in the case $u=0$, if $N \cover{f,w} M$, then $f(N)
\subset \inter M$ and  $w=1$.
%\textbf{PZ: Think  whether to keep $w$ any or just set the definition for linear $A$}
Also note that in the above definition $s(N)$ and $s(M)$ could be
different, see \cite[Def. 2.2]{W2}.

%\begin{rem}
%Observe, that since for any norm in $\mathbb{R}^n$ the closed unit
%ball is homeomorphic to $[-1,1]^n$, therefore for h-sets and
%covering relations we may use different norms in different
%contexts.
%\end{rem}

\begin{rem}If the map $A$  in condition 2 of Def.~\ref{def:covw} is a linear
map, then  condition~(\ref{eq:mapaway}) implies that
\begin{displaymath}
  \deg(A,\overline {B_u}(0,1),0)=\pm 1.
\end{displaymath}
Hence condition (\ref{eq:deg-A}) is fulfilled with $w =\pm 1$.
  In fact, this is the most common situation in the applications of
  covering relations.
\end{rem}

Often we will not need the value of $w$ in the symbol $N \cover{f,w} M$ and then we omit it by simply writing $N \cover{f} M$.
If $f$ is known from the context, we can even drop it and just write $N \cover{} M$.

\subsection{h-sets with a product structure }

In view of our application to the toy model system we define a h-set with a product structure.

\begin{definition}
\label{def:hset-product}
Assume that for our ambient space $\mathbb{R}^n$ we have a decomposition into vector spaces
\begin{eqnarray*}
 \mathbb{R}^n &=& X \oplus Y , \\
  X &=& X_1 \oplus X_2 \oplus \dots \oplus X_{u'}, \\
  Y &=& Y_1 \oplus Y_2 \oplus \dots \oplus Y_{s'},
\end{eqnarray*}
where $X_i$ and $Y_j$ are equipped with some norms. According to this decomposition we will represent points $p \in \mathbb{R}^n$
as $p=(x,y)$, where $x=(x_1,\dots,x_{u'})\in X$ and $y = (y_1,\dots,y_{s'}) \in Y$.

We define a norm on $X \oplus Y$
by
\[
  \|(x,y)\|=\max \left( \max_{j=1,\dots,u'} |x_j| , \max_{i=1,\dots,s'}  |y_i| \right),
\]
which induces norms on $X$ and $Y$ which are used below.

Given collections of positive numbers $\Gamma_i$, for $i=1,\dots,s'$, and $R_j$, for $j=1,\dots,u'$, as well as a center point $(x^c,y^c) \in X \oplus Y$, we define a \emph{h-set $N$ with the product structure}  as follows:
\begin{itemize}
\item $|N|$, the \emph{support} of $N$, which is given by
\begin{eqnarray*}
   |N|&=&\left\{(x,y) \in X \oplus Y, \, |x_j- x_j^c| \leq R_j, \, j=1,\dots,u' |y_i-y_i^c| \leq \Gamma_i, \, i=1,\dots,s' \right\},
\end{eqnarray*}
\item $u(N)=\sum_{j=1}^{u'}  \dim X_i $, $s(N)=\sum_{i=1}^{s'}  \dim Y_i $,
\item $c_N: \mathbb{R}^n \to X \oplus Y$ is defined by
   \begin{eqnarray*}
       \pi_{X_j}(x,y)=R_j^{-1}\left(x_j-x_j^c\right), j=1,\dots,u',\\
       \pi_{Y_i}(x,y)=\Gamma_i^{-1}\left(y_i-y_i^c\right), \quad i=1,\dots,s'.
    \end{eqnarray*}
\end{itemize}
\end{definition}
Note that with this definition $X$ represents the exit directions while $Y$ represents the exit directions.

We now state a theorem that will be used later to establish covering relations between two sets h with a product structure. $B_{X_j}(p,r)$ and $B_{Y_i}(p,r)$ denote open balls in $X_j$ or $Y_i$ with center $p$ and radius~$r$.
\begin{theorem}
\label{thm:cv-prod-struct}
 Assume that $N$ and $M$ are two h-sets with the product structure (for $M$ we denote $\{X_j\}$ and $\{Y_i\}$ and other constants with a tilde).

 Assume that $\tilde{u}'=u'$ and $\dim X_j=\dim \tilde{X}_j$ for $j=1,\dots,u'$.

 Assume that we have a continuous map $f:|N| \to \tilde{X} \oplus \tilde{Y}$ and a continuous homotopy $H:[0,1]\times |N| \to \tilde{X}\oplus \tilde{Y} $ such that the inequalities
 \begin{eqnarray}
   \|\pi_{\tilde{Y}_i}H_t(N) - \tilde{y}^c_i\| &<& \tilde{\Gamma}_i, \quad i=1,\dots,\tilde{s}',   \label{eq:entry-Yi} \\
   \|\pi_{\tilde{X}_j}H_t(x,y) - \tilde{x}^c_j \| &>& \tilde{R}_j, \, \mbox{if $(x,y) \in N^-$, $\|x_j\|=R_j$} \quad j=1,\dots,\tilde{u},  \label{eq:exit-Xj}
 \end{eqnarray}
 hold for all $t \in [0,1]$, as well as
 \[
 H_0=f,
 \]
 and $H_1$ is affine diagonal on exit directions in the sense that
  \begin{eqnarray}
   \pi_{\tilde{X}_j}H_1(x,y)=\mathcal{A}_j(x_j):=A_j(x_j - x_j^c) + a_j,  \label{eq:H1-aff-diag}
 \end{eqnarray}
is satisfied for each $j=1,\dots,u'$, where $A_j: X_j \to \tilde{X}_j$ is a linear isomorphism, $a_j \in \tilde{X}_j$,
and one of the following conditions
 \begin{eqnarray}
   \|a_j - \tilde{x}^c_j\| &<& \tilde{R}_j,  \label{eq:aj-close} \\
   \tilde{x}_j^c &\in& \mathcal{A}_j(B_{X_i}(x_j^c,R_j)).  \label{eq:xjc-in-image}
 \end{eqnarray}
is fulfilled.

 Then
 \[
   N \cover{f,w} M,
 \]
 where $w=\sgn \Pi_{j=1}^{u'} \det A_j$.
\end{theorem}
\textbf{Proof:}
Let $h_t= c_M \circ H_t \circ c_N^{-1}$, for $t \in [0,1]$, the homotopy $H$ expressed in internal coordinates of $N$ on input
and of $M$ on output. Using the formulas for $c_N$ and $c_M$ (see Def.~\ref{def:hset-product}),  conditions (\ref{eq:exit-Xj}) and (\ref{eq:entry-Yi}) become
 \begin{eqnarray}
   \|\pi_{\tilde{Y}_i}h_t(N_c)\| &<& 1, \quad i=1,\dots,\tilde{s}' ,  \label{eq:entry-Yic} \\
   \|\pi_{\tilde{X}_j}h_t(x,y) \| &>& 1, \, \mbox{if $(x,y) \in N^-_c$, $\|x_j\|=1$} \quad j=1,\dots,\tilde{u}'. \label{eq:exit-Xjc}
 \end{eqnarray}
It is easy to see that (\ref{eq:entry-Yic}) implies the entry condition (\ref{eq:hc3}), while the exit condition (\ref{eq:hc2}) follows from condition (\ref{eq:exit-Xjc}).

From (\ref{eq:H1-aff-diag}) it follows that  $h_1$ is affine diagonal on exit directions, so for each $j=1,\dots,u'$
 \[
   \pi_{\tilde{X}_j}h_1(x,y)= \mathcal{L}x_j =  L_j x_j + b_j = \frac{R_j}{\tilde{R}_j}A_jx_j + \tilde{R}_j^{-1}(a_j-\tilde{x}^c_j)
 \]
holds,  where $L_j: X_j \to \tilde{X}_j$ is a linear isomorphism  and one of the following conditions (compare (\ref{eq:aj-close},\ref{eq:xjc-in-image}))
  \begin{eqnarray}
   \|b_j \| &<& 1,  \label{eq:bj-close} \\
   0 &\in& \mathcal{L}_j(B_{X_i}(0,1))  \label{eq:0-in-image}
 \end{eqnarray}
 is satisfied. As $h$ is not yet a good enough homotopy, we need to modify it so that $h_1$ satisfies conditions from point 2 in Def.~\ref{def:covrel}.

First we achieve that $\pi_Y h_1(x,y)=0$, which is a part of condition (\ref{eq:hc4}).
For this we define $C: [0,1] \times (\tilde{X} \oplus \tilde{Y}) \to \tilde{X} \oplus \tilde{Y}$  by
\[
  C(t,x,y)=(x,(1-t)y),
\]
i.e., $C$ is a deformation retraction on $0$ in the $y$-direction.
Now the new homotopy
\[
  g(t,x,y)= C(t,h(t,x,y))
\]
obviously satisfies exit condition (\ref{eq:hc3}) and entry condition (\ref{eq:hc2}), and moreover, we have
\begin{eqnarray*}
  \pi_{\tilde{Y}} g_1(x,y)&=&0, \\
    \pi_{\tilde{X}_j}g_1(x,y)&=& \mathcal{L}_j (x_j).
\end{eqnarray*}

Our next step is to construct $\tilde{h}$, a homotopy connecting $g_1$ with a linear map of the desired form, i.e., satisfying conditions (\ref{eq:hc4}) and (\ref{eq:mapaway}) by setting
\begin{eqnarray*}
   \pi_{\tilde{Y}} \tilde{h}(t,x,y)&=&0, \\
    \pi_{\tilde{X}_j}\tilde{h}(t,x,y)&=& L_j x_j + (1-t)b_j.
\end{eqnarray*}
It is obvious that homotopy $\tilde{h}$ satisfies entry condition  (\ref{eq:hc3}) and that $\tilde{h}_1$ is in the form required by (\ref{eq:hc4}).
We will show that the exit conditions (\ref{eq:hc2}) and (\ref{eq:mapaway}) are also satisfied.
For this it is enough to show that
\begin{equation}
   \|L_j x_j + (1-t)b_j\| > 1, \quad \mbox{for all $\|x_j\| =1$ and $t\in [0,1]$}. \label{eq:Ljhom-map-away}
\end{equation}

Let us fix $j$ and denote
\[
  K =  \overline{B}_{\tilde{X}_j}(0,1), \qquad Z=  L_j  B_{X_j}(0,1).
\]

Since $L_j$ is a linear isomorphism, the set $b_j + \partial Z$ separates $\tilde{X_j}$ into two connected components, one bounded
$b_j + Z$ and the other unbounded $X_j \setminus (b_j +\bar{Z})$. Moreover, from (\ref{eq:exit-Xjc}) for $h_1$ we know that
\[
  b_j + \partial Z \subset \tilde{X}_j \setminus K,
\]
and from the fact that one of  conditions (\ref{eq:bj-close}) and (\ref{eq:0-in-image}) is satisfied we have
\[
  b_j \in K \qquad  \mbox{or} \qquad 0 \in b_j + Z.
\]
In both cases we have that  $(b_j + Z) \cap K \neq \emptyset $.
Therefore since $K$ is connected we obtain
\begin{equation}
  K \subset b_j + Z.  \label{eq:K-subset-bjA}
\end{equation}

To prove (\ref{eq:Ljhom-map-away}) we need to show that for each $\alpha \in [0,1]$
\begin{equation}
   K \subset \alpha b_j + Z,  \label{eq:K-t-subset-bjA}
\end{equation}
is fulfilled, because this implies that for all $t \in [0,1]$ we have
\[
   \overline{B}_{\tilde{X}_j}(0,1) \subset (1-t)b_j +  L_j  B_{X_j}(0,1),
\]
which gives (\ref{eq:Ljhom-map-away}).

In order to prove (\ref{eq:K-t-subset-bjA}) we observe  that $Z$ and $K$ are convex and
\[
  K=-K, \quad Z = - Z.
\]
Hence from  (\ref{eq:K-subset-bjA}) it follows that
\begin{eqnarray}
  K=-K &\subset& -b_j - Z = -b_j + Z, \notag \\
  K &\subset& (-b_j + Z) \cap (b_j + Z). \label{eq:L-sub-intersect}
\end{eqnarray}

Now let us take $x\in K$. Then by (\ref{eq:L-sub-intersect}) there exists $z_1,z_2 \in Z$ such that
\[
  x=-b_j + z_1=b_j + z_2.
\]

For any $\theta \in [0,1]$ we have that
\[
  x=\theta(-b_j + z_1) + (1-\theta)(b_j + z_2)= (1-2\theta)b_j + \theta  z_1 + (1-\theta) z_2
\]
holds. Let us fix $\alpha \in [0,1]$. Then for $\theta=\frac{1-\alpha}{2}$ and $z_3=\theta  z_1 + (1-\theta) z_2 \in Z$ we have
that $x \in \alpha b_j + Z$. This finishes the proof of (\ref{eq:K-t-subset-bjA}).

\qed

\subsection{The operation of dropping some exit dimensions}

In this subsection we define what we mean by dropping some exit dimensions in a $h$-set.
\begin{definition}
\label{def:dropping-exit-dim}
Assume that we have a decomposition $\mathbb{R}^n=\mathbb{R}^{u}  \oplus \mathbb{R}^t \oplus   \mathbb{R}^{s}$ and the norm for $(x_1,x_2,x_3) \in \mathbb{R}^{u}  \oplus \mathbb{R}^t \oplus  \mathbb{R}^{s}$ is $\|(x_1,x_2,x_3)\|=\max (\|x_1\|,\|x_2\|,\|x_3\|)$.

Assume that $N$ is an h-set, with $u(N)=u+t$ and $s(N)=s$.  In view of the norm on $\mathbb{R}^n$ we have
\[
  c_{N}(|N|)=\left(\overline{B}_{u} \oplus \overline{B}_{t} \right) \oplus \overline{B}_{s},
\]
where the parentheses enclose the exit directions.

Let us denote by $V$ the subspace $\{0\} \times \mathbb{R}^t \times \{0\}$ and let $\pi_V$ be the projection onto $V$ (i.e., on the $\mathbb{R}^t$ component).

%\textbf{PZ: replace "contracted h-set" by something else  like: cross\-section ?  }
For  $v \in B_t$  we define a contracted h-set $R_{V,v} N$
as a subset in $N$ as follows.

Let
\[
  R_{V,v}N_c = \{z \in N_c, \ \pi_V z=v  \}
\]
then we set
\begin{eqnarray*}
  u(R_{V,v}N)&=&u, \quad s(R_{V,v}N)=s, \\
 |R_{V,v}(N)|&=&c_N^{-1}( R_{V,v}N_c) , \\
 \dim(R_{V,v}N) &=& n-t,\\
   R_{V,v}(N) ^\pm_c&=&N_c^\pm \cap R_{V,v}(N)_c \\
   R_{V,v}(N)^\pm&=&c_N^{-1}\left(  R_{V,v}(N)^\pm \right).
\end{eqnarray*}
\end{definition}
Roughly speaking,  $R_{V,v}N$ is obtained from $N$ by contracting in some exit directions.
An important feature is that  on $R_{V,v}N$ we use coordinates inherited from $N$. This is a technical issue, which facilitates the proof of our main shadowing result in the next section.

The definition given here differs significantly from an
analogous definition in  \cite{DSZ} (Def. 6).
In \cite{DSZ} the operation of dropping exit dimensions was simply a relabeling of some exit direction as an entry one, which then led to some artificial (but negligible) difficulties in our constructions. Although both definitions can be used in our construction of  chain of covering relations, the one given here appears to have more geometric appeal and avoids the aforementioned artificial complications.

We can also define covering relation for a contracted h-set.
\begin{definition}
\label{def:covw-contr}
Assume that $N,M$ are $h$-sets, such that $u(N)=u+t$, $u(M)=u$ and $s(N)=s-t$, $s(M)=s$, and assume that $t,V,v$ are such that $R_{V,v}(N)$ makes sense (see Def.~\ref{def:dropping-exit-dim}) Let $f:R_{V,v}N \to {\mathbb R}^n$ be a
continuous map and $f_c= c_M \circ f \circ c_N^{-1}: R_{V,v}(N)_c \to
{\mathbb R}^u \times {\mathbb R}^s$. Let $w$ be a nonzero integer.

We say that \emph{$ R_{V,v}N$ $f$-covers $M$ with degree $w$}, in symbols
\begin{displaymath}
   R_{V,v}N \cover{f,w} M
\end{displaymath}
iff the following conditions
are satisfied
\begin{description}
\item[1.] There exists a continuous homotopy $h:[0,1]\times R_{V,v}(N)_c \to {\mathbb R}^u \times {\mathbb R}^s$,
   such that the following conditions hold true
   \begin{eqnarray*}
      h_0&=&f_c,  \\%\label{eq:hcch1}
      h([0,1],R_{V,v}N_c^-) \cap M_c &=& \emptyset , \\% \label{eq:hcch2}
      h([0,1],R_{V,v}N_c) \cap M_c^+ &=& \emptyset .%\label{eq:hcch3}
   \end{eqnarray*}
\item[2.] If $u >0$, then there exists a  map $A:{\mathbb R}^u \to {\mathbb
R}^u$ such that
   \begin{eqnarray*}
    h_1(p,v,q)&=&(A(p),0), \mbox{ for $p \in \overline{B_u}(0,1)$ and $q \in
    \overline{B_{s-t}}(0,1)$,}\\%\label{eq:hcch4}
      A(\partial B_u(0,1)) &\subset & {\mathbb R}^u \setminus
      \overline{B_u}(0,1). % \label{eq:mapawaych}
   \end{eqnarray*}
  Moreover, we require that
\[
  \deg(A,\overline {B_u}(0,1),0)=w. %\label{eq:deg-Ach}
\]
\end{description}
\end{definition}

\begin{figure}[hbt!]
\centering
\includegraphics[width=0.45\textwidth]{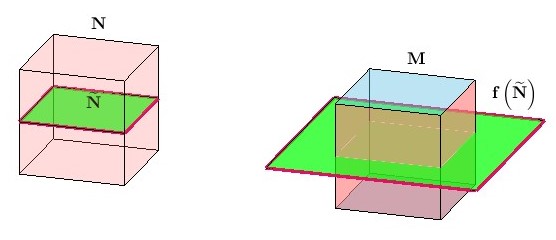}\hspace{0.09\textwidth}
\includegraphics[width=0.45\textwidth]{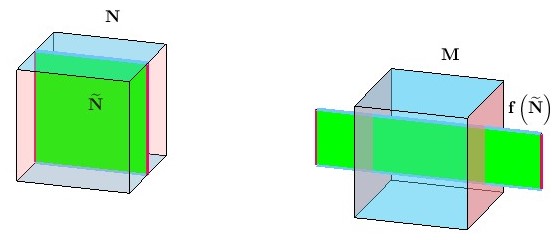}
\caption{Illustration of the covering relation $R_{V,v}N=\widetilde{N} \cover{f} M$, where $N$ and $M$ are 3-dimensional h-sets.
On the left picture $u(N)=3$ and $s(N)=0$, the exit set (in red) $N^-$ of $N$ is its whole boundary $\partial N$, and $\widetilde{N}$ is the horizontal rectangle $z=0$; on the other hand, $u(M)=2$ and $s(M)=1$, the entry set (in blue) $M^+$ consists of the horizontal faces of $\partial M$ whereas the exit set (in red) $M^-$ of the lateral faces.
On the right picture $u(N)=2$, $s(N)=1$, the exit set (in red) $N^-$ consists of the lateral faces of $\partial N$ and the entry set (in blue) $N^-$ of the horizontal faces; on the other hand, $u(M)=1$, $s(M)=2$, the entry set (in blue) $M^+$ consists of the top, down, front and back faces whereas the exit set (in red) $M^-$ of the left and right faces.}

%On the left, $u(N)=u(M)=3$, $s(N)=s(M)=0$, the exit sets (in red) of $N^+$ and $M^+$ are the whole boundaries $\partial N$ and $\partial M$, respectively, and $\widetilde{N}$ is the horizontal rectangle $z=0$.
%On the right, $u(N)=u(M)=2$, $s(N)=s(M)=1$, the exit sets (in red) $N^+$ and $M^+$ are the lateral boundaries of $\partial N$ and $\partial M$, the entry sets (in blue) $N^-$ and $M^-$ are the horizontal boundaries, and $\widetilde{N}$ is the vertical rectangle $x=0$.}
\label{fig:Cubes}
\end{figure}

For an illustration of covering relations with a contracted h-set see Figure~\ref{fig:Cubes}. Note that the above definition~\ref{def:covw-contr} almost coincides with the definition~\ref{def:covrel} of a covering relation between two h-sets.
We just have a lower-dimensional set on the left of the relation $\cover{f}$, but which has the same number of exit directions. Despite this,
a contracted h-set cannot be on the right side of $\cover{f}$.  The typical setting in which it will be used is as follows: assume
\begin{eqnarray*}
   N_0 \cover{f} N_1, \quad R_{V,v}(N_1) \cover{f} N_2,
\end{eqnarray*}
and we would like to establish the existence of a point $x \in N_0$ such that $f(x) \in N_1$ and $f^2(x) \in N_2$.

\subsection{The mechanism of dropping dimensions---the main topological theorem }
%\label{sec:topThm}

\begin{theorem}
\label{thm:top-shadowing}

Assume that
\begin{eqnarray*}
  N_{0,0} \cover{f_{0,0}} N_{0,1} \cover{f_{0,1}} &\cdots& \cover{f_{0,i_0}} N_{0,i_0+1}, \\
  R_{V_0,\eta_0}N_{0,i_0+1}=  N_{1,0} \cover{f_{1,0}} N_{1,1} \cover{f_{1,1}} &\cdots& \cover{f_{1,i_1}} N_{1,i_1+1}, \\
   R_{V_1,\eta_1}N_{1,i_1+1}=  N_{2,0} \cover{f_{2,0}} N_{2,1} \cover{f_{2,1}} &\cdots& \cover{f_{2,i_2}} N_{2,i_2+1}, \\
   &\dots& \\
    R_{V_{L-1,\eta_{L-1}}}N_{L-1,i_{L-1}+1}=  N_{L,0} \cover{f_{L,0}} N_{L,1} \cover{f_{L,1}} &\cdots& \cover{f_{L,i_L}} N_{L,i_L+1}.
\end{eqnarray*}

Then there exist $q_0,\dots,q_L$, such that
\begin{eqnarray*}
  q_k &\in& N_{k,0},  \\
  q_{k+1}&=& f_{k,i_k}\circ \cdots \circ f_{k,1}\circ f_{k,0}(q_k), \quad k=0,\dots,L-1, \\
  f_{k,j}\circ \cdots \circ f_{k,1} \circ f_{k,0}(q_k) &\in& N_{k,j+1}, \quad  j=0,\dots,i_k, \quad k=0,\dots,L
\end{eqnarray*}

\end{theorem}
\textbf{Proof:}
Without any loss of the generality we can assume that $c_{N_{i,j}}=\mbox{id}$ for all $i,j$ appearing in the assumption of our theorem
and we can order our coordinates so that all our h-sets belong to a unit ball in
\[
  \mathcal{W}=V_0 \oplus V_1 \oplus V_2 \oplus \dots \oplus V_{L-1} \oplus W_u \oplus W_s,
\]
where each of $V_k$, $W_u$, $W_s$ are some $\mathbb{R}^d$ (with possibly different $d$),
such that for $N_{k,j}$ with $k=0,\dots,L$ and $j=0,\dots, i_k+1$, $V_k \oplus V_{k+1}\dots \oplus V_{L-1} \oplus W_u$ are exit directions
and $V_0 \oplus \dots \oplus  V_{k-1} \oplus W_s$ are entry directions.

We will use the following notation: $z\in N_{k,j}$ will be represented by a pair $(x,y)=(x(z),y(z))$, where
 $x \in V_k \oplus V_{k+1}\dots \oplus V_{L-1} \oplus W_u $, i.e. $x$ belongs to the "exit" subspace
 and $y \in V_0 \oplus \dots \oplus  V_{k-1} \oplus W_s$, i.e., $y$ belongs to the entry subspace. Observe that in this notation
   $x(z) \in W_u$ for $z \in N_{L,j}$ and $y(z) \in W_s$ for $z \in N_{0,j}$. $y_l(z)$ for $l<k$ and $x_l(z)$ for $l\geq k$ will denote the projection onto $V_l$. By $x_u(z)$ we will denote the projection on $W_u$ and by $y_s(z)$ we denote the projection on $W_s$. Observe that $x(z)=x_u(z)$ for $z \in N_{L,j}$ and $y(z)=y_s(z)$ for $z \in N_{0,j}$.

We will prove the following statement which implies our assertion:

For any $\bar{y}_s \in B_{W_s}(0,1)$ and $\bar{x}_u \in B_{W_u}(0,1)$, there exists $q_{k,j} \in N_{k,j}$ for $k=0,\dots,L$ and $j=0,1\dots,i_k+1$, such that
\begin{eqnarray}
  q_{k,j} &\in& N_{k,j}, \quad k=0,\dots,L \quad j=0,1\dots,i_k+1    \label{eq:sys-1}\\
  f_{k,j}(q_{k,j})&=&q_{k,j+1}, \quad  k=0,\dots,L \quad j=0,1\dots,i_k,  \label{eq:sys-2} \\
  q_{k+1,0}&=& q_{k,i_k+1}, \quad k=0,\dots,L-1  \label{eq:sys-3} \\
   y(q_{0,0})&=&\bar{y}_s,   \label{eq:sys-4} \\
   x_k(q_{k,i_k+1}) &=& \eta_k, \quad k=0,\dots,L-1  \label{eq:sys-5} \\
  x(q_{L,i_L+1})&=& \bar{x}_u  \label{eq:sys-6}
\end{eqnarray}
Equations (\ref{eq:sys-1},\ref{eq:sys-2},\ref{eq:sys-3}) simply mean that we travel through all sets $N_{k,j}$ in the desired direction.
Equation (\ref{eq:sys-4}) means that on the $W_s$ direction we can fix any value at the beginning, equation (\ref{eq:sys-5}) means that when dropping the
$V_k$ direction we can fix any  value of $x_k(q_{k+1,0})$ and (\ref{eq:sys-6}) implies that after the travel we can have arbitrary value on the $W_u$ direction.

In the sequel we will identify $q_{k+1,0}$ with $q_{k,i_k+1}$ but we will use both interchangeably depending on the context  to make notation and argument more readable

We consider the following system of equations
\begin{eqnarray}
  f_{k,j}(q_{k,j}) - q_{k,j+1}&=&0, \quad  k=0,\dots,L \quad j=0,1\dots,i_k,   \label{eq:sys-r1} \\
  x_k(q_{k,i_k+1}) - \eta_k&=&0, \quad k=0,\dots,L-1, \label{eq:sys-r2} \\
  y(q_{0,0})-\bar{y}_s &=& 0, \label{eq:sys-r3} \\
   x(q_{L,i_L+1}) - \bar{x}_u &=& 0. \label{eq:sys-r4}
\end{eqnarray}

We would like to prove that there exists a solution of this system in the set
\[
 D =\prod_{\substack{k=0,\dots,L\\j=0,\dots,i_k}} N_{k,j} \times N_{L,i_L+1}.
\]
Notice that since $q_{k+1,0}$ and $q_{k,i_k+1}$ are identified, for $k<L$ we do not have $N_{k,i_k+1}$ in the cartesian product defining $D$, but since $N_{L+1,0}$ does not exist we must include $N_{L,i_L+1}$.

Let us denote the r.h.s. of system (\ref{eq:sys-r1},\ref{eq:sys-r2},\ref{eq:sys-r3},\ref{eq:sys-r4}) by $F$.  We look for $q \in D$ such that
$F(q)=0$.

Denote by $w$ the dimension of $\mathcal{W}$.
Note that the dimension of $D$ is equal to $ \left(L+1+\sum_{k=0}^L i_k\right)w$.  Let us now count the number of equations in $F$.
\begin{itemize}
\item the number of equations in  (\ref{eq:sys-r1}) is equal to $w \sum_{k=0}^L (1+i_k)= w\left(L + \sum_{k=0}^{L}i_k\right)$
\item the number of equations is  (\ref{eq:sys-r2}) is equal to $\sum_{k=0}^{L-1} \mbox{dim} V_k$
\item in (\ref{eq:sys-r3}) we have $\mbox{dim} W_s$ equations
\item in (\ref{eq:sys-r4}) we have $\mbox{dim} W_u$ equations
\end{itemize}
Since $w=\sum_{k=0}^{L-1} \mbox{dim} V_k + \mbox{dim} W_s + \mbox{dim} W_u$ we see that the numbers of equations and variables in the equation $F(q)=0$, with $q \in D$, coincide.

We will prove that the system of equations  $F(q)=0$  has a
solution in $D$, by using the homotopy argument to show that the
local Brouwer degree $\deg (F,\inter D,0)$ is non-zero.

Let $h_{k,j}$  be the homotopies  of the covering relations from the assumptions of our theorem.

We embed $F$ into a one-parameter family of maps (a homotopy)
$H_t$ as follows (which we write as the system of equations $H_t(q)=0$):
\begin{eqnarray}
 h_{t,k,j}(q_{k,j}) - q_{k,j+1}&=&0, \quad  k=0,\dots,L \quad j=0,1\dots,i_k,   \label{eq:sys-hr1} \\
  x_k(q_{k,i_k+1}) - (1-t)\eta_k&=&0, \quad k=0,\dots,L-1, \label{eq:sys-hr2} \\
  y(q_{0,0})-(1-t)\bar{y}_s &=& 0, \label{eq:sys-hr3} \\
   x(q_{L,i_L+1}) - (1-t)\bar{x}_u &=& 0. \label{eq:sys-hr4}
 \end{eqnarray}
It is easy to see that  $H_0(q)=F(q)$.

We now show that if $q \in \partial D$ then $H_t(q) \neq 0$ is satisfied for all $t \in [0,1]$.
This will imply that $\deg(H_t,D,0)$ is
defined for all $t \in [0,1]$ and does not depend on $t$.

Take $q \in \partial D$. Then $q_{k,j} \in \partial N_{k,j}$ for some $k,j$ appearing in the definition of $D$.
We have the following four cases:
\begin{itemize}
\item $q_{0,0} \in \partial N_{0,0}$.  If $q_{0,0}\in N_{0,0}^+$ then  $\|y(q_{0,0})\|=1 > \|\bar{y}_s\|$, so  equation  (\ref{eq:sys-hr3}) is not satisfied for all $t \in [0,1]$. If $q_{0,0}\in N_{0,0}^-$, from the exit condition (\ref{eq:hc2}) in the covering relation $N_{0,0} \cover{f_{0,0}} N_{0,1}$ it follows that
    $h_{t,0,0}(q_{0,0}) \notin N_{0,1}$, and therefore the $(0,0)$-th equation in (\ref{eq:sys-hr1}) does not hold for any $q_{0,1} \in N_{0,1}$.

\item $q_{k,0} \in \partial N_{k,0}$ for $k>0$. Let us recall that $N_{k,0}=R_{V_{k-1}}N_{k-1,i_{k-1}+1}$, both with the same support, and if we denote
     \begin{eqnarray*}
     Z_u=V_k \oplus \dots V_{L-1} \oplus W_u, \quad Z_s=V_0 \oplus \dots  \oplus V_{k-2} \oplus W_s,
     \end{eqnarray*}
     then the exit directions are  $Z_u$  in $N_{k,0}$ and $Z_u \oplus V_{k-1}$ in $N_{k-1,i_{k-1}+1}$. The entry directions are
        $Z_s\oplus V_{k-1}$ for $N_{k,0}$ and $Z_s$ for $N_{k-1,i_{k-1}+1}$.

     If $q_{k,0} \in \partial N_{k,0}$, we have one of the following three cases
      \begin{itemize}
        \item $\|\pi_{Z_u} q_{k,0}\|=1$, hence $q_{k,0} \in N^-_{k,0}$. From the exit condition (\ref{eq:hc2}) in the covering relation $N_{k,0} \cover{f_{k,0}} N_{k,1}$ it follows that
        $h_{t,k,0}(q_{k,0}) \notin N_{k,1}$ so the $(k,0)$-th equation in (\ref{eq:sys-hr1}) is not satisfied.

        \item $\|\pi_{V_{k-1}} q_{k,0}\|=1$. We will show that in this case equation (\ref{eq:sys-hr2}) for $k-1$ is not satisfied. Indeed  in view of our identification
           of $q_{k,0}$ with $q_{k-1,i_{k-1}+1}$, equation (\ref{eq:sys-hr2}) becomes
              $x_{k-1}(q_{k,0})-(1-t)\eta_{k-1}=0$.
              But we have  $\|x_{k-1}(q_{k,0})\|=1 > \|(1-t)\eta_{k-1}\|$.

        \item $\|\pi_{Z_s} q_{k,0}\|=1$. Then $q_{k,0} \in N^+_{k-1,i_{k-1}+1}$, and from the entry condition (\ref{eq:hc3}) in relation $N_{k-1,i_{k-1}} \cover{f_{k-1,i_{k-1}}} N_{k-1,i_{k-1}+1} $ we know that
            \begin{equation*}
            h_{t,k-1,i_{k-1}}(q_{k-1,i_{k-1}}) \neq q_{k,0}, \quad \forall q_{k-1,i_{k-1}} \in N_{k-1,i_{k-1}}.
       \end{equation*}
       Hence the $(k-1,i_{k-1})$-th equation in (\ref{eq:sys-hr1}) is not satisfied.
      \end{itemize}

\item $q_{L,i_L+1} \in \partial N_{L,i_L+1}$. If $q_{L,i_L+1} \in N^+_{L,i_L+1}$, then from the entry condition (\ref{eq:hc3}) in the covering relation $N_{L,i_L} \cover{f_{L,i_L}} N_{L,i_{L}+1} $ we know that $h_{t,L,i_L}(q_{L,i_L}) \neq q_{L,i_L+1}$     for all $q_{L,i_L} \in N_{L,i_L}$.
    If $q_{L,i_L+1} \in N^-_{L,i_L+1}$ then equation (\ref{eq:sys-hr4}) is not satisfied, because
    $\|q_{L,i_L+1}\|=1 > \|\bar{x}_u\|$.

\item $q_{k,j} \in \partial N_{k,j}$ with $1 \leq j \leq i_k$. Then we have one of the following two cases
  \begin{itemize}
   \item $q_{k,j} \in N^+_{k,j}$.
   From the covering relation $N_{k,j-1} \cover{f_{k,j-1}} N_{k,j} $ we know that $h_{t,k,j-1}(q_{k,j-1}) \neq q_{k,j}$
   for all $q_{k,j-1} \in N_{k,j-1}$ (this follows from the entry condition (\ref{eq:hc3})), hence $(k,j-1)$-th equation in (\ref{eq:sys-hr1}) is not satisfied.
   \item $q_{k,j} \in N^-_{k,j}$.
   From the exit condition (\ref{eq:hc2}) in the covering relation $N_{k,j} \cover{f_{k,j}} N_{k,j+1}$ it follows that
   $h_{t,k,j}(q_{k,j}) \notin N_{k,j+1}$ so the $(k,j)$-th equation in (\ref{eq:sys-hr1}) is not satisfied.
\end{itemize}
\end{itemize}

We have thus proved that $\deg(H_t,\inter D,0)$ is defined. By the
homotopy invariance we have
\[
  \deg(F,\inter D,0)=\deg(H_1,\inter D,0). \label{eq:deg-cont}
\]

Observe that $H_1(q)=0$ is the following system of \emph{linear
equations}
\begin{eqnarray*}
  A_{k,j}x(q_{k,j}) - x(q_{k,j+1})&=&0, \quad  k=0,\dots,L \quad j=0,1\dots,i_k,\\%    \label{eq:sys-h1-x} \\
   - y(q_{k,j+1})&=&0, \quad  k=0,\dots,L \quad j=0,1\dots,i_k, \\%  \label{eq:sys-h1-y} \\
  x_k(q_{k,i_k+1})&=&0, \quad k=0,\dots,L-1, \\%\label{eq:sys-h1-xk} \\
  y(q_{0,0}) &=& 0,\\% \label{eq:sys-h1-ys} \\
   x(q_{L,i_L+1}) &=& 0. %\label{eq:sys-h1-xuend}
\end{eqnarray*}
where $A_{k,j}$ is a linear map which appears at
the end of the homotopy $h_{k,j}$.

 It is immediate that $q=0$ is a solution of this system. This must be a unique solution, because otherwise the solution
set will be a linear subspace which would intersect $\partial D$. But from the previous part of the proof we know that there is no solution
for this system on $\partial D$. Therefore the determinant of the matrix defining this system  is not zero and its sign is equal to
 $\deg(H_1,\inter D,0)$.
This and \eqref{eq:deg-cont} implies that
\begin{equation}
  \deg(F,\inter D,0)=\pm 1,
\end{equation}
hence there exists a solution of equation $F(q)=0$ in $D$. This
finishes the proof.
\qed

%% file: enclosure.tex
\section{Enclosures of solutions for the local system in the polynomial normal form}
\label{sec:encl-full-system}

The purpose of this section is to exhibit bounds, more precisely \emph{enclosures)}, for the solutions of the polynomial system~\eqref{eq:NF} of Lemma~\ref{lem:PolyNormForm},
which is the normal form for our TMS on the $j$-th chart. Without any loss of generality we will assume $\lambda=1$, simply by scaling the time.

We fix $j$ to be the index of the torus on which our chart is centered.
We use the notation for monomials introduced in section~\ref{sub:polynomialNFc}, where there is also the notion of a resonant monomial~\eqref{eq:resonant-monomial} for the
saddle variables.

Throughout this section we will fix compact neighborhoods $Z\pm$ of the origin for the variables $z_\pm=\left(x_\pm,y_\pm\right)$ with the norms
$\abs{z\pm}=\abs{x_\pm}+\abs{y_\pm}$ as well as a compact neighborhood $C_*$ of the origin for the variable $c_*=\left(c_1,\dots,c_{j-2}, c_{j+2},\dots,c_n\right)$
with the norm $\abs{c_*}=\sum_{\ell \neq j,j\pm 1}\abs{c_\ell }$, merging them in a fix compact neighborhood $Z=Z_-\times Z_+\times C_*$ of the variable $z:=(z_-,z_+,c_*)$ with the norm $\abs{z}=\abs{z_-}+\abs{z_+}+\abs{c_*}$.
Finally, we will denote $\alpha=[-1,1]$.

\subsection{The enclosure of solutions for saddle directions}
%\label{subsec:full-encl}

\begin{theorem}
\label{thm:enclo-hyp-dir}
Consider an autonomous system depending on a parameter $\rho \in [0,1]$
\begin{eqnarray}
  \dot{x}_-&=&  x_- + \rho f_{x_-}(z), \label{eq:sys-res-0-full-xy-1} \\
  \dot{y}_-&=& -y_- + \rho f_{y_-}(z),  \notag \\
  \dot{x}_+&=&  x_+ + \rho f_{x_+}(z) ,  \notag \\
  \dot{y}_+&=& -y_+ + \rho f_{y_+}(z), \label{eq:sys-res-0-full-xy-4} \\
  \dot{c}_*&=&  f_{c_*}(\rho,z), \notag %\label{eq:sys-res-0-full-c}
\end{eqnarray}
for $z=(x_-,y_-,x_+,y_+,c_*)=(z_-,z_+,c_*)\in Z=Z_-\times Z_+\times C_*$, with initial conditions satisfying
\begin{equation}
  x_-(0) = a_0e^{-2T}, \quad y_-(0) =\eta, \quad x_+(0) =  d_0 e^{-T}, \quad y_+(0) =0,
 \label{eq:sys-ic}
\end{equation}
where $|d_0|,|\eta| \leq  3\sigma/2$ and $a_0 \in [-T^k/2,T^k/2]$ for some integer $k \geq 1$.

%\noindent\textbf{A: Put here the EXACT hypotheses of Lemma~\ref{lem:PolyNormForm}}

Assume that for some $k_c \geq 0$ and $B \leq  B_0$
\begin{equation}
   \abs{c_\ell (t)}\leq B T^{k_c} e^{-T}, \quad \ell \neq j,j\pm 1, \qquad  t \in [0,T].  \label{eq:it-sys-ic}
\end{equation}
\textbf{Comment: we have $B$ and $B_0$ because of the following, once we have $T> T_0(B_0,\dots)$ big enough for our assertion to hold, we would like to have a freedom to decrease $B$ if necessary. }

We assume that $f_v(z)$, for $v \in \{x_-,y_-,x_+,y_+\}$, can be split as
\[
  f_v(z)=\sum_{m \in M_{1,v}}  g_{v,m}(z)z^m  + \sum_{m \in M_{2,v}} g_{v,m}(z)z^m,  %\label{eq:f-full-enclo}
\]
where $g_{v,m}(z)$ are bounded continuous functions
\[
  \abs{g_{v,m}(z)} \leq K, \qquad \forall v,m,\quad \forall z \in Z,  % \label{eq:K}
\]
$M_{1,v},M_{2,v}$ are \emph{finite} sets of indices,
and any $z^m$, with $m=\left(m_{x_-},m_{y_- },m_{x_+},m_{y_+},m_c\right)$, is a resonant monomial~\eqref{eq:resonant-monomial} for the saddle variables,
satisfying, on the one hand $m_s:=m_{x_-}+m_{y_- }+m_{x_+ } +m_{y_+}\geq 3$ and $m_c =0$,  if $m\in M_{ 1, v }$,
and, on the other hand, $m_s=1$ and $m_c \geq 3$, if $m\in M_{2,v }$.
%and any $z^m$ is a resonant monomial~\eqref{eq:resonant-monomial} for the saddle variables,
%satisfying on the one hand $m_s:=m_{x_-}+m_{y_- }+m_{x_+ } +m_{y_+}\geq 2$ and $m_c =0$  if $m=\left(m_{x_-},m_{y_- },m_{x_+},m_{y_+},m_c\right) \in M_{ 1, v }$,
%and on the other hand $m_s=1$ and $m_c \geq 3$ if $m\in M_{2,v }$.

%We assume that
%\begin{description}
%  \item[A3] the following bad terms are excluded, i.e. such $w_J$ do not appear in (\ref{eq:f-full-enclo})
%    \begin{itemize}
%      \item in $f_{y_+}$:   illegal term:  $O(c_*)y_-$.
%    \end{itemize}
%\end{description}

\textbf{$A\geq 3$ is needed later for covering relations;  and $A\geq 2$ is needed in the estimates in the proof }

Take $A \geq 3$. Then there exists $\sigma_0=\sigma_0(A,B_0,k,k_c)>0$   \textbf{PZ: check whether $\sigma_0$ depends on $A$?} such that  for all $0<\sigma \leq \sigma_0$  there exists $T_0=T_0(\sigma,A,B_0,k,k_c)$   such that  for all $T\geq T_0$ the following enclosures hold for all $t \in [0,T]$ and for all $\rho \in [0,1]$
\textbf{PZ: does $T_0$ depends on $G_{1,v}$,$G_{2,v}$ - definitely yes }
\begin{equation}
\label{eq:fen}
\begin{split}
   z(t)=(x_-(t),y_-(t),x_+(t),y_+(t),c_*(t)) &\in  Z \\
  x_-(t) &\in e^{-2T}e^{t} \left(a_0 + \alpha T^{k}/A    \right),  \\
   y_-(t) &\in e^{-t}\left(\eta + \alpha  \sigma K /A \right), \\
  x_+(t) &\in  e^{-T}e^t\left(d_0 + \alpha \sigma/A \right),\\
  y_+(t) &\in \alpha e^{-t} e^{-T}  K \sigma t/A.
\end{split}
\end{equation}
%\begin{eqnarray}
%   (y_-(t),x_-(t),y_+(t),x_+(t),c_*(t)) &\in&  Z \times Z_c, \notag \\
%   y_-(t) &\in& e^{-t}\left(\eta + \alpha  \sigma K /A \right), \label{eq:fen-y-}\\
%  x_-(t) &\in& e^{-2T}e^{t} \left(a_0 + \alpha T^{k}/A    \right), \label{eq:fen-x-} \\
%  y_+(t) &\in& \alpha e^{-t} e^{-T}  K \sigma t/A, \label{eq:fen-y+} \\
%  x_+(t) &\in&  e^{-T}e^t\left(d_0 + \alpha \sigma/A \right). \label{eq:fen-x+}
%\end{eqnarray}
\end{theorem}

%The proof of this theorem is technical and can be found in Appendices~\ref{sec:proof-estm} (the conceptual part) and \ref{sec:good-terms} (the computational part).

%Using the concrete hypotheses above of this theorem, its proof is straightforward but somewhat technical, so it is postponed to Appendix~\ref{sec:proof-estm} (the conceptual part) and Appendix~\ref{sec:good-terms} (the computational part).

Using the concrete hypotheses of this theorem, its proof is straightforward but somewhat technical, so it is postponed to Appendix~\ref{sec:proof-estm} (the conceptual part) and Appendix~\ref {sec:good-terms} (the computational part).

\subsection{The estimates along the center direction}
%\label{subsec:tms-center-estm}
The goal of this subsection is to complete the estimates for the flow near $\mathbb{T}_j$  by providing bounds for $c_\ell $ for $\ell \neq j,j\pm 1$.

In the saddle directions we assume the following enclosures on $t \in [0,T]$ (compare with Theorem~\ref{thm:enclo-hyp-dir})
\begin{equation}
\label{eq:v-bnd}
\begin{split}
  x_-(t) &\in \alpha e^{-2T}e^t T^{k_1},\\
  y_-(t) &\in \alpha e^{-t} d, \\
  x_+(t)  &\in  \alpha e^{-T} e^{t} d,\\
  y_+(t) &\in \alpha  e^{-T}e^{-t} T^{k_1},
\end{split}
\end{equation}
where
\begin{eqnarray}
d &\leq& 2\sigma  \label{eq:dleq2sigma}, \\
k_1 &>& k. \label{eq:k1-cdir-thm}
\end{eqnarray}
\textbf{PZ: why to keep (\ref{eq:dleq2sigma}) and not write  $d=2\sigma$ and then in fact remove $d$. Analogously $k_1=k+1$ or even $k$}

Observe that for $T$ large enough we have $\left(z_-(t),z_+(t)\right) \in Z_-\times Z_+$.
%\textbf{A: Why? Here?}

On the center direction we assume the differential equation (compare with \eqref{eq:NF})
\begin{equation}
 \dot{c_\ell }= i c_\ell (\nu_\ell + \rho g_\ell (z)) , \qquad \ell \neq j,j\pm 1,  \label{eq:sys-res-0-full-center-1}
\end{equation}
where $\rho \in [0,1]$ and $\nu_\ell  \in \mathbb{R}$.

%\textbf{PZ: for a moment I was thinking that it is necessary to add to rhs of (\ref{eq:c-eq-estm}) term \\
%$\bar{c}_j O(|y_+|+|y_-|+|x_-|+|x_+| + |c_*|)$. But this is not need as}
%\begin{eqnarray*}
%\bar{c}_j O(|y_+|+|y_-|+|x_-|+|x_+| + |c_*|) = \\
%c_j \frac{\bar{c}_j}{c_j}O(|y_+|+|y_-|+|x_-|+|x_+| + |c_*|)= \\
%c_j O(|y_+|+|y_-|+|x_-|+|x_+| + |c_*|)
%\end{eqnarray*}

We assume that there exists a constant $G >0$ such that for all $\ell$ the following bound holds
\begin{equation}
 \abs{g_\ell (z)} \leq G\abs{z}, \qquad z \in  Z .   \label{eq:c-eq-estm}
\end{equation}
%\noindent\textbf{A: If \boldmath$\abs{g_\ell (z)} \leq G\abs{z}^2$ then $ |c_\ell (t)|^2 < |c_\ell (0)|^2 \exp(5 \cdot G d^2)$}

The next theorem is the adaptation of \cite{CK} estimates to our setting.
\begin{theorem}
\label{thm:center-estm}
Assume the enclosures~\eqref{eq:v-bnd} for the saddle variables and the evolution of the center variables $c_\ell $ for $\ell \neq j,j\pm 1$
%(\ref{eq:y-bnd}--\ref{eq:x+bnd})~
are given by~(\ref{eq:sys-res-0-full-center-1}) under estimate~(\ref{eq:c-eq-estm}) and initial conditions
\begin{equation}
  c_\ell (0)=u_\ell  e^{-T} , \quad u_\ell  \in [-T^{k_c},T^{k_c}].  \label{eq:center-ic-size}
\end{equation}

Then there exists $T_0(k_c,k_1,G, n,d)$, such that for $T \geq T_0$ the following inequalities hold if $c_\ell (0) \neq 0$

\begin{equation}
   |c_\ell (0)|^2 \exp(-5 \cdot G d)   < |c_\ell (t)|^2 < |c_\ell (0)|^2 \exp(5 \cdot G d), \quad t \in [0,T]. \label{eq:cj-lup-estm}
\end{equation}
\end{theorem}
\begin{rem}
Note that although the sizes of $u_\ell $ may be different, all of them are bounded by $T^{k_c}$.
\end{rem}

\noindent
\textbf{Proof:}
All relevant estimates in the proof are increasing with $\rho$, so that it is enough to  consider $\rho=1$.

From (\ref{eq:sys-res-0-full-center-1}) we obtain
\begin{eqnarray*}
  \frac{d}{dt}|c_\ell (t)|^2=\dot{c}_\ell  \overline{c}_\ell  + c_\ell  \dot{\overline{c}}_\ell = - 2 |c_\ell |^2 \im g_\ell ((z(t)),
\end{eqnarray*}
so that
\begin{equation*}
  -2|c_\ell (t)|^2 \cdot |g_\ell (z(t))| \leq \frac{d}{dt}|c_\ell (t)|^2 \leq 2|c_\ell (t)|^2 \cdot |g_\ell (z(t))|,
\end{equation*}
and
\begin{equation}
|c_\ell (0)|^2 \exp\left(-2 \int_0^t |g_\ell (z(s))|ds\right)  \leq  |c_\ell (t)|^2 \leq |c_\ell (0)|^2 \exp\left(2 \int_0^t |g_\ell (z(s))|ds\right). \label{eq:cj-low-upp-estm}
\end{equation}

Therefore, in view of~(\ref{eq:c-eq-estm}), it is enough  to show  that we  have a finite bound
for
\[
\int_0^T\abs{z(s)} ds=\int_0^T\left(|x_-(s)|+ |y_-(s)| + |x_+(s)|+|y_+(s)|+ \sum_{\ell \neq j,j\pm 1}\abs{c_\ell (s)}\right) ds.
\]
As the main step to achieve this goal we will show that there exists $B>1$ such that
\begin{equation}
   |c_\ell (t)| \leqslant B |u_\ell | e^{-T}  \label{eq:cj-apriori}
\end{equation}
holds for all $t \in [0,T]$.
We will use the continuation argument.

\noindent\textbf{Concrete reference for the continuation argument or continuation method or continuation-type argument}

Let us fix any $B>1$, whose precise value will be found during the proof, indeed in~\eqref{eq:central-B2}, so that the bound~\eqref{eq:cj-apriori} is satisfied for $t \in [0,T']$, where $T'>0$. %\textbf{PZ: this uses that fact that we work in finite dimension. To work in infinite dimension we will use something like the Picard iteration, for this we need the uniqueness for this system}
We will show that $T' \geq T$.

From~\eqref{eq:v-bnd} and~\eqref{eq:cj-apriori}, we have for  any $t \in [0,\min(T',T)]$
\begin{equation}
\label{eq:int}
\begin{split}
  \int_0^{t} |x_-(s)| ds &\leq  \int_0^{t} e^{s} e^{-2T} T^{k_1}ds < e^{-T}T^{k_1},\\%  \label{eq:int-x-} \\
    \int_0^t |y_-(s)|ds &\leq d, \\%\label{eq:int-y-} \\
  \int_0^{t} |x_+(s)| ds &\leq d e^{-T} \int_0^t e^{s} ds < d e^{-T} e^T=d, \\% \label{eq:int-x+}\\
  \int_0^{t} |y_+(s)| ds &\leq \int_0^{t} e^{-s}e^{-T} T^{k_1} ds < e^{-T} T^{k_1},\\% \label{eq:int-y+} \\
 \int_0^{t} |c_\ell (s)| ds  &\leq \int_0^t B |u_\ell | e^{-T}ds \leq B |u_\ell | T e^{-T}, \quad  \ell \neq j,j\pm 1,% \label{eq:int-cj}
 \end{split}
 \end{equation}

By combining the inequalities of~\eqref{eq:int} we  obtain for any $t \in [0,\min(T',T)]$
\[
  %\int_0^{t} \left(\sum_{j \neq j_0,j \neq j_0 \pm 1}|c_j(s)|\right) + |x_+(s)|+ |x_-(s)| + |y_+(s)|+ |y_-(s)| ds
  \int_0^T\abs{z(s)} ds \leq
  2T^{k_1}e^{-T}+n B   T^{k_c+1} e^{-T} +  2d =E+2d,
\]
where we have introduced
\[
E=2T^{k_1}e^{-T}+n B   T^{k_c+1} e^{-T}.
\]
Using~\eqref{eq:c-eq-estm}, for any $t \in [0,\min(T',T)]$ and $\ell \neq j,j\pm 1$ we have
\begin{equation}
  \int_0^{t} |g_\ell (z(s))|ds \leq G E + 2Gd.  \label{eq:estm-int-g}
\end{equation}

For the continuation argument we need to satisfy
\[
   |c_\ell (T')| < B |u_\ell | e^{-T},  \quad \forall T' \in [0,T] \quad \forall \ell \neq j,j\pm 1.  %\label{eq:cj-cont-arg}
\]

By (\ref{eq:center-ic-size},\ref{eq:cj-low-upp-estm},\ref{eq:estm-int-g}) it is enough to have %\textbf{PZ: here we use again the finite-dimensionality}
\begin{equation*}
  |u_\ell |^2 e^{-2T} \exp\left( 2GE + 4Gd\right) < |u_\ell |^2 B^2 e^{-2T},
\end{equation*}
which is equivalent to
\begin{equation*}
  \exp\left( 2GE + 4Gd\right) <  B^2.
\end{equation*}
It is clear that this can be achieved if
we take $B$ such that
\begin{equation}
  \exp\left(5  G  d \right) = B^2 \label{eq:central-B2}
\end{equation}
and then we take $T \geq T_0(k_c,k_1,G, n,d)$ large enough to have
\[
E=2T^{k_1}e^{-T}+n B   T^{k_c+1} e^{-T}<d/2.
\]

Finally, from (\ref{eq:cj-low-upp-estm}) we have also obtained the required estimate~\eqref{eq:center-ic-size}.
\qed

\subsection{Full system---the saddle and center directions}

\begin{theorem}
\label{thm:ful-system}
Consider the system of ODEs  depending on parameter $\rho$
\begin{eqnarray*}
  \dot{x}_-&=&  x_- + \rho f_{x_-}(z),   \\
 \dot{y}_-&=& -y_- + \rho f_{y_-}(z),     \\
  \dot{x}_+&=&  x_+ + \rho f_{x_+}(z) ,  \\
  \dot{y}_+&=& -y_+ + \rho f_{y_+}(z),  \\
  \dot{c_\ell }&=& i c_\ell \left(\nu_\ell + \rho g_\ell (z)\right) , \qquad \ell \neq j,j\pm 1.
\end{eqnarray*}
We assume that the terms  $f_v(z)$ and $g_\ell (z)$ appearing in the equations satisfy the hypotheses listed in Theorem~\ref{thm:enclo-hyp-dir} for the saddle
directions and (\ref{eq:sys-res-0-full-center-1},\ref{eq:c-eq-estm}) for the center directions.

Assume that there exist integers $k \geq 1$ and $k_c \geq 0$ such that
the initial conditions satisfy  (\ref{eq:sys-ic}) in the saddle directions and (\ref{eq:center-ic-size}) in the center direction.

Take $A \geq 3$. Then there exists $\sigma_0>0$  such that for all $0<\sigma \leq \sigma_0$ there is $T_0(\sigma)$ such that for all $T\geq T_0$ the following
conditions hold for  $t \in [0,T]$ and all $\rho \in [0,1]$
%\begin{eqnarray}
% (y_-(t),x_-(t),y_+(t),x_+(t),c_*(t)) &\in&  Z \times C_*, \notag
%\end{eqnarray}
%and
%\begin{eqnarray}
%   y_-(t) &\in& e^{-t}\left(\eta + \alpha  \sigma K /A \right),\notag \\
%  x_-(t) &\in& e^{-2T}e^{t} \left(a_0 + \alpha T^{k}/A   \right) \notag\\
%  y_+(t) &\in& \alpha e^{-t} e^{-T}  K \sigma t/A, \notag\\
%  x_+(t) &\in&  e^{-T}e^t\left(d_0 + \alpha \sigma/A \right),   \label{eq:fencl-x+} \\
%   |c_j(0)|^2 \exp(-10 \cdot G \sigma)   &<& |c_j(t)|^2 < |c_j(0)|^2 \exp(10 \cdot G \sigma),  \label{eq:fencl-cj}
%    \mbox{if} \, c_j(0) \neq 0, \, j \neq j_0,j_0 \pm 1.  %\notag
%\end{eqnarray}
%
\begin{align}
 z(t)=(y_-(t),x_-(t),y_+(t),x_+(t),c_*(t)) &\in  Z , \notag\\
   y_-(t) &\in e^{-t}\left(\eta + \alpha  \sigma K /A \right),\notag \\
  x_-(t) &\in e^{-2T}e^{t} \left(a_0 + \alpha T^{k}/A   \right) \notag\\
  y_+(t) &\in \alpha e^{-t} e^{-T}  K \sigma t/A, \notag\\
  x_+(t) &\in  e^{-T}e^t\left(d_0 + \alpha \sigma/A \right),   \label{eq:fencl-x+} \\
   |c_\ell (0)|^2 \exp(-10 \cdot G \sigma)   < |c_\ell (t)|^2 &< |c_\ell (0)|^2 \exp(10 \cdot G \sigma),  \label{eq:fencl-cj}
    \text{ if } c_\ell (0) \neq 0, \, \ell \neq j,j\pm 1.  %\notag
\end{align}

\end{theorem}
\textbf{Proof:}
Our assertion  combines   Theorems~\ref{thm:enclo-hyp-dir} and~\ref{thm:center-estm}. These theorems are not independent, hence some care
should be taken.

For Theorem~\ref{thm:enclo-hyp-dir} we take $B_0= \exp\left(5  G d/2 \right)$ (see (\ref{eq:cj-lup-estm})) and
 $\sigma=1/2$, so $d=2\sigma=1$  (see (\ref{eq:dleq2sigma})).
Note that we have the freedom of decreasing $\sigma$ (and $d$) and this $B_0$ will still be good.

For this $B_0$ we apply Theorem~\ref{thm:enclo-hyp-dir} to obtain enclosures for $(x_{\pm}(t),y_{\pm}(t))$ for $t \in [0,T]$, $\sigma$ small enough and $T$ big enough, with the freedom  to further increase $T$ if needed.

Then we see that assumptions of Theorem~\ref{thm:center-estm} are satisfied with $k_1$ satisfying (\ref{eq:k1-cdir-thm})
and $T>T_0$ large enough.

Once this $T_0$ is found, we can do a formal continuation argument or the Picard iteration to get the final result.

\qed

%% file: scheme.tex
\section{Proof of the main theorems---The construction of covering relations}
\label{sec:scheme}

The goal of this section is to construct a sequence of covering relations for the generalized toy model system introduced in Section~\ref{sec:genToyModel}  in order to follow
 the heteroclinic chain $\mathbb{T}_0 \to \mathbb{T}_1 \to \cdots \to \mathbb{T}_n$, which will give the proof of  Theorems~\ref{thm:toymodel} and~\ref{thm:toymodel-gen}.% and~\ref{thm:localModel-tran}.

We will denote by $\varphi(t,x_0)=\varphi_t(x_0)$ the (local) flow induced by the generalized toy model system, that is, the shift during a time $t$ along the orbit with the initial condition. $ x(0)=x_0 $.

The h-sets in the covering relations will be defined over the  $j$-th \emph{normal formal coordinates} charts $\{(z_-,z_+,c_*)\}$, with $z_\pm=(x_\pm, y_\pm)$ obtained
via Lemma~\ref{lem:PolyNormForm}, in which our local system is in polynomial normal form. In these coordinates the
heteroclinic orbits $H_\inn$ and $H_\out$ are locally straight and given by
$H_{\inn}=\{(z_-,0,0):  x_-=0, y_- >0\}$, and $H_{\out}=\{(0,z_+,0):  y_+=0, x_+>0\}$.

We will denote by $\mathcal{C}_j$ the transformation that expresses the $j$-th \emph{normal formal coordinates} in terms of the coordinates of the $j$-th chart
introduced in~(\ref{eq:j-chart-polar}, \ref{eq:j-chart-ck}).

We fix $\sigma>0$ (some small number whose size is decided by Theorem~\ref{thm:ful-system}). For any $j$, denote $A_j$ the point on the heteroclinic orbit from
$\mathcal{T}_{j-1}$ to $\mathcal{T}_{j}$ such that
in the $j$-th  normal form coordinates satisfies
\begin{equation}
 y_-(A_j)=\sigma, \qquad (x_-,z_+, c_*)(A_j)=0 \label{eq:Aj}
\end{equation}
and $B_j$   be the point on the heteroclinic orbit from
$\mathcal{T}_{j}$ to $\mathcal{T}_{j+1}$ such that
in the same $j$-th  normal form coordinates satisfies
\begin{equation}
x_+(B_j)=\sigma, \quad (z_-,y_+,c_*)(B_j)=0. \label{eq:Bj}
\end{equation}
Observe that there exists $t_j >0$ such that $\varphi(t_j,B_j)=A_{j+1}$ for each $j$. We assume that we have a uniform upper bound for $t_j$ (in fact if we have translational
symmetry of our lattice then $t_j$ does not depend on $j$).

Travelling along the heteroclinic requires also taking into account several coordinate changes.  In our setting
it will be the composition of the following maps:
\begin{itemize}
\item $\mathcal{C}_j^{-1}$ in the neighborhood of $B_j$, which maps $j$-th normal form coordinates  to  the  $j$-th chart coordinates.
\item propagation along the heteroclinic from
$\T_{j-1}$ to $\T_{j}$ by $\varphi_{t_j}$
\item switching to $(j+1)$-th chart, i.e. application of $\J_{j\to j+1}$ \eqref{eq:transitionj} introduced in Section~\ref{subsubsec:cc-charts}.
\item changing coordinate to the $(j+1)$-th normal form coordinates by application of $\mathcal{C}_{j+1}$.
\end{itemize}
Accordingly we define a transition  map $\T_j$ by
\begin{equation*}
  \T_j =\mathcal{C}_{j+1} \circ \J_{j\to j+1}  \circ \varphi_{t_j} \circ \mathcal{C}_j^{-1}. %\label{eq:Tj-def}
\end{equation*}
Since we compose maps with the block-diagonal structure, the map $\T_j$ also has this structure, with the only nontrivial blocks being the same as for the map $ \J_{j\to j+1} $,
which are listed in Section~\ref{subsubsec:cc-charts}.

\subsection{Some heuristic about the sizes}
%\label{subsec:cov-sch-heu}

When defining our h-sets
for heuristic reasons we will deal with three different sizes: \emph{macro}, \emph{micro} and \emph{nano}.
 The macro size will be $O(1)$, the micro size $w(T)e^{-T}$ and the nano size $w(T) e^{-2T}$, where $w(T)$ are some polynomials in the variable $T$ (possibly different for each variable). By increasing $T$ we will be able
 to follow the heteroclinic chain as close as we desire, keeping the relative sizes of the variables.
 This idea comes from \cite{CK}.

Compared to~\cite{DSZ} the strategy we adopt here is different, to be able to deal with the nonlinearities present both in the transition along the heteroclinic orbits and in the coordinate changes. In particular this will mean that the exit variables will never be of nano size (as was the case in~\cite{DSZ} for the $y_+$ variable).

The realization of our strategy must take into account the following issues.
\begin{itemize}
\item We declare that $x_-$ is an entry direction, despite being unstable.
\item In the center directions we may also have some growth or decay that might be undesirable for the past modes and for the future modes, respectively.
\end{itemize}

To overcome these difficulties we must carefully choose the relative sizes of the variables. We will follow the same principles as in \cite{DSZ}, which we recall below.

In the evolution related to passing near $\mathbb{T}_j$ in the $j$-th chart, at the end in the exit section $x_+=\sigma$ we want $y_+$ to be of nano size and all other variables of micro size. To achieve these sizes we have to impose in the entry section $y_-=\sigma$ that $x_-$ be of nano size and the same for $y_+$ as it will then become $x_-$ in the next  $(j+1)$-th chart, while all other variables should be of micro size.

Regarding the behavior of the center directions, it turns out that it is possible to keep the micro size throughout this transition for the $c_\ell$ variables for $\ell \leq j-2$ or $\ell \geq j+2 $, because the possible decay or growth is limited by a constant, which does not depend on $T$  (see assertion (\ref{eq:fencl-cj}) in Theorem~\ref{thm:ful-system}).

\subsection{On h-sets and covering relations}

We will have four types of  $h$-sets in the $j$-th normal form coordinates chart,  $N^j_{\inn}$, $\widetilde{N}^j_{\inn}$, $N^j_{\out}$, $\widetilde{N}^j_{\out}$  (objects with tilde are   contracted h-sets in the sense of Def.~\ref{def:dropping-exit-dim}), chosen such that
 the following covering relations are satisfied
\begin{eqnarray}
  R_{<y_+>,0}N^j_{\inn}=\widetilde{N}^j_{\inn} &\cover{\varphi_T}& N^j_{\out}, \quad j=0,\dots,n\label{eq:tms-nin-nout}\\
  R_{<x_+>,\sigma}N^j_{\out}=\widetilde{N}^j_{\out}&\cover{\mathcal{T}_j}& N^{j+1}_{\inn},\quad j=0,\dots,n-1 \label{eq:tms-tran-hetero-cv}.
\end{eqnarray}
 Recall that (see  Definition~\ref{def:dropping-exit-dim}) $R_{<x_+>,\sigma} N^{j}_{\out}$ means that we drop the coordinate direction $x_+$ and we set $x_+=\sigma$.

All our h-sets are defined as  h-sets  with the product structure (see Def.~\ref{def:hset-product}). The decomposition giving the product structure are always in coordinate directions  or groups of variables.

 We will use the following conventions
\begin{itemize}
\item $\gamma(variable)$ - will be used for radius  in the entry directions,
\item $r(variable)$ - will be used for radius  in the exit directions.
\end{itemize}

By $c_p$ (the past modes) we will denote the collection $\{c_\ell\}_{\ell \leq j-2}$ and by $c_f$ (the future modes) we will denote the collection $\{c_\ell\}_{\ell \geq j+2}$.

 On $c_p$ and $c_f$ we use the sup norm, i.e. $\|c_p\|=\sup_{\ell \leq j-2} |c_\ell|$.
For blocks $z_s=(x_s,y_s)$ with $s \in \{+,-\}$ we will always use the norm $|z_s|=\max(|x_s|,|y_s|)$. \textbf{PZ: check the norms defined earlier by Amadeu, make adjustemnt if needed}

%Observe that when we are dropping some directions, the sizes in these directions will be set to zero.

\subsection{Construction of the h-sets}

%  Quite often we will write
%$z_-=(x_-,y_-)$ and $z_+=(x_+,y_+)$ as these blocks of coordinates are natural objects to consider when dealing with $\T_j$ or $J_{j \to j+1}$.
%Especially this makes sense if both variables in $z_\pm$ are of the same diameter and are both either entry or exit variables.  In such
%situation we will use $\gamma^j_?(z_?)$ or $r^j_{?}(z_?)$.

\subsubsection{ $N^{j}_{\inn}$ and $\widetilde{N}^{j}_{\inn}$}
\label{subsubsec:sizeNin}

\begin{itemize}
\item[$N^{j}_{\inn}$] is  centered at $A_j$ (see (\ref{eq:Aj})).
  \begin{itemize}
   \item  Its entry variables are $c_p$, $z_-=(x_-,y_-)$,
     \item its exit variables are $z_+=(x_+,y_+)$, $c_f$.
  \end{itemize}
\item[$\widetilde{N}^{j}_{\inn}$] is obtained from   $N^{j}_{\inn}$ by dropping $y_+$ and setting it to $0$, i.e., $\widetilde{N}^j_{\inn}=R_{<y_+>,0}N^j_{\inn}$.
    \begin{itemize}
     \item  Its entry variables are $c_p$, $z_-=(x_-,y_-)$, $y_+$,
     \item  its exit variables are $x_+$, $c_f$.
  \end{itemize}
\end{itemize}

The parameters of $N^{j}_{\inn}$ and $\widetilde{N}^{j}_{\inn}$ only differ in the $y_+$-direction.
In $N^{j}_{\inn}$, $y_+$ is an exit variable of micro size
\[
 |y_+| \leq r^j_{\inn}(y_+) e^{-T}, \quad r^j_{\inn}(y_+)=\frac{3}{2}\sigma, \quad \mbox{(micro)}, % \label{eq:rjy+in}
\]
whereas in $\widetilde{N}^j_{\inn}$ it is an entry variable:
\[
 y_+=0,    \quad  \gamma^j_{\inn}(y_+)=0   \quad \mbox{(nano)}
\]

For the remaining  entry  directions we set
\begin{align*}
 |c_p| &\leq \gamma^j_{\inn}(c_p) e^{-T}, \hspace{4em}\mbox{(micro)} \\
  y_-&\in\sigma + \alpha \gamma^j_{\inn}(y_-)e^{-2T}, \hspace{1em}  \mbox{(macro variable, nano diameter)}, \\
  |x_-|&\leq \gamma^j_{\inn}(x_-)e^{-2T},  \hspace{3em} \mbox{(nano)}. %  \qquad  \gamma^j_{\inn}(x_-)\approx 0.
\end{align*}
%It does not matter whether the diameter of $y_-$ is nano or micro at this stage, however it happens that the nano-size is fine.

For the  remaining  exit directions we set:
\begin{align}
  |x_+| &\leq r^j_{\inn}(x_+) e^{-T},  r^j_{\inn}(x_+)=\frac{3}{2} \sigma,\  \mbox{(micro)} \label{eq:rjinx+}\\
  |c_f| &\leq r^j_{\inn}(c_f)e^{-T}. \hspace{2em} \hphantom{r^j_{\inn}(x_+)=\frac{3}{2}}\mbox{(micro)}\notag
\end{align}

Since the block $z_+=(x_+,y_+)$ consists of two exit variables and as such enters into the covering relation $\widetilde{N}^j_{\out} \cover{\T_j} N^{ j+1}_{\inn}$ then we set
\begin{equation}
 r^j_{\inn}(z_+)=r^j_{\inn}(x_+)=r^j_{\inn}(y_+)=\frac{3 \sigma}{2}.  \label{eq:rjz+in}
\end{equation}

\subsubsection{$N^j_{\out}$ and $\widetilde{N}^j_{\out}$}
%\label{sssec:sizesNjs}

\begin{itemize}
\item[$N^j_{\out}$] is  centered at $B_j$ (see (\ref{eq:Bj})).
   \begin{itemize}
     \item  Its entry variables are $c_p$, $z_-=(x_-,y_-)$, $y_+$,
     \item  its exit variables are $x_+$, $c_f$.
  \end{itemize}
\item[$\widetilde{N}^j_{\out}$] is obtained from $N^j_{\out}$ by dropping $x_+$ and setting $x_+=\sigma$, i.e.,  $\widetilde{N}^j_{\out}=R_{<x_+>,\sigma}N^j_{\out}$.
   \begin{itemize}
     \item  Its entry variables are $c_p$, $z_-=(x_-,y_-)$, $z_+=(y_+,x_+)$,
     \item  its exit variables are $c_f$.
  \end{itemize}
\end{itemize}

The parameters of $N^j_{\out}$ and $\widetilde{N}^j_{\out}=R_{<x_+>,\sigma}N^j_{\out}$ only differ in $x_+$-direction.
In $N^j_{\out}$, $x_+$ is an exit variable of macro size
\begin{equation}
 x_+ \in \sigma + \alpha r^j_{\out}(x_+),    \ r^j_{\out}(x_+) = \sigma/100.  \qquad \mbox{(macro)}   \label{eq:rjsx+}
\end{equation}
the denominator in $r^j_{\out}(x_+)$ being some large number, while in $\widetilde{N}^j_{\out}$ the direction $x_+$ is dropped and we have
\begin{equation}
 x_+ =\sigma, \quad \gamma^j_{\out}(x_+) = 0  \quad \mbox{(nano)}.  \label{eq:njs-x+-zero}
\end{equation}

For the remaining entry directions:
\begin{eqnarray}
   |z_-|&\leq&\gamma^j_{\out}(z_-)e^{-T},  \quad \mbox{(micro)} \notag \\
  |c_p| &\leq& \gamma^j_{\out}(c_p) e^{-T},  \quad\ \mbox{(micro)} \notag \\
   |y_+| &\leq& \gamma^j_{\out}(y_+) e^{-2T}.   \quad \mbox{(nano)}  \label{eq:nsj-y+-nano}
\end{eqnarray}

For the  remaining  exit directions:
\begin{eqnarray*}
  |c_f| &\leq& r^j_{\out}(c_f)e^{-T}.  \quad \mbox{(micro)}
\end{eqnarray*}

The  sizes not given explicitly above  will be determined when we will examine conditions for covering relations.

\subsection{Covering relations}

\subsubsection{Constants for the maps $\varphi_T$ and $\T_j$}
\label{subsec:const-maps}

In this section we define several constants (of Lipschitz type) related to the maps $\varphi_T$ and $\T_j$.
We set $A=3$, fix $\sigma$ in Theorem~\ref{thm:ful-system} and we still have the freedom to increase $T$.
We set (see (\ref{eq:fencl-cj}) in Theorem~\ref{thm:ful-system})
\[
  L_c(\varphi_T)=\exp(10 G \sigma ).  %\label{eq:Lc}
\]
We assume that $\T_j$ is defined in a ball $B(B_j,r_{\T})$ with $r_{\T}$ with the same $r_{\T}$ for all $j$.
We define the following constants for $\T_j$, the same for all $j$ of our transition chain, over the non-trivial blocks
listed in Section~\ref{subsubsec:cc-charts}:
\begin{eqnarray*}
  L(\T) &>& l(\T) \geq \sup_j \left\|\frac{\partial \pi_{c_{j-1}} \T_{j} }{\partial z_-}(B_{j})\right\|,\\% \label{eq:L(J)1} \\
  L(\T) &>& l(\T) \geq \sup_j\left\|\frac{\partial \pi_{\tilde{z}_{-}} \T_{j}}{\partial z_+}(B_j)\right\|,\\  %\label{eq:L(J)2} \\
  L(\T) &>& l(\T) \geq \sup_j\left\|\left(\frac{\partial \pi_{\tilde{z}_{+}} \T_{j} }{\partial c_{j+2}}(B_{j})\right)^{-1}\right\| ,%\label{eq:L(J)3}
\end{eqnarray*}
and
\begin{eqnarray*}
   L_c(\T) &>& l_c(\T) \geq \sup_j\left\|\frac{\partial \pi_{\tilde{c}_{\ell}} \T_{j} }{\partial c_\ell}(B_{j})\right\|, \quad \ell < j-1, \\% \label{eq:Lc(J)1}
  L_c(\T) &>& l_c(\T) \geq  \sup_j\left\|\left(\frac{\partial \pi_{\tilde{c}_{\ell}} \T_{j} }{\partial c_{\ell}}(B_{j})\right)^{-1}\right\|, \quad \ell>j+2. % \label{eq:Lc(J)2}
\end{eqnarray*}

Summarizing, the constants $L(\T)$, $l(\T)$ are defined over the non-zero blocks involving the saddle directions, either in the domain or in the image, and the constants $L_c(\T)$, $l_c(\T)$ over the non-zero blocks in the past or future modes $c_\ell$.

Now a bound for the second derivative of $\T_{j}$  for all $j$:
\begin{equation}
  D_2(\T) = 2 \sup_j \max_{p \in \bar{B}_j(B_j,r_{\T})} \left\| D^2 \T_{j} (p) \right\|,  \label{eq:D2T}
\end{equation}
where we use the norm $\|(x_\pm,y_{\pm},c_*)\|=\max \left(|x_\pm|,|y_\pm|,|c_\ell|_{\ell \neq j\pm 1,j}\right)$, and finally we set
\begin{equation}
   L_{ct}=L_c(\varphi_T) L_c(\T).  \label{eq:Lct}
\end{equation}

\subsubsection{Sizes for $N^j_{\inn}$}
%\label{subsubsec:sizesNjin}

We define
\begin{equation}
k_0=1, \quad k_{c,0}=0, \label{eq:init-k0-kc0}
\end{equation}
and
\begin{eqnarray}
k_{j+1}&=&2k_j+1, \quad k_{c,j+1}=2k_{c,j} +1. \label{eq:kjkcj}
\end{eqnarray}
%\noindenttext{A: by induction \boldmath$k_j=2^{j+1}-1$. If $k_{c,j+1}=2 k_{c,j} +1$, $k_{c,j}=k_{j-1}$.}

It is easy to see that $k_j=2^{j+1}-1$ and $k_{c,j}=2^{j}-1$. This growth of powers of $T$ is apparently avoided in~\cite{GK}.

%\noindenttext{A: Can this be avoided if we can put \boldmath$O_3(z)$ in $\dot{c}_\ell$ in \eqref{eq:NF}?}

We will set the radii $\gamma^j_{\inn}(\dots)$, $r^j_{\inn}(\dots)$ in the following form
\begin{eqnarray}
  \gamma^j_{\inn}(c_p)&=&T^{k_{cj}},  \label{eq:gjin-cp} \\
   \gamma^j_{\inn}(y_-)&=&T^{k_{j}}, \label{eq:gjin-y-} \\
    \gamma^j_{\inn}(x_-)&=&T^{k_{j}}/2, \label{eq:gjin-x-} \\
    r^j_{\inn}(c_f)&=&   L_{ct}^{-j} r_{\inn}^0(c_f), \label{eq:rjin-cf} \\
r_{\inn}^0(c_f)&=&\left\{
  \begin{array}{ll}
     L_{ct}^{n-1}( L(\T) L_c(\varphi_T))\frac{3}{2}\sigma, & \hbox{if $L_{ct} \geq 1$;} \\
    ( L(\T) L_c(\varphi_T))\frac{3}{2}\sigma, & \hbox{if $L_{ct} < 1$.}
  \end{array}
\right. \label{eq:rjin-cf0} % \\
 %    r_{\inn}^0(c_f)&=& L_{ct}^{n-1}( L(\T) L_c(\varphi_T))\frac{3}{2}\sigma, \quad \mbox{if $L_{ct} \geq 1$}, \label{eq:rjin-cf01}  \\
 %   r_{\inn}^0(c_f) &=& ( L(\T) L_c(\varphi_T))\frac{3}{2}\sigma, \quad \mbox{if $L_{ct} < 1$}. \label{eq:rjin-cf02}
\end{eqnarray}

This together with the values of $r^j_{\inn}(x_+)$, $r^j_{\inn}(y_+)$ and $\gamma^j_{\inn}(y_+)$ already defined in Section~\ref{subsubsec:sizeNin}
 determines all the parameters for $N^j_{\inn}$ and $\widetilde{N}^j_{\inn}$.

Note that
 \[
    r_{\inn}^j(c_f) \geq ( L_c(\T) L_c(\varphi_T)) \frac{3 \sigma}{2}, \quad j=0,\dots,n-1. %\label{eq:rincf-lb}
 \]

%\noindenttext{A: check $r_{\inn}^0(c_f)$ and this last bound}

\subsubsection{Covering relation $\widetilde{N}^j_{\inn} \cover{\varphi_T} N^j_{\out}$}

First, using the radii introduced in~(\ref{eq:gjin-cp}--\ref{eq:rjin-cf0}), it is easy to check that  $N^j_{\inn}$ is contained in the set on which the estimates from Theorem~\ref{thm:ful-system} are applicable
with $k=k_j$ and $k_c=k_{c,j}$.
% provided that  all these numbers are bounded.

\begin{lemma}
%\label{lem:cv-in-out}
There exists a constant $K_1=K_1(K,\sigma)$, such that if
\begin{eqnarray}
  \gamma^j_{\out}(z_-) &=& K_1 T^{k_j},  \label{eq:gjsz-} \\
   \gamma^j_{\out}(y_+)&=&K_1 T. \label{eq:gjsy+} \\
 \gamma^j_{\out}(c_p) &=& L_c(\varphi_T) T^{k_{cj}}, \label{eq:gjscp}  \\
  r^j_{\out}(c_f) &=& \frac{1}{L_c(\varphi_T)} r^j_{\inn}(c_f). \label{eq:rjscg}
\end{eqnarray}
then for $T$ large enough
\begin{equation*}
\widetilde{N}^j_{\inn} \cover{\varphi_T} N^j_{\out}, \quad j=0,\dots,n.
\end{equation*}
\end{lemma}
\textbf{Proof:}
We use Theorem~\ref{thm:ful-system} to estimate $\varphi_T$ on $\widetilde{N}^j_{\inn}$ and Theorem~\ref{thm:cv-prod-struct} to establish the
covering relation.

Using $\theta=1-\rho$ as a parameter in Theorem~\ref{thm:ful-system} we obtain a homotopy $H_\rho$ which connects $\varphi_T$ with the shift by $T$
for the flow induced by the linearization of our ODE. %This is a linear map, which is diagonal in the sense required by Theorem~\ref{thm:cv-prod-struct}  \textbf{check whether this makes sense and write it better}.

The entry conditions (\ref{eq:entry-Yi}) are implied by the following inequalities (the rightmost symbol on each line below is the variable for which we check the entry condition):
\begin{eqnarray*}
  L_c(\varphi_T) \gamma^j_{\inn}(c_p)e^{-T}  &\leq& \gamma^j_{\out}(c_p)e^{-T}, \ \qquad c_p \\
   e^{-T}\left(\sigma + \gamma^j_{\inn}(y_-)e^{-2T} + \sigma K/A\right)   &\leq& \gamma^j_{\out}(z_-)e^{-T}, \qquad y_-\\
  e^{-T}T^{k_j} &\leq& \gamma^j_{\out}(z_-)e^{-T}, \qquad x_-\\
 %  e^{-2T}  K \sigma \left( T + T^{p(y_+)} \right)  &\leq& \gamma^j_{\out}(y_+) e^{-2T}, \quad y_+
  e^{-2T}  K \sigma T/A   &\leq& \gamma^j_{\out}(y_+) e^{-2T}. \ \ \quad y_+
\end{eqnarray*}

Now let us consider the exit directions. First we note that $H_1$ in the output directions has the desired affine diagonal form:
\begin{eqnarray*}
  x_+(H_1(x,y))&=& e^T x_+ , \\
  c_\ell(H_1(x,y))&=& e^{i\nu_\ell T}c_\ell, \quad \ell \neq j,j\pm 1.
\end{eqnarray*}

First, for the direction $x_+$, from~(\ref{eq:fencl-x+}) we have
%\begin{eqnarray*}
%  x_+(H_\rho (x,y)) > r^j_{\inn}(x_+) - \sigma/A, \, \mbox{for} \, (x,y) \in \left(\widetilde{N}^j_{\inn}\right)^-,\ x_+=R_{x_+}=r^j_{\inn}(x_+) e^{-T}, \\
%  x_+(H_\rho(x,y)) < -r^j_{\inn}(x_+) + \sigma/A, \, \mbox{for} \, (x,y) \in \left(\widetilde{N}^j_{\inn}\right)^-,\ x_+=-R_{x_+}=-r^j_{\inn}(x_+) e^{-T}.
%\end{eqnarray*}
\begin{align*}
  x_+(H_\rho (x,y)) &> r^j_{\inn}(x_+) - \sigma/A,\ \mbox{for} \, (x,y) \in \left(\widetilde{N}^j_{\inn}\right)^-,&& x_+=R_{x_+}=r^j_{\inn}(x_+) e^{-T}, \\
  x_+(H_\rho(x,y)) &< -r^j_{\inn}(x_+) + \sigma/A, \, \mbox{for} \, (x,y) \in \left(\widetilde{N}^j_{\inn}\right)^-,&&\ x_+=-R_{x_+}=-r^j_{\inn}(x_+) e^{-T}.
\end{align*}
Hence in order to satisfy (\ref{eq:exit-Xj}) we need that (observe that $\tilde{x}^c$ for direction $x_+$ is equal to $\sigma$)
\begin{equation}
  -r^j_{\inn}(x_+) + \sigma/A  < \sigma \pm r^j_{\out}(x_+) <   r^j_{\inn}(x_+) - \sigma/A.  \label{eq:cvin-x+}
\end{equation}
We see that this condition is satisfied if (compare  (\ref{eq:rjinx+},\ref{eq:rjsx+}))
 $3\sigma/2 > (1+1/A + 1/100) \sigma$.
Notice that (\ref{eq:cvin-x+}) implies condition (\ref{eq:xjc-in-image}).

Now, for $c_f=(c_{j+2},c_{j+3},\dots)$, % For each $\ell$, for $c_\ell$ we have the same bounds and the same arguments apply.
let us fix $\ell \geq j+2$. From (\ref{eq:fencl-cj}) it follows that for $p \in \left(\widetilde{N}^j_{\inn}\right)^-$ such that $c_\ell(p)=r_{\inn}^j(c_f)e^{-T}$
we have
\begin{equation*}
  |c_\ell(H_\rho (p))| > \frac{1}{L_c(\varphi_T)} |c_\ell(p)| = \frac{1}{L_c(\varphi_T)}  r_{\inn}^j(c_f)e^{-T}.
\end{equation*}
Hence if
\[
  \frac{1}{L_c(\varphi_T)} r_{\inn}^j(c_f) \geq r_{\out}^j(c_f),
\]
then condition (\ref{eq:exit-Xj}) is satisfied for the $c_f$ directions.

It is easy to see that (\ref{eq:xjc-in-image}) also holds, because  the center point in $c_\ell$-direction in  both h-sets is located at zero
and is mapped onto itself.

From the above discussion and
in view of (\ref{eq:gjin-cp}--\ref{eq:rjin-cf}), in order to have the covering relation it is enough
that for a large enough constant $K_1=K_1(K,\sigma)$ we define $\gamma^j_{\out}(z_-)$, $\gamma^j_{\out}(y_+)$, $\gamma^j_{\out}(c_p)$ and $ r^j_{\out}(c_f)$
as in (\ref{eq:gjsz-}--\ref{eq:rjscg}).
\qed

%Since in Section~\ref{sssec:sizesNjs} we have already fixed values for $r^j_{\out}(x_+)$ and $\gamma^j_{\out}(x_+)$ we now have all parameters defining $N^j_{\out}$ and
%$\widetilde{N}^j_{\out}$.

\subsubsection{Covering relation $ \widetilde{N}^j_{\out} \cover{\T_{j}} N^{j+1}_{\inn}$}

First note that from our construction of $N^j_{\out}$ it follows that for for $T$ large enough (recall that we have only a finite number of steps, so $k_j$ and $k_{c,j}$ are bounded), $N^j_{\out} \subset B(B_{j},r_{\T})$, which means that this h-set is contained in the domain of $\T_{j }$ and thus the various Lipschitz constants defined for this map in subsection~\ref{subsec:const-maps} can be applied. The desired covering relation is a consequence of the following Lemma.

\begin{lemma}
 For $T$ large enough the following covering relation holds
\begin{equation*}
\widetilde{N}^j_{\out} \cover{\T_{j}} N^{j+1}_{\inn}, \quad j=0,\dots,n-1.
\end{equation*}
\end{lemma}
\textbf{Proof:}
We have
\[
  \T_{j} (B_{j}+p)= A_{j+1} + D\T_{j} (B_{j+1})p + R(B_{j},p),
\]
where $R(B_{j},p)=O(|p|^2)$ is the remainder term.

We will use Theorem~\ref{thm:cv-prod-struct} with the homotopy $H$ given by
\[
  H(\rho,B_{j}+p)= A_{j+1} + D\T_{j} (B_{j})p + (1-\rho) R(B_{j},p).
\]

We see that
\begin{eqnarray*}
  H_0&=&\T_j, \\
  H_1(\rho,B_j+p) &=& A_{j+1} + D \T_j(B_j)p.
\end{eqnarray*}
Since we know from our assumptions that $D \T_j(B_j)$ is a block diagonal matrix and the center of $\widetilde{N}^j_{\out}$ ($=B_j$) is mapped onto the center of $N^{ j+1}_{\inn}$ ($=A_{j+1}$), then conditions (\ref{eq:H1-aff-diag},\ref{eq:xjc-in-image}) are satisfied.

From the Taylor expansion with the second order remainder~\eqref{eq:D2T} we have
\[
\|R(B_j,p)\| \leq D_2(\T)T^{2k_j(s)}e^{-2T},
\]
where $k_j(s)$ is the highest power in $T$ that appears just in  the micro-sized expressions $w(T)e^{-T}$ in $N^j_{\out}$, since the powers of the nano-terms (i.e., in the $x_+$ and $y_+$) do not contribute, and so from (\ref{eq:init-k0-kc0},\ref{eq:kjkcj})
\[
  k_j(s)=\max(k_j,k_{c,j})=k_j.  %\label{eq:k(s)-trans}
\]

%Due to the block diagonal structure of $D\T_j$ it makes sense to think in terms of blocks $z_-=(x_-,y_-)$ and $z_+=(x_+,y_+)$. Observe that in %both sets $\widetilde{N}^j_{\out}$ and $N^j_{\out}$ these blocks are entry directions.

Note that in the $z_+$ direction in $\widetilde{N}^j_{\out}$ only the variable $y_+$ has a non-zero diameter (see~\eqref{eq:njs-x+-zero}), therefore by~\eqref{eq:nsj-y+-nano} the block $z_+$ is of nano-size and by~\eqref{eq:gjsy+}) we can set
\begin{equation}
  \gamma^j_{\out}(z_+)=\max (\gamma^j_{\out}(x_+),\gamma^j_{\out}(y_+))= K_1 T.  \label{eq:gjsz+tran}
\end{equation}
In the rest of the proof the variables with tildes will refer to the coordinates used in $N^{j+1}_{\inn}$, i.e., expressed in the normal form coordinates of the $(j+1)$-th chart.

The conditions to be verified for the entry directions are conditions \eqref{eq:entry-Yi} for the blocks $c_{\ell} \mapsto \tilde{c}_{\ell}$ $\ell \leq j-2$, $z_- \mapsto \tilde{c}_{j-1}$, $z_+ \mapsto \tilde{z}_-$), which are implied, respectively, by the following inequalities
\begin{eqnarray*}
  l_c(\T) \gamma^j_{\out}(c_p)e^{-T} + D_2(\T)T^{2k_j}e^{-2T} &<& \gamma^{j+1}_{\inn}(c_p)e^{-T},\\% \label{eq:cov-J-cp}  \\
  l(\T)  \gamma^j_{\out}(z_-) e^{-T} + D_2(\T)T^{2k_j}e^{-2T} &<& \gamma^{j+1}_{\inn}(c_p)e^{-T}, \\% \label{eq:cov-J-xy-} \\
  l(\T) \gamma^j_{\out}(z_+)e^{-2T} + D_2(\T)T^{2k_j}e^{-2T} &<&  \gamma^{j+1}_{\inn}(z_-)e^{-2T}. % \label{eq:cov-J-xy+}
\end{eqnarray*}

For the exit directions $c_{j+1} \mapsto \tilde{z}_+$ and $c_\ell \mapsto \tilde{c}_\ell$ for $\ell \geq j+2$, they are implied by
\begin{eqnarray*}
  \frac{1}{l(\T)} r_{\out}^{j}(c_f)e^{-T} &>&D_2(\T)T^{2k_j}e^{-2T} +r^{j+1}_{\inn}(z_+)e^{-T}, \\
    \frac{1}{l_c(\T)} r_{\out}^{j}(c_f)e^{-T} &>&D_2(\T)T^{2k_j}e^{-2T} + r^{j+1}_{\inn}(c_f)e^{-T}.
\end{eqnarray*}
Let us now see that the above conditions are obtained for $T$ large enough from the following ones, where we drop the terms $D_2(\T)T^{2k_j}e^{-2T}$ in equations
for micro-sized variables, replace $l(\T)$, $l_c(\T)$ by $L(\T)$, $L_c(\T)$ and change strong inequalities to weak ones:
\begin{eqnarray}
  L_c(\T) \gamma^j_{\out}(c_p)  &\leq& \gamma^{j+1}_{\inn}(c_p),   \label{eq:cJ-cp1}  \\
  L(\T)  \gamma^j_{\out}(z_-)   &\leq& \gamma^{j+1}_{\inn}(c_p), \label{eq:cJ-cp2} \\
  L(\T) \gamma^j_{\out}(z_+)  + D_2(\T)T^{2k_j} &\leq& \gamma^{j+1}_{\inn}(x_-),  \label{eq:cJ-x-}  \\
   L(\T) \gamma^j_{\out}(z_+)  + D_2(\T)T^{2k_j} &\leq& \gamma^{j+1}_{\inn}(y_-), \label{eq:cJ-y-} \\
   r_{\out}^{j}(c_f) &\geq&  L(\T)r^{j+1}_{\inn}(z_+), \label{eq:cJ-z+} \\
   r_{\out}^{j}(c_f) &\geq& L_c(\T) r^{j+1}_{\inn}(c_f). \label{eq:cJ-cf}
\end{eqnarray}

Let us examine first the inequalities involving $\gamma^{j+1}_{\inn}(c_p)$, i.e., (\ref{eq:cJ-cp1},\ref{eq:cJ-cp2}).
We need  the following inequalities to be
satisfied (see   (\ref{eq:gjsz-},\ref{eq:gjscp}))
\begin{eqnarray*}
  L_c(\T) L_c(\varphi_T) T^{k_{c,j}} &\leq&  \gamma^{j+1}_{\inn}(c_p), \\
  L(\T) K_1 T^{k_j} &\leq&  \gamma^{j+1}_{\inn}(c_p),
\end{eqnarray*}
which follow directly from (\ref{eq:kjkcj},\ref{eq:gjin-cp}) for $T$ large enough.

Next we continue with (\ref{eq:cJ-x-},\ref{eq:cJ-y-}). From (\ref{eq:gjsz+tran}) it follows that we need the inequalities
\begin{eqnarray*}
      K_2 T^{2k_j} &<& \gamma^{j+1}_{\inn}(x_-),\\
       K_2 T^{2k_j} &<& \gamma^{j+1}_{\inn}(y_-),
\end{eqnarray*}
to be satisfied for a constant $K_2$. Using (\ref{eq:kjkcj},\ref{eq:gjin-x-},\ref{eq:gjin-y-}) it follows that the above inequalities hold for $T$ large enough.

Now we focus on the inequality (\ref{eq:cJ-cf}).
After plugging in (\ref{eq:rjscg}) and (\ref{eq:Lct}) we get
\begin{eqnarray*}
     r_{\inn}^j(c_f) \geq L_c(\varphi_T) L_c(\T) r^{j+1}_{\inn}(c_f) = L_{ct}r^{j+1}_ {\inn}(c_f),
\end{eqnarray*}
which holds due to (\ref{eq:rjin-cf}).

We are left with (\ref{eq:cJ-z+}).
From (\ref{eq:rjscg},\ref{eq:rjin-cf}) we know that
\begin{eqnarray*}
  r_{\out}^j(c_f) = \frac{1}{ L_c(\varphi_T)} r_{\inn}^j(c_f)=\frac{L_{ct}^{-j}}{ L_c(\varphi_T)} r_{\inn}^0(c_f).
\end{eqnarray*}
Now, taking into account equation (\ref{eq:rjz+in}), condition (\ref{eq:cJ-z+}) becomes
\begin{eqnarray*}
   \frac{L_{ct}^{-j}}{ L_c(\varphi_T)} r_{\inn}^0(c_f) \geq L(\T) \frac{3}{2}\sigma,
\end{eqnarray*}
so that
\begin{equation}
   r_{\inn}^0(c_f) \geq L_{ct}^{j} ( L(\T) L_c(\varphi_T)) \frac{3}{2}\sigma. \label{eq:cJ-z+-mod}
\end{equation}
Notice that we require \eqref{eq:cJ-z+-mod} to be satisfied for all $j=0,\dots,n-1$.

Now depending whether $L_{ct} \lessgtr 1$ we need to make a different choice for $r_{\inn}^0(c_f)$.
If $L_{ct} \geq 1$, we take
\[
   r_{\inn}^0(c_f)= L_{ct}^{n-1}( L(\T) L_c(\varphi_T))\frac{3}{2}\sigma,
\]
whereas if $L_{ct} < 1$, then we take
\[
r_{\inn}^0(c_f) = ( L(\T) L_c(\varphi_T))\frac{3}{2}\sigma,
\]
which agrees with (\ref{eq:rjin-cf0}).
\qed

\subsection{Conclusion of the proof of Theorems~\ref{thm:toymodel} and \ref{thm:toymodel-gen}}.

We apply Theorem~\ref{thm:top-shadowing} to chain of covering relations (\ref{eq:tms-nin-nout},\ref{eq:tms-tran-hetero-cv})  to obtain the
desired orbit shadowing our non-transversal heteroclinic chain.

%% file: appendix.tex
\appendix

\input auxLemmas.tex

\input hyp-enclo.tex

\input res-terms-estm.tex

%% file: auxLemmas.tex
\section{Some formulas and Lemmas}

Let $\alpha=[-1,1]$.

%\begin{lemma}
%\label{lem:int-eat}
%\[
%  \int_0^t t e^{at} dt = \frac{e^{at}}{a}\left(t - \frac{1}{a}\right)+\frac{1}{a^2}=\frac{1-e^{at}}{a^2} + \frac{t e^{at}}{a}. % \label{eq:int-te}
%\]
%If $a<0$, then
%\[
%  \int_0^t t e^{at} dt < \frac{1}{a^2}, %\label{eq:ite-a-neg}
%\]
%if $a>0$, then
%\[
%  \int_0^t t e^{at} dt < \frac{t e^{at}}{a}, %\label{eq:ite-a-pos}
%\]
%\end{lemma}

%\begin{lemma}
%If $a>0$ and $n>0$, then
%\begin{equation}
%  \int_0^t t^n e^{at} dt = O\left(\frac{t^n e^{at}}{a}\right), %\label{eq:itne-a-pos}
%\end{equation}
%\end{lemma}

\begin{lemma}
\label{lem:ek>0}
If $\lambda>0$ and $T>0$, $\displaystyle e^{-\lambda T}\frac{e^{\lambda t}-1}{\lambda} \leq t$ for $t \in [0,T]$.
 \end{lemma}
\textbf{Proof:}
For $\displaystyle f(t)=t -  e^{-\lambda T}\frac{e^{\lambda t}-1}{\lambda}$, $f(0)=0$ and $\displaystyle f'(t)=1-e^{\lambda(t-T)}>0$ imply $f(t)\geq 0$, for $t \in [0,T]$.
\qed

%\begin{lemma}
%\label{lem:ek>0}
%Assume that $\lambda>0$ and $T>0$, then
%\[
%  e^{-\lambda T}\frac{e^{\lambda t}-1}{\lambda} \leq t, \quad t \in [0,T].
%\]
%\end{lemma}
%\textbf{Proof:}
%Let us denote $f(t)=t -  e^{-\lambda T}\frac{e^{\lambda t}-1}{\lambda}$. Then $f(0)=0$ and $f'(t)=1-e^{\lambda(t-T)}>0$ for $t \in [0,T)$. Hence $f(t)\geq 0$ for $t \in [0,T]$.
%\qed

\begin{lemma}
\label{lem:ek<0}
Assume $\lambda>0$ and $t>0$. Then $\displaystyle \frac{1-e^{-\lambda t}}{\lambda} < \frac{1}{\lambda}$ and $\displaystyle \frac{1-e^{-\lambda t}}{\lambda} < t$.
\end{lemma}
\textbf{Proof:} The first inequality is clear and the second follows from $\displaystyle \left(\frac{1-e^{-\lambda t}}{\lambda}\right)'=e^{-\lambda t} \in (0,1)$ for $t>0$.
\qed

%\begin{lemma}
%\label{lem:ek<0}
%Assume that $\lambda>0$, then for $ t \in (0,\infty)$  holds
%\begin{eqnarray*}
%  \frac{1-e^{-\lambda t}}{\lambda} \leq t, \\
%   \frac{1-e^{-\lambda t}}{\lambda} < \frac{1}{\lambda}.
%\end{eqnarray*}
%\end{lemma}
%\textbf{Proof:} Observe that $\left(\frac{1-e^{-\lambda t}}{\lambda}\right)'=e^{-\lambda t} \in (0,1)$ for $t \in (0,\infty)$. From this we immediately get our assertions.
%\qed

\begin{lemma}
\label{lem:estm-stable-unstable}
Let $\lambda\neq \lambda_i \in \mathbb{R}$ and $D, D_i\geq 0$ for $i$ in a finite set $I$. Assume that
 \[
    \dot y \in \lambda y + \alpha D e^{\lambda t} + \alpha \sum_{i \in I} D_i e^{\lambda_i t}.
\]
Then the following inclusion is satisfied for $t \in [0,T]$:
\[
  y(t) \in  e^{\lambda t} y(0) + \alpha e^{\lambda t} \left( D t +  \sum_{i \in I}D_i \frac{e^{(\lambda_i-\lambda)t} -1}{\lambda_i-\lambda}\right)
       \subset  e^{\lambda t} y(0) + E_{y}(t),
%       \alpha t e^{\lambda t}\left( D  +  \sum_{\substack{i \in I \\ \lambda_i < \lambda}}D_i  +
%           \sum_{\substack{i \in I \\  \lambda_i > \lambda}}D_i e^{(\lambda_i - \lambda)T} \right),
\]
where
\begin{equation}
\label{Enclosurey}
E_{y}(t)=\alpha t e^{\lambda t}\left( D  +  \sum_{\substack{i \in I \\ \lambda_i < \lambda}}D_i  +
           \sum_{\substack{i \in I \\  \lambda_i > \lambda}}D_i e^{(\lambda_i - \lambda)T} \right)
\end{equation}  
is the \emph{enclosure} to the linear solution $e^{\lambda t} y(0)$ originated by the term $\alpha D e^{\lambda t} + \alpha \sum_{i \in I} D_i e^{\lambda_i t}.$          
\end{lemma}
\textbf{Proof:}
For $\lambda_i < \lambda$ and $t \in [0,\infty)$, by Lemma~\ref{lem:ek<0} we have
\begin{equation*}
  \frac{e^{(\lambda_i-\lambda)t} -1}{\lambda_i-\lambda} \leq t,
\end{equation*}
whereas for $\lambda_i>\lambda$ $t \in [0,T]$, by Lemma~\ref{lem:ek>0} we get
\begin{equation*}
  \frac{e^{(\lambda_i-\lambda)t} -1}{\lambda_i-\lambda} \leq e^{(\lambda_i - \lambda)T}t.
\end{equation*}
\qed

%% file: hyp-enclo.tex
\section{Proof of Theorem~\ref{thm:enclo-hyp-dir}}
\label{sec:proof-estm}

For the proof it will not be necessary to assume that the monomials $z^m$ of the system of Theorem~\ref{thm:enclo-hyp-dir} satisfy $m\in M_v:= M_{1,v}\cup M_{ 2, v }$. Nevertheless, this restriction is justified by Lemma~\ref{lem:PolyNormForm}, and allows to shorten the proof, thus avoiding considering about 30 different types of monomials.
All relevant bounds in the proof increase with $\rho$, so it suffices to consider $\rho=1$ .

\subsection{The set $\mathcal{W}$}

Let us introduce the set $\mathcal{W}$ which will be very convenient to establish the enclosures of Theorem~\ref{thm:enclo-hyp-dir}. We will denote by $\mathcal{W}$ the set of functions $(x_\pm(t),y_\pm(t))$ satisfying  (\ref{eq:it-bnd-x-}--\ref{eq:it-bnd-y+}) for $t \in [0,T]$:
\begin{eqnarray}
  x_-(t) &\in&  \mathcal{W}_{x_-}(t):= \alpha e^{-2T}e^{t} T^k,  \label{eq:it-bnd-x-} \\
 y_-(t) &\in&  \mathcal{W}_{y_-}(t):=\alpha e^{-t}\cdot 2 \cdot \sigma , \label{eq:it-bnd-y-} \\
  x_+(t) &\in&  \mathcal{W}_{x_+}(t):= \alpha e^{-T}e^t 2\sigma. \label{eq:it-bnd-x+}\\
    y_+(t) &\in&  \mathcal{W}_{y_+}(t):= \alpha e^{-t} e^{-T}  K \sigma ( t + e^{-q(y_+) T}),   \label{eq:it-bnd-y+}
\end{eqnarray}
where $q(y_+)\geq 1$.

%\textbf{PZ: this $e^{-q(y_+)T}$ in $\mathcal{W}_{y_+}$ is included to make sure that $y_+(0)=0 \in \inter \mathcal{W}_{y_+}(0)$}

%We will use $\mathcal{W}(t)$ to denote the set $\{w(t) \ | \  w \in \mathcal{W} \}$.
%It is easy to see that  $\mathcal{W}(t) \subset Z$ for $t \in [0,T]$.
%Moreover, note that $\mathcal{W}(0)$ contains the initial conditions given by (\ref{eq:sys-ic}) with some margin. This is done  to leave space for various estimates.
%We will use the  continuation method~\cite{Tao} to prove that the solutions of~(\ref{eq:sys-res-0-full-xy-1}--\ref{eq:sys-res-0-full-xy-4})
%satisfying  the initial conditions~(\ref{eq:sys-ic}) are in $\mathcal{W}$ for $t \in [0,T]$.
%
%Let us define
%\[
%  d(y_-)=0, \quad d(x_-)=2, \quad d(y_+)=1, \quad d(x_-)=1. %\label{eq:def-d()}
%\]
%$d(v)$ appears in $\mathcal{W}_{v}(t)$ in the factor $e^{-d(v)T}$.

We will use the notation $\mathcal{W}(t)$ to denote the set $\{w(t):w \in \mathcal{W} \}$.
It is easy to see that $\mathcal{W}(t) \subset Z$ for $t \in [0,T]$.
Moreover, note that $\mathcal{W}(0)$ contains the initial conditions given by~\eqref{eq:sys-ic} with some extra margin. This is done to leave room for several estimates.
We will use the continuation method (see, for instance,~\cite{Tao}) to show that the solutions of~(\ref{eq:sys-res-0-full-xy-1}--\ref{eq:sys-res-0-full-xy-4})
that satisfy the initial conditions~\eqref{eq:sys-ic} are contained in $\mathcal{W}$ for $t \in [0,T]$.

We will denote by $d(v)$ the exponent that appears in $\mathcal{W}_{v}(t)$ in the factor $e^{-d(v)T}$:
\[
   d(x_-)=2, \quad d(y_-)=0, \quad d(x_+)=1, \quad d(y_+)=1.
%\label{eq:def-d()}
\]

\subsection{The contribution of each monomial to the solutions of (\ref{eq:sys-res-0-full-xy-1}--\ref{eq:sys-res-0-full-xy-4})}

\begin{definition}
Given a monomial $z^m$ and a variable $v \in \{x_\pm,y_\pm\}$, an \emph{enclosure $E_{v,m}(t)$} is an interval in $\mathbb{R}$ which contains the enclosure~\eqref{Enclosurey}
obtained from Lemma~\ref{lem:estm-stable-unstable} originated by the $v$-component of $w^m(t)$ for a function $w\in\mathcal{W}$.
\end{definition}
Quite often we will also write $E_{v,z^m}$ instead of $E_{v,m}$.

\subsubsection{Some examples of computation of $E_{v,m}$}

To illustrate the definition of an enclosure $E_{v,m}$, let us first compute  $E_{v,m}$ for  the variable $x_-$ and the monomial $y_- x_+^2$.
From estimates~\eqref{eq:it-bnd-y-} and~\eqref{eq:it-bnd-x+}
for the functions in  $\mathcal{W}$ we have
\[
(y_- x_+^2) (t) \in \alpha e^{t} e^{-2T} 8 \sigma^3,
\]
and from enclosure~\eqref{Enclosurey} of Lemma~\ref{lem:estm-stable-unstable} we get
\begin{equation}
\label{eq:Exminus}
  E_{x_-, y_- x_+^2} (t)=\alpha t e^{t} e^{-2T} 8 \sigma^3.
\end{equation}

Analogously, for the variable $y_+$ and the monomial $y_-^2 x_+$, again from~\eqref{eq:it-bnd-y-} and~\eqref{eq:it-bnd-x+}
\[
(y_-^2 x_+) (t) \in \alpha e^{-t} e^{-T} 8 \sigma^3,
\]
and from enclosure~\eqref{Enclosurey} of Lemma~\ref{lem:estm-stable-unstable} we get
\begin{equation}
\label{eq:Eyplus}
  E_{y_+, y_-^2 x_+} (t)=\alpha t e^{-t} e^{-T} 8 \sigma^3.
\end{equation}

One more, for  the variable $y_-$ and the monomial $x_- x_+$,
from~\eqref{eq:it-bnd-x-} and~\eqref{eq:it-bnd-x+}
we have
\[
(x_- x_+) (t) \in \alpha e^{2 t} e^{-3T} 2 \sigma   \cdot T^k,
\]
and from enclosure~\eqref{Enclosurey} of Lemma~\ref{lem:estm-stable-unstable} we get
\[
  E_{y_-,x_-x_+} (t)=e^{-t}\alpha \sigma \cdot 2  \cdot T^k \cdot e^{-3T}\frac{e^{3t}-1}{3}  \subset
  e^{-t}\alpha  \sigma T^k  e^{3(t-T)}.
\]

\subsection{The notion of suitable monomial for the set $\mathcal{W}$}

We  define the notion of suitable monomial for the variable $v \in \{x_\pm,y_\pm\}$ as follows
\begin{definition}
%\label{def:good-mono}
 $z^m$   is  \emph{suitable monomial} for the variable $v \in \{x_\pm,y_\pm\}$ (with respecto to the set $\mathcal{W}$) if  for any $C>0$
there exist $\sigma_0$, such that for any $0 <\sigma \leq \sigma_0$  there exists $T_0=T_0(\sigma)$  such that the following holds for $T \geq T_0$ and $t \in [0,T]$
 \begin{equation}\label{eq:gm-v}
       E_{v,m}(t) \subset  \frac{1}{C} \mathcal{W}_{v}(t).
      \end{equation}

 %\begin{itemize}
 %\item for $v=y_-$
 %     \begin{equation}\label{eq:gm-y-}
 %      E_{y_-,J}(t) \subset  \frac{1}{C} \mathcal{W}_{y_-}(t)   % \alpha \frac{e^{-t}}{C}\left(\sigma/2 + \sigma  K T^{p(y_-)} e^{w(y_-)(t-T)} \right).
 %     \end{equation}

% \item for $v=x_-$

 %     \begin{eqnarray}\label{eq:gm-x-}
 %         E_{x_-,J}(t) &\subset&  \frac{1}{C} \mathcal{W}_{x_-}(t)   % \alpha \frac{1}{C} e^{-2T}e^t \left( T^k/2 + K \sigma T^{p(x_-)} \right).
  %    \end{eqnarray}

% \item for $v=y_+$
%      \begin{equation}\label{eq:gm-y+}
%       E_{y_+,J}(t) \subset  \frac{1}{C} \mathcal{W}_{y_+}(t) % \alpha \frac{1}{C} e^{-t}e^{-T} K \sigma \left( t + T^{p(y_+)} e^{w(y_+)(t-T)}  \right).
%      \end{equation}

%\item for $v=x_+$
%      \begin{equation}\label{eq:gm-x+}
%          E_{x_+,J}(t) \subset  \frac{1}{C} \mathcal{W}_{x_+}(t)  % \alpha e^t e^{-T}\frac{\sigma}{2C}.
%      \end{equation}

% \end{itemize}
\end{definition}

\begin{definition}
\label{def:good-set-monomials}
We say that system~(\ref{eq:sys-res-0-full-xy-1}--\ref{eq:sys-res-0-full-xy-4}) (or more precisely, the vector field $f=(f_{x_-},f_{y_-},f_{x_+},f_{y_+})$) is \emph{suitable} (for $\mathcal{W}$) if for each $m \in M_v :=M_{1,v} \cup M_{2,v}$ the monomial $z^m$ is suitable for the variable $v$ (with respecto to $\mathcal{W}$).

\end{definition}

\subsection{Theorem on estimates for a suitable system (\ref{eq:sys-res-0-full-xy-1}--\ref{eq:sys-res-0-full-xy-4})}

\begin{theorem}
\label{thm:hype-ode-estm}
Assume that system~(\ref{eq:sys-res-0-full-xy-1}--\ref{eq:sys-res-0-full-xy-4}) is suitable in the sense of Definition \ref{def:good-set-monomials}.
Then there exist $\sigma_0>0$ such that for any $0<\sigma \leq \sigma_0$ there exists $T_0=T_0(\sigma)>0$, such that
 for  $T \geq T_0$ the solutions of (\ref{eq:sys-res-0-full-xy-1}--\ref{eq:sys-res-0-full-xy-4}) with initial conditions (\ref{eq:it-sys-ic})
 belong to  $\mathcal{W}(t)$  for  $t \in [0,T]$ and satisfy~\eqref{eq:fen}.

\end{theorem}

\noindent
\textbf{Proof:}
Denote $K=\sup_{(z,c) \in Z} |g_{v,m}(z,c)|$ and $\#M_v$ the number of elements in $M_v$.

Since any monomial $z^m$ with $m$ in $M_{x_\pm}$ or $M_{y_\pm}$ is suitable, there exists $\sigma_0 >0$  such that
 for $0<\sigma\leq \sigma_0$  there exists $T_0=T_0(\sigma)$ so that, for  $T\geq T_0$,  conditions (\ref{eq:gm-v}) are satisfied for $v \in \{y_\pm,x_\pm\}$ with
 \[
 C \geq A\cdot K \cdot \max\{\#M_{x_-},\#M_{y_-},\#M_{x_+},\#M_{y_+}\}.
 \]

For the proof let us fix initial conditions~(\ref{eq:it-sys-ic})
for $x_\pm,y_{\pm}$ and take $T \geq T_0$, $\sigma \leq \sigma_0$ and $c_*(t) \in C^*$ for $t \in [0,T]$.
By the continuation-type argument \cite{Tao}  one can easily show that the solution of our system is in $\mathcal{W}$. For instance, let us illustrate it for $x_-$-component.

Assume that $x_\pm(t),y_{\pm}(t) $ are in $\mathcal{W}(t)$  for $t \in [0,T']$, $T' \leq T$, then
from condition~(\ref{eq:gm-v}) with $v=x_-$   it follows that
\begin{eqnarray*}
  x_-(t) &\in&
   e^{t} x_-(0) + \sum_{m \in M_{x_-}}G_{x_-,m}  E_{x_-,m}(t)  \subset e^{t} e^{-2T} a_0 + \#M_{x_-}\cdot K\cdot  \frac{1}{C}  \cdot \mathcal{W}_{x_-}(t)   \\
   &\subset&  e^{t} e^{-2T} a_0 + \alpha e^{t} e^{-2T} \frac{T^{k}}{A}= \alpha e^t e^{-2T} \frac{T^k}{2} +  \frac{1}{A} \mathcal{W}_{x_-}(t)=\left(\frac{1}{2}+\frac{1}{A}\right) \mathcal{W}_{x_-}(t)
   \subset \inter\mathcal{W}_{x_-}(t).
\end{eqnarray*}
Analogously for the other variables. Hence, by the continuation argument, the solution is in $\mathcal{W}$.

Notice that the first inclusion in the last line above establishes the first enclosure of~\eqref{eq:fen}. The same happens with the other variables, and thus the enclosure~\eqref{eq:fen}  is established.
\qed

%% file: res-terms-estm.tex
\section{Part 2 of the proof of Theorem~\ref{thm:enclo-hyp-dir}---the estimates}
\label{sec:good-terms}
Our goal is to check that all the monomials $z^m$ of system~(\ref{eq:sys-res-0-full-xy-1}-\ref{eq:sys-res-0-full-xy-4}) are suitable monomials to apply
Theorem~\ref{thm:hype-ode-estm} and thus prove Theorem~\ref{thm:enclo-hyp-dir}. Thanks to the Normal Form Lemma~\ref{lem:PolyNormForm} all of them satisfy that $m \in M_{1,v} \cup M_{2,v}$ for $v \in \{x_\pm,y_\pm\}$ so we simply need to check all the various types of monomials.

%\textbf{Maybe remove these comments?}
%In our search of the vector fields so that estimates from Theorem~\ref{thm:enclo-hyp-dir} are satisfied we restrict ourselves to
%systems satisfying assumptions \textbf{A1} and \textbf{A2} (see Section~\ref{sec:encl-full-system}. According with Def.~\ref{def:illegal} terms  incompatible with those %assumptions are  \emph{illegal}. Additionally we assume that
%only resonant polynomials are present in (\ref{eq:f-full-enclo}).

\subsection{Very suitable, potentially suitable and unsuitable monomials}
Recall that in Section~\ref{sub:polynomialNFc} we already introduced $\lambda_{x_\pm}=1$ and $\lambda_{y_\pm}=-1$.
In addition, for each monomial $z^m$  and $v \in \{x_\pm,y_\pm\}$ we introduce the following constants
\begin{eqnarray}
\lambda_m&=&m_{x_-}-m_{y_-}+ m_{x_+}-m_{y_+}, \label{eq:lamb-m} \\
 \kappa_{m}&=& m_c+2 m_{x_-} + m_{x_+}+ m_{y_+}, \label{eq:kappa-m} \\
 \theta_m &=& m_{x_+}+m_{y_-}+m_{y_+}, \notag \\ %\label{eq:theta-m}
 s_m &=& k m_{x_-}+ m_{y_+} + k_c m_c. \notag
\end{eqnarray}

We now can write a general estimate for a monomial in the set $\mathcal{W}$ for $t \in [0,T]$ (we assume that $T\geq 1$ in order to have $e^{-q(y_+)T} \leq T$):
\begin{eqnarray}
z^m(t) &\in& \alpha e^{\lambda_m t} e^{-T \kappa_m} \sigma^{\theta_m} 2^{m_{y_-}+m_{x_+}}  K^{m_{y_+}} B^{m_c}
       T^{k m_{x_-}+ k_c m_c}  \left(  e^{-q(y_+)T} + t  \right)^{m_{y_+}}  \nonumber  \\
  &\subset& \alpha e^{\lambda_m t} e^{-T \kappa_m} (2\sigma)^{\theta_m}   K^{m_{y_+}} B^{m_c}
       T^{s_m}. \label{eq:zm-estm}
\end{eqnarray}

From \eqref{Enclosurey} of Lemma~\ref{lem:estm-stable-unstable} and (\ref{eq:zm-estm}),  for $t \in [0,T]$
it follows that
\begin{equation}
\label{eq:e-est}
 E_{v,m}(t)\subset\alpha e^{\lambda_v t} e^{-a(v,m) T} (2\sigma)^{\theta_m}   K^{m_{y_+}} B^{m_c}
       T^{s_m}   t ,
\end{equation}
where
\begin{equation}
\label{eq:av}
a(v,m)=\begin{cases} \kappa_m -\lambda_m +\lambda_v,& \lambda_m > \lambda_v, \\
                    \kappa_m, & \lambda_m \leq \lambda_v.
\end{cases}
\end{equation}

We now introduce what we mean by \emph{very suitable}, \emph{potentialy suitable} and \emph{unsuitable} monomials.
For a monomial to be \emph{very suitable} we want the constant $a(v,m)$ defined in~\eqref{eq:av} which appears in the factor $e^{-a(v)T}$ of $E_{v,m}(t)$ in \eqref{eq:e-est}
to be greater than $d(v)$, if we have the equality then the monomial $z^m$ will be called \emph{potentially suitable} and will require further analysis to check if it gives rise to a suitable monomial, and otherwise the monomial will be called \emph{unsuitable}.
Their explicit determination clearly depends on the values of $\lambda_m$ and $\lambda_v$.

%\subsubsection{Case $\lambda_m > \lambda_v$}
%
%In view of (\ref{eq:e-est->}) we introduce the following definition
%\begin{eqnarray*}
%  \kappa_m - \lambda_m+\lambda_v &>& d(v), \quad \mbox{the term is \correct },\\  %\label{eq:gbb->g}\\
%    \kappa_m - \lambda_m+\lambda_v &=& d(v), \quad \mbox{the term is borderline},\\%\label{eq:gbb->bo} \\
%      \kappa_m - \lambda_m+\lambda_v &<& d(v), \quad \mbox{the term is bad}. %\label{eq:gbb->bad}
%\end{eqnarray*}
%
%Due (\ref{eq:e-est->}) we see that if term is \correct for variable $v$, then there exists $T_0$ such that for all $T\geq T_0$  condition (\ref{eq:gm-v}) holds. $T_0$ depends on $K$,  monomial $m$ and $C$. Therefore this monomial is good.
%
%Moreover, if the term is bad, then no such $T_0$ exists.
%
%
%\subsubsection{$\lambda_m \leq \lambda_v$}
%
%
%
%In view of (\ref{eq:e-estm-<=}) propose the following definition
%\begin{eqnarray}
%  \kappa_m  &>& d(v), \quad \mbox{the term is \correct },  \label{eq:gbb<g}\\
%    \kappa_m  &=& d(v), \quad \mbox{the term is borderline}, \label{eq:gbb<bo} \\
%      \kappa_m  &<& d(v), \quad \mbox{the term is bad}. \label{eq:gbb<bad}
%\end{eqnarray}
%Due to (\ref{eq:e-estm-<=}) we see that  if a term is \correct for variable $v$, then there exists $T_0$ such that for all $T\geq T_0$  condition (\ref{eq:gm-v}) holds. $T_0$ depends on $K$,  monomial $m$ and $C$. Therefore this monomial is good.
%
%
%Moreover, if the term is bad, then no such $T_0$ exists.

\subsection{All monomials $z^m$ with $m\in M_{2,v}$ are very suitable}
In this case $m_c\geq 3$. Notice that $d(v) \leq 2$ and $d(v) - \lambda_v \leq 2$. By \eqref{eq:lamb-m} and \eqref{eq:kappa-m},
\begin{eqnarray*}
 \kappa_m - \lambda_m +\lambda_v&=& m_c + m_{y_-} + m_{x_-} + 2m_{y_+}+\lambda_v \geq m_c+d(v)-2>d(v)\\
  \kappa_m &\geq& m_c> d(v),
\end{eqnarray*}
so $a(v,m)>d(v)$ and $z^m$ is very suitable.

\subsection{All monomials $z^m$ with $m\in M_{1,v}$ are suitable}
In this case $m_c=0$, so $\lambda_m=\lambda_v=m_{x_-}-m_{y_-}+ m_{x_+}-m_{y_+}$ with $m_{x_-}+m_{y_-}+ m_{x_+}+m_{y_+}\geq 2$ and $a(m,v)=\kappa_m=2 m_{x_-} + m_{x_+}+ m_{y_+}$.
Now we will check that $\kappa_m> d(v)$ holds for all monomials $z^m$ with $m\in M_{1,v}$ and $v \in \{x_-,y_-,x_ + ,y_+\}$, except for two monomials for which $\kappa_m= d(v)$ and which will therefore require further study.

\subsubsection{Variable $x_-$}
We have $\lambda_{x_-}=1$ and $d(x_-)=2$. Therefore $m_{x_-}+ m_{x_+}=1+m_{y_-}+m_{y_+}$ and $1+2(m_{y_-}+ m_{y_+})\geq 2$ so that $m_{y_-}+ m_{y_+}\geq 1$. Then
\[
\kappa_m=m_{x_-} + m_{x_+}+m_{x_-}+ m_{y_+}=1+m_{y_-} + m_{y_+}+m_{x_-}+ m_{y_+}\geq 2+m_{x_+}+ m_{y_+}\geq 2=d(x_-),
\]
so that there are no unsuitable monomials. All the monomials are very suitable  and, in particular, suitable, except those satisfying $\kappa_m=2$, which only takes place when $m_{y_-}+ m_{y_+}=1$ and $m_{x_-}=m_{y_+}=0$ and therefore $m_{y_-}=1$ and $m_{x_+}=2$, so
there is only one potentially suitable monomial $y_- x_+^2$, for which we already obtained
in~\eqref{eq:Exminus} its enclosure
\[
  E_{x_-, y_- x_+^2} (t)=\alpha t e^{t} e^{-2T} 8 \sigma^3,
\]
which for small $\sigma$ and $t\in[0,T]$ is included in $\alpha e^t e^{-2T}   T^{k}/C$, because by our assumption   $k \geq 1$.

\subsubsection{Variable $y_-$}

We have $\lambda_{y_-}=-1$ and $d(y_-)=0$. Therefore $m_{y_-}+m_{y_+}=1+m_{x_-}+ m_{x_+}$ and $1+2(m_{x_-}+ m_{x_+})\geq 2$ so that $m_{x_-}+ m_{x_+}>0$ and consequently
$a(m,y_-)=\kappa_m>0=d(y_-)$, so $z^m$ is very suitable and, in particular, suitable.

\subsubsection{Variable $x_+$}
We have $\lambda_{x_+}=1$ and $d(x_+)=1$. Therefore $m_{x_-}+m_{x_+}=1+m_{y_-}+ m_{y_+}$ and $1+2(m_{y_-}+ m_{y_+})\geq 2$ so that $m_{y_-}+ m_{y_+}\geq 1$ and consequently
$m_{x_-}+m_{x_+}\geq 2$ so that $a(m,y_-)=\kappa_m\geq 2 >1=d(x_+)$, so $z^m$ is very suitable and, in particular, suitable.

\subsubsection{Variable $y_+$}
%This case is very similar to the previous one, exchanging $x\longleftrightarrow y$ and $+\longleftrightarrow -$.
We have $\lambda_{y_+}=-1$ and $d(y_+)=1$. Therefore $m_{y_-}+ m_{y_+}=1+m_{x_-}+m_{x_+}$ and $1+2(m_{x_-}+ m_{x_+})\geq 2$ so that $m_{x_-}+ m_{x_+}\geq 1$. Then
\[
\kappa_m=m_{x_-} + m_{x_+}+m_{x_-}+ m_{y_+}\geq 1+m_{x_-}+ m_{y_+}\geq 1=d(y_+),
\]
so that there are no unsuitable monomials. All the monomials are very suitable and, in particular, suitable, except those satisfying $\kappa_m=1$, which only takes place when $m_{x_-}+ m_{x_+}=1$ and $m_{x_-}=m_{y_+}=0$ and therefore $m_{x_+}=1$ and $m_{y_-}=2$, so
there is only one potentially suitable monomial $y_- x_+^2$, for which we already obtained
in~\eqref{eq:Eyplus} its enclosure
\[
  E_{y_+, y_-^2 x_+} (t)=\alpha t e^{-t} e^{-T} 8 \sigma^3,
\]
which for small $\sigma$ and $t\in[0,T]$ is included in $\alpha t e^{-t} e^{-T} /C$.

\subsection{Conclusion}
We have shown that all the monomials $z^m$ with $m \in M_{1,v} \cup M_{2,v}$ for $v \in \{x_\pm,y_\pm\}$ are suitable monomials. Therefore from Theorem~\ref{thm:hype-ode-estm} we obtain assertion~\eqref{eq:fen} of Theorem~\ref{thm:enclo-hyp-dir}.